%% file: main_arxiv.tex
\documentclass[a4paper, 10pt ]{article}
\usepackage[a4paper,left=1.25cm,right=1.25cm,top=2.5cm,bottom=2.5cm]{geometry}
\usepackage{hyperref}
\hypersetup{colorlinks,linkcolor={blue},citecolor={blue},urlcolor={black}}
\usepackage{pgfplots}
\usepackage{pgfplotstable}
\usepackage{mathtools}
\usepackage{esint}
\usepackage{multicol}
\usepackage{comment}
\usepackage{booktabs}
\pgfplotsset{compat=1.5}
\usepackage{amssymb}
\usepackage{url}
\usepackage{bm}
\usepackage{caption}
\usepackage{overpic}
\usepackage{pdfsync}
\usepackage{float}
\usepackage{tabularx}
\usepackage{enumerate}
\usepackage{array}
\usepackage{xspace}
\usepackage{tikz}
\usepackage{tikz-cd}
\usepackage{tikzsymbols}
\usetikzlibrary{calc,trees,positioning,arrows,chains,shapes.geometric,%
decorations.pathreplacing,decorations.pathmorphing,shapes,%
matrix,shapes.symbols, decorations.markings, patterns,fit}

\usepackage{siunitx}
\usepackage{subcaption}
\usepackage{authblk}

\usepackage[draft,inline,marginclue]{fixme}
\usepackage{mathrsfs}
\usepackage{color, colortbl}
\usepackage{multirow}

\usepackage{amsthm}
\usepackage{amsmath}
\usepackage{stmaryrd}
\usepackage[ruled,vlined,linesnumbered]{algorithm2e}
\usepackage{graphicx}
\usepackage{booktabs}
\usepackage{cleveref}
\usepackage{verbatim}
\usepackage{arydshln}
\usepackage{etoolbox}
\usepackage{nomencl}
\usepackage[inkscapeformat=png]{svg}

\FXRegisterAuthor{gs}{adt}{\color{red}GS}
\FXRegisterAuthor{fr}{afr}{\color{blue}FR}
\FXRegisterAuthor{gp}{agp}{\color{magenta}GP}

\theoremstyle{remark}

\DeclareMathOperator*{\argmin}{\arg\!\min}

\newcolumntype{C}[1]{>{\centering\arraybackslash}m{#1}}

\definecolor{Gray}{gray}{0.9}

\makenomenclature

\providetoggle{nomsort}
\settoggle{nomsort}{true} 

\makeatletter
\iftoggle{nomsort}{%
    \let\old@@@nomenclature=\@@@nomenclature        
        \newcounter{@nomcount} \setcounter{@nomcount}{0}%
        \renewcommand\the@nomcount{\two@digits{\value{@nomcount}}}
        \def\@@@nomenclature[#1]#2#3{
          \addtocounter{@nomcount}{1}%
        \def\@tempa{#2}\def\@tempb{#3}%
          \protected@write\@nomenclaturefile{}%
          {\string\nomenclatureentry{\the@nomcount\nom@verb\@tempa @[{\nom@verb\@tempa}]%
          \begingroup\nom@verb\@tempb\protect\nomeqref{\theequation}%
          |nompageref}{\thepage}}%
          \endgroup
          \@esphack}%
      }{}
\makeatother

\begin{document}

\title{Generative Models
for the Deformation of Industrial Shapes with Linear Geometric Constraints: model order and\\ parameter space reductions}

\author[1]{Guglielmo~Padula\footnote{guglielmo.padula@sissa.it}}
\author[1]{Francesco~Romor\footnote{francesco.romor@sissa.it}}
\author[2]{Giovanni~Stabile\footnote{giovanni.stabile@uniurb.it}}
\author[1]{Gianluigi~Rozza\footnote{gianluigi.rozza@sissa.it}}

\affil[1]{Mathematics Area, mathLab, SISSA, via Bonomea 265, I-34136 Trieste,
Italy}
\affil[2]{Department of Pure and Applied Sciences, Informatics and Mathematics Section, University of Urbino Carlo Bo, Piazza della Repubblica, 13, I-61029 Urbino, Italy}

\maketitle

\begin{abstract}
  \input{./sections/abstract.tex}
\end{abstract}

\tableofcontents
\listoffixmes

\section{Introduction}
\label{sec:intro}
\input{./sections/intro.tex}

\section{Constrained Free Form Deformation}
\label{sec:methods}
\input{./sections/methods.tex}

\section{Constrained Generative Models}
\label{sec:constrained generative models}
\input{./sections/cgm.tex}

\section{Surrogate modelling}
\label{sec:mor}
\input{./sections/mor.tex}

\section{Numerical results}
\label{sec:results}
\input{./sections/results.tex}

\section{Discussion}
\label{sec:discussion}
\input{./sections/discussion.tex}

\section{Conclusions and perspectives}
\label{sec:conclusions}
\input{./sections/conclusions.tex}

\section*{Acknowledgements}
This work was partially funded by European Union Funding for Research and Innovation — Horizon 2020 Program —
in the framework of European Research Council Executive Agency: H2020 ERC CoG 2015 AROMA-CFD project 681447
“Advanced Reduced Order Methods with Applications in Computational Fluid Dynamics” P.I. Professor Gianluigi Rozza.
We also acknowledge the PRIN 2017 “Numerical Analysis for Full and Reduced Order Methods for the eﬃcient and accurate
solution of complex systems governed by Partial Diﬀerential Equations” (NA-FROM-PDEs). We thank the Stanford University Computer Graphics Laboratory for having made available the STL file of the Stanford bunny.

\newpage
\appendix
\input{./sections/appendix.tex}

\bibliographystyle{abbrv}
\bibliography{biblio}


\end{document}

%% file: sections/abstract.tex
Real-world applications of computational fluid dynamics often involve the evaluation of quantities of interest for several distinct geometries that define the computational domain or are embedded inside it. For example, design optimization studies require the realization of response surfaces from the parameters that determine the geometrical deformations to relevant outputs to be optimized. In this context, a crucial aspect to be addressed is represented by the limited resources at disposal to computationally generate different geometries or to physically obtain them from direct measurements. This is the case for patient-specific biomedical applications for example. When additional linear geometrical constraints need to be imposed, the computational costs increase substantially. Such constraints include total volume conservation, barycenter location and fixed moments of inertia. We develop a new paradigm that employs generative models from machine learning to efficiently sample new geometries with linear constraints. A consequence of our approach is the reduction of the parameter space from the original geometrical parametrization to a low-dimensional latent space of the generative models. Crucial is the assessment of the quality of the distribution of the constrained geometries obtained with respect to physical and geometrical quantities of interest. Non-intrusive model order reduction is enhanced since smaller parametric spaces are considered. We test our methodology on two academic test cases: a mixed Poisson problem on the 3d Stanford bunny with fixed barycenter deformations and the multiphase turbulent incompressible Navier-Stokes equations for the Duisburg test case with fixed volume deformations of the naval hull.

%% file: sections/intro.tex
In the last years, there has been an increasing interest in using Deep
Neural Networks to approximate distributions of 3D objects \cite{tewari,wing,cheng,qtan,wen} via generative models. In different contexts, new methods that define deformation maps of 3d objects with the preservation of some geometrical properties, like the volume \cite{hirota,vonfunck,cervero,hahmann,eisemberg}, are studied. The possibility to enforce geometrical constraints on the domains of computational fluid dynamics simulations is of great interest, especially for industrial and real-world applications. For example, in naval engineering, the fast generation of geometries of naval ships' hulls such that the volume of the submerged part is preserved is fundamental for hydrodynamic stability. On the other hand, biomedical applications that involve patient-specific numerical models often struggle to obtain new valid geometries that also preserve some geometrical properties of interest since the experimental data are scarce.

We try to close this gap by implementing constrained generative models (cGMs) that are able to reproduce distributions of constrained free from deformations (cFFD)~\cite{sederberg}. The main advantage of our novel methodology is that the computational costs are substantially reduced in the predictive phase with respect to classical numerical methods that define volume-preserving maps for example. Another great advantage is the dimension reduction of the space of parameters from the original space of deformations to the latent space of generative models.

\paragraph{Constrain preserving deformations}
One of the first works that introduce volume preserving deformations in computer graphics is Hirota et al.
\cite{hirota} who carried out a variation of Free Form
Deformation that conserves the volume of the original domain exploiting an augmented Lagrangian formulation. This argument has been developed further by Hahmann et al.\cite{hahmann}, that take
a restriction on the possible deformations and obtain an explicit solution
using a sequence of three quadratic programming problems. Von funck \cite{vonfunck} et al. (2006) propose a method for creating
deformations using path line integration of the mesh vertices over divergence-free vector fields, thus preserving the volume. Eisemberg et al. \cite{eisemberg} extend this method for shape interpolation, using the eigenvectors of the Von Neumann-Laplace equation; the interpolation is done using the Karhunen-Lo\'eve expansion. Cerverò et al. \cite{cervero} obtain a volume preserving method using a cage-based
deformation scheme with generalized barycenter coordinates. They also propose a
measure of the local stress of the deformed volume.

\paragraph{Generative models for 3D meshes}
In the last years, there has been an increasing number of new generative models' architectures for 3D mesh deformation. Qtan et al.\cite{qtan} develop a Variational Autoencoder in which the input data is encoded into a rotation invariant mesh difference representation, which is based on the concept of deformation gradient and it is rotation invariant. They also propose an extended model in which the prior distribution can be altered. Ranjan et al.\cite{ranjan} develop a Convolutional Mesh Autoencoder (CoMA) which is based on Chebyscev Spectral Convolution and on a new sampling method which is based on the hierarchical mesh representation. Rana Hanocka et al. \cite{hanocka} introduce some new original convolutional and pooling layers for mesh data points based on their intrinsic geodesic connections. Yuan et al. \cite{yuan} extend the previous works, by implementing a variational autoencoder that works for meshes with the same connectivity but supported on different geometries. Hahner et al. \cite{hahner} develop an autoencoder model for semiregular meshes with different sizes, which have a locally regular connectivity and hierarchical meshing. They are also able to apply the same autoencoder to different datasets. Cheng et al.\cite{cheng} develop a new mesh convolutional GAN model based on \cite{began}.\\

The work is organized as follows. In section~\ref{sec:methods}, we describe the classical free form deformation method along with its variant to impose linear and multilinear constraints in subsection~\ref{subsec:cffd}. In section~\ref{sec:constrained generative models}, the generative models we are going to use are presented. In subsection~\ref{subsec:cgm}, our novel approach to enforce linear and multilinear constraints on generative models is introduced. Since we also show how model order reduction can be performed efficiently in this context, a brief overview of the non-intrusive proper orthogonal decomposition with interpolation method (PODI) is described in section~\ref{sec:mor} along with radial basis functions interpolation to deform the computational meshes. We test the new methodology on two academic benchmarks in section~\ref{sec:results}: a mixed Poisson problem for the 3d Stanford bunny~\cite{turk} \textbf{SB} with fixed barycenter deformations in section~\ref{subsec:bunny} and the incompressible Reynolds' Averaged Navier-Stokes equations with volume-preserving deformations of the Duisburg's test case~\cite{white2019numerical} \textbf{HB} naval bulb in section~\ref{subsec:bulb}.

%% file: sections/methods.tex
Free Form Deformation (FFD) was introduced in~\cite{sederberg}. It is successfully employed in optimal shape design~\cite{tezzele2018dimension} as a technique to geometrically parametrize the domain through a basis of Bernstein polynomials and a set of bounding control points. When applied to deform computational meshes, some challenges related to FFD, include the preservation of the regularity of the mesh for relatively large deformations, the possibly high dimensional parameter space employed and the loss of conservation of geometrical or physical quantities of interest such as the volume. In this work, we address the last two points.

\subsection{Free Form Deformation}
We consider the FFD of 3d shapes in $\mathbb{R}^3$. Let $D=[0,1]\times[0,1]\times[0,1]\subset\mathbb{R}^3$ be the bounding box inside which the deformation is performed and $P=\{\mathbf{P}_{ijk}\}_{i,j,k=1}^{m, n, o}=\{\left[\frac{i}{m},\frac{j}{n},\frac{k}{o}\right]\}_{i,j,k=1}^{m, n, o}$ the lattice of control points with $m,n, o>0$.

The bases employed for the interpolation are the trivariate Bernstein polynomials $\{B_{ijk}^{\lambda\mu\nu}\}_{i,j,k=0}^{\lambda\mu\nu}\subset\mathbb{P}^{\boldsymbol{\alpha}}(D)$ of polynomial degree $\boldsymbol{\alpha}=(\lambda,\mu,\nu)\in\mathbb{N}^3$ and support $D\subset\mathbb{R}^3$ defined from the one-dimensional basis $\{b_{s}^{\kappa}\}_{s=0}^{\kappa}\subset\mathbb{P}^{\kappa}([0,1])$:
\begin{align}
    b_{s}^{\kappa}(x)&=\binom{\kappa}{s} x^s(1-x)^{\kappa-s},\qquad \{b_{s}^{\kappa}\}_{s=0}^{\kappa}\subset\mathbb{P}^{\kappa}([0,1]),\\
    B_{ijk}^{\lambda\mu\nu}(u, v, w) &= b_{i}^{\lambda}(u)b_{j}^{\mu}(v)b_{k}^{\nu}(w),\qquad \{B_{ijk}^{\lambda\mu\nu}\}_{\nu=0}^{n}\subset\mathbb{P}^{\boldsymbol{\alpha}}(D).
\end{align}
Other bases such as B-splines can be employed~\cite{hahmann}.

Given the set of displacements of the control points $\delta P=\{\delta \mathbf{P}_{ijk}\}_{i,j,k=1}^{m, n, o}$, the deformation map $T_P:D\subset\mathbb{R}^3\rightarrow\mathbb{R}^3$ is defined as follows:
\begin{equation}
    \label{eq:FFD_map}
    T_{Q}(u,v,w)= (u,v,w)+\sum_{i,j,k=0}^{m,n,o}b_{i}^{m}(u)b_{j}^{n}(v)b_{k}^{o}(w)\delta \mathbf{P}_{ijk},\quad\forall (u,v,w)\in D\subset\mathbb{R}^3
\end{equation}
mapping the domain $D=[0,1]\times[0,1]\times[0,1]\subset\mathbb{R}^3$ onto the bounding box of control points $T_P(D)=K$. The displacement field defined through the displacements of the control points $\delta P$ is interpolated inside $D\subset\mathbb{R}^3$ by the second term of Equation~\eqref{eq:FFD_map}.

To apply the FFD deformation, represented by the map $T_P$, to arbitrary 3d meshes, point clouds (CAD or STL files) or general subdomains $\Omega\subset\mathbb{R}^3$, an affine map $\varphi:\mathbb{R}^3\rightarrow\mathbb{R}^3$ is evaluated to map the subdomain $D\subset\mathbb{R}^3$ into a parallelepiped $\Tilde{K}\subset\mathbb{R}^3$ that intersects such 3d meshes, point clouds or general subdomains, $\Tilde{K}\cap\Omega\neq\varnothing$. In this way, we define the FFD map as the composition $\Tilde{T}_P=\varphi \circ T_P\circ\varphi^{-1}:\Tilde{K}\cap D\subset\mathbb{R}^3\rightarrow\mathbb{R}^3$:
\begin{center}
    \begin{tikzcd}
        \Tilde{K}\cap \Omega\subset\mathbb{R}^3 \ar[d,"\varphi^{-1}"] & \varphi(K)\subset\mathbb{R}^3 \\
        D\subset\mathbb{R}^3 \ar[r,"T_P"] & T(D)=K\subset\mathbb{R}^3 \ar[u,"\varphi"] \\
    \end{tikzcd}
\end{center}
notice that, generally, the deformation can affect the boundaries $\partial\left(\Tilde{K}\cap \Omega\right)$ of the intersection $\Tilde{K}\cap \Omega\subset\mathbb{R}^3$ and not only its interior. We will thus employ it to deform 3d objects that will be embedded afterward into a computational mesh, see section~\ref{sec:mor}. The control points $P\subset D\subset\mathbb{R}^3$ are mapped by $\varphi$ into $\Tilde{P}=\{\varphi(\mathbf{P}_{ijk})\}_{i,j,k=1}^{m, n, o}\subset\Tilde{K}\cap\Omega\subset\mathbb{R}^3$ and deformed from $\Tilde{P}\subset\Omega$ to $\tilde{P}+\delta \tilde{P}$, that is $\forall i\in\{0,\dots,m\},j\in\{0,\dots,n\},k\in\{0,\dots,o\}$:
\begin{equation}
    \Tilde{T}_P(\Tilde{\mathbf{P}}_{ijk})=\Tilde{\mathbf{P}}_{ijk}+\delta \Tilde{\mathbf{P}}_{ijk},\qquad
    \delta\Tilde{\mathbf{P}}_{ijk}=\sum_{i,j,k=0}^{m,n,o}b_{i}^{m}(u)b_{j}^{n}(v)b_{k}^{o}(w)A_{\varphi}(\delta \mathbf{P}_{ijk}),
\end{equation}
where $A_{\varphi}\in\mathbb{R}^3\times 3$ is the matrix of the affine transformation $\varphi(\mathbf{x})=A_{\varphi}\mathbf{x}+b_{\varphi}$ for all $\mathbf{x}\in\mathbb{R}^3$. In practice, we only need to fix the control points $\Tilde{P}$ and their deformations $\delta\Tilde{P}$. For this basic form of FFD we use the open-source Python package~\cite{tezzele2021pygem}.

Fixed the control points $\Tilde{P}$ and the displacements $\delta\Tilde{P}$, the map $\Tilde{T}_P:\Tilde{K}\cap D\subset\mathbb{R}^3\rightarrow\mathbb{R}^3$ can also be defined similarly to what has been done previously. We make explicit the dependency on the set of displacements $\delta P$ with the notation
\begin{equation}
    \Tilde{T}_P(\mathbf{Q}, \{\delta\mathbf{P}_{ijk}\}_{ijk=0}^{m,n,o})=\mathbf{Q}+\sum_{i,j,k=0}^{m,n,o}b_{i}^{m}(Q_x)b_{j}^{n}(Q_y)b_{k}^{o}(Q_z)A_{\varphi}(\delta \mathbf{P}_{ijk}),\qquad \Tilde{T}_P(\mathbf{Q}, \delta P):\mathbb{R}^3\times\mathbb{R}^{(m\cdot n\cdot o)\times 3}\rightarrow\mathbb{R}^3
    \label{eq:dep_displacements}
\end{equation}
with $\mathbf{Q}=(Q_x, Q_y, Q_z)\in\mathbb{R}^3$.

\subsection{Free Form Deformation with multilinear constraints}
\label{subsec:cffd}
One of the main objectives of this work is to approximate distributions of point clouds along with some of their geometrical properties such as the volume or the barycenter if they are kept constant. If we want to employ generative models for this purpose, in section~\ref{sec:constrained generative models}, we must first collect a dataset of geometries that satisfy the constraint we want to impose.

Regarding volume preserving deformations, many strategies involve the definition of divergence-free vector fields~\cite{vonfunck, eisemberg}. Other volume-preserving methods are directly applied to FFD, like \cite{hirota} with an augmented Lagrangian nonlinear multi-level formulation and \cite{hahmann} with a least-squares formulation. We will follow the last work for the simple closed-form enforcement of linear and multilinear constraints. Multiple linear constraints can be enforced at the same time as well with this method. In general, for the whole methodology presented in section~\ref{sec:constrained generative models} to work, only the constraint itself needs to be linear, how the training dataset representing the distribution to be approximated with GM is obtained is not relevant. In section~\ref{sec:constrained generative models} are reported other linear and multilinear constraints that could be imposed apart from the volume and the barycenter that are effectively conserved with our procedure as shown in the results section~\ref{sec:results}.

Let us suppose we are given a constraint that is linear with respect to the displacements $\delta P=\{\delta \mathbf{P}_{ijk}\}_{i,j,k=1}^{m, n, o}$. Chosen a subset of $N$ points $Q\in\mathbb{R}^{N\times 3}$ of the undeformed subdomain $Q=\{\mathbf{Q}_i\}_{i=0}^{N}\subset\Omega\subset\mathbb{R}^3$  $\Omega\subset\mathbb{R}^3$ with $\mathbf{Q}_i\in\mathbb{R}^3,\ \forall i\in\{0,\dots,N\}$, we can write the linear constraint as
\begin{subequations}
    \begin{equation}
        \mathbf{c}=A_c\text{vec}(\Tilde{T}_P(Q, \delta P))
    \end{equation}
    where the displacements $\delta P$ of the control points $P$ have been made explicit through the notation of equation~\eqref{eq:dep_displacements} and $\text{vec}(\Tilde{T}_P)\in\mathbb{R}^{3N}$ is the rowwise vectorization of the matrix $\Tilde{T}_P\in\mathbb{R}^{N\times 3}$. We have introduced the vector $\mathbf{c}\in\mathbb{R}^n$ and the matrix $A_c\in\mathbb{R}^{n\times 3N}$ representing the linear constraint. We can expand the previous equation in
    \begin{equation}
        \mathbf{c}=A_c\text{vec}(\Tilde{T}_P(Q, \delta P))=A_c \text{vec}(Q)+A_c \text{vec}(\delta Q),\qquad
        \delta\mathbf{Q}_l = \sum_{i,j,k=0}^{m,n,o}b_{i}^{m}(Q^l_x)b_{j}^{n}(Q^l_y)b_{k}^{o}(Q^l_z)A_{\varphi}(\delta \mathbf{P}_{ijk})
    \end{equation}
\end{subequations}
where $\{\delta\mathbf{Q}_l\}_{l=1}^{N}=\delta Q=\in\mathbb{R}^{N\times 3}$ with $\delta\mathbf{Q}_l\in\mathbb{R}^3$ and $\mathbf{Q}_l=(Q^l_x, Q^l_y, Q^l_z)\in\mathbb{R}^3$, for all $l\in\{1,\dots.N\}$. The linearity with respect to $\delta P =\{\delta\mathbf{P}_{ijk}\}_{i,j,k=0}^{m, n, o}$ is thus made explicit.

A perturbation of the displacement of the control points $\delta d_{\text{cFFD}}=\{\delta\mathbf{d}_{ijk}\}_{i,j,k=1}^{m, n, o}\subset\mathbb{R}^3$ is found solving the least-squares problem
\begin{equation}
    \label{eq:cFFD}
    \delta d_{\text{cFFD}} = \argmin_{\delta d=\{\delta\mathbf{d}_{ijk}\}_{i,j,k=1}^{m, n, o}\subset\mathbb{R}^3}\ \lVert \delta d\rVert_2 \qquad\text{such that}\qquad\mathbf{c}=\mathbf{A}_c\text{vec}(\Tilde{T}_P(\mathbf{Q}, \delta P+\delta d)),
\end{equation}
that can be effectively solved in closed form.

If the constraints are linear in each component (x, y, z) they can be imposed component-wise finding first the x-components $\delta d_{\text{cFFD}, x}$ of the perturbations $\delta d_{\text{cFFD}}=\{\delta\mathbf{d}_{ijk}\}_{i,j,k=1}^{m, n, o}\subset\mathbb{R}^3$ and subsequently the y- and z-components,  $\delta d_{\text{cFFD}, y}$ and  $\delta d_{\text{cFFD}, z}$, respectively.

This strategy is successfully applied for constraints on the volumes of triangulations that define 3d objects in~\cite{hahmann}. Other linear constraints that can be imposed are the position of the barycenter (linear) and the surface area (bi-linear).

Sometimes additional constraints on the position of the control points must be enforced: for example for the Duisburg test case \textbf{HB} of section~\ref{subsec:bulb} we deform only the region of the bulb of a whole ship's hull. This region is extracted from the hull's STL file and two straight cuts are thus introduced. To keep the displacements null close to these cuts additional constraints must be enforced on the cFFD deformations. We do so with the employment of a weight matrix $M=\{\omega_{ijk}\}_{i,j,k=1}^{m, n, o}\subset\mathbb{R}\cap\{x\geq 0\}$ that multiplies $\delta d$:
\begin{equation}
    \label{eq:McFFD}
    \delta d_{\text{cFFD}} = \argmin_{\delta d=\{\delta\mathbf{d}_{ijk}\}_{i,j,k=1}^{m, n, o}\subset\mathbb{R}^3}\ \lVert M\delta d\rVert_2 \qquad\text{such that}\qquad\mathbf{c}=\mathbf{A}_c\text{vec}(\Tilde{T}_P(\mathbf{Q}, \delta P+\delta d)),
\end{equation}
with $M\delta d = \{\omega_{ijk}\delta\mathbf{d}_{ijk}\}_{i,j,k=1}^{m, n, o}$.

%% file: sections/cgm.tex
In the following sections, we introduce the architectures we employ for generative modelling (GM) and the novel constrained generative models (cGMs).

\subsection{Generative modelling}
\label{subsec:gm}
Generative models (GMs)~\cite{goodfellow2016deep} have been successfully applied for computer graphics' tasks such as 3d objects generation. Their employment in surrogate modelling is especially beneficial when the computational costs needed to create a new geometry from real observations is reduced. One disadvantage is that depending on the complexity of the distribution to be approximated, their training requires a large amount of data. In this work, we focus on studying the possible employment of GMs for the reduction of the space of the parameters that specify geometrical deformations. In our case, the parameters are the displacements of cFFD. Another objective is to reduce the computational costs needed to generate new linear constrained geometries with cFFD.

We denote with $p_{X}(x)$ the probability density of the distribution of 3d objects and with $X:(A, \mathcal{A}, P)\rightarrow\mathbb{R}^{M}$ the random variable representing it, with $M$ the number of points or degrees of freedom of the 3d mesh. The triple $(A, \mathcal{A}, P)$ corresponds to the space of events, sigma algebra and probability measure of our setting. The variable $x$ may represent 3d point clouds, 3d meshes or general subdomains of $\mathbb{R}^3$. We are going to focus on generative models that approximate $p_{X}(x)$ with $p_{\theta}(x)$ where $\theta$ are parameters to be found and that can be factorized through latent variables $z$ such that 
\begin{equation}
    p_{\theta}(x)=\int_{Z} p_{\theta}(x|z)p_{\theta}(z)dz
\end{equation}
where $z$ is typically low dimensional with respect to $x$, with dimension $R\ll M$. These models are called latent variable models\cite{tomzack}.
Given a set of training samples of the distribution associated to $X$ obtained with cFFD (section~\ref{subsec:cffd}), the objective is to sample new 3d objects with quality comparable to the one generated using cFFD. It is crucial that the new samples satisfy exactly the linear or multilinear constraints imposed with cFFD.

As the subject of generative models is well-known in the literature we summarize briefly the architectures we are going to employ.

\paragraph{Simple Autoencoder}
\label{par:ae}
As a toy model for pure benchmarking, we first implement a very simple Autoencoder \cite{ae}, composed of two parts: an Encoder $Enc_{\psi}:\mathbb{R}^{M}\rightarrow \mathbb{R}^{R}$ which encodes the mesh in the latent space, and a Decoder $Dec_{\theta}: \mathbb{R}^{R} \rightarrow \mathbb{R}^{M}$ that takes a point of the latent space and returns it to the data space. 
Autoencoders were originally created for learning the latent representations of the data and have many applications \cite{ae}, like denoising. Some recent research has empirically shown  \cite{ghose,daly} that they can indeed be used as generative models.
They are trained on the $L^2$ loss
\begin{equation}
    L_{\theta,\psi}(x)=||x-Dec_{\theta}(Enc_{\psi}(x))||_{2},
\end{equation}
the associated generative model is
\begin{equation}
    p_{\theta}(x|z)=Dec_{\theta}(z)
\end{equation}
where $z$ is sampled from a multivariate normal distribution.

\paragraph{Beta Variational Autoencoder}
\label{par:vae}

A Variational Autoencoder \cite{varauto} is a probabilistic version of the autoencoder. It is based on the concept of Evidence Variational Lower Bound (ELBO) that can be employed in place of the log-likelihood $p_{\theta}(x)$ when it is not directly computable. Given a variational distribution $q_{\psi}(z|x)$ the ELBO is 

\begin{equation}
    ELBO_{\psi,\theta}(x)= E_{q_\phi(z \mid x)}\left[\log \left(\frac{p_\theta(x, z)}{q_\phi(z \mid x)}\right)\right]= E_{q_\phi(z \mid x)}\left[\log \left(p_\theta(x \mid z)\right)\right]-K L\left(q_\phi(z \mid x) \| p_\theta(z)\right)
\end{equation}
it follows that
\begin{equation}
    ELBO_{\psi,\theta}(x)\le \log p_{\theta}(x).
\end{equation}

Variational autoencoders model $q_{\psi}(z|x)$ and $p_{\theta}(x|z)$ as two normal distributions:
\begin{equation}
    q_{\psi}(z|x)=\mathcal{N}(a_{\psi}(x),b_{\psi}(x))
\end{equation}
where $a_{\psi},b_{\psi}:\mathbb{R}^{M}\rightarrow \mathbb{R}^{R}$ are encoders and
\begin{equation}
    p_{\theta}(x|z)=\mathcal{N}(Dec_{\theta}(z),\sigma),
\end{equation}
Notice that, the $KL$ term in the ELBO is a regularizing term on the variational distribution.
Depending on the strength of $p_{\theta}(x|z)$ the Variational autoencoders can suffer from two problems: if it is dominant then there is no regularizing effect and so the samples from $p_{\theta}(z)$ will have poor quality. Otherwise, the $KL$ will go to $0$, so $a_{\psi},b_{\psi}$ will degenerate to constants and the model will start behaving autoregressively. To solve these two problems, an $\alpha$ term is added to the ELBO:
\begin{equation}
    L_{\theta,\phi}=E_{q_\phi(z \mid x)}\left[\log \left(p_\theta(x \mid z)\right)\right]-\alpha K L\left(q_\phi(z \mid x) \| p_\theta(z)\right),
\end{equation}
this technique is also used in \cite{qtan}.

\paragraph{Adversarial Autoencoder}
\label{par:aae}

Another class of Generative Models is Generative Adversarial Networks which have the same structure of an autoencoder
\begin{equation}
    p_{\theta}(x|z)=Dec_{\theta}(z)
\end{equation} 
but they are trained with the help of a discriminator $D_{\lambda}$ in the following min-max game:
\begin{equation}
    \min _\theta \max _\lambda V(\theta, \lambda)=E_{p_X(x)}\left[\log \left(D_\lambda(x)\right)\right]+E_{p_\theta(z)}\left[\log \left(1-D_\lambda\left(Dec_\theta(z)\right)\right)\right].
\end{equation}
In this context, the decoder is also called generator.
The discriminator is used to distinguish between the data and the generator samples while the generator wants to fool it. A problem of generative adversarial networks is that the min-max game can bring instabilities, for the details of this we refer to \cite{arjovsky}. A solution is to combine them with variational autoencoders to get the following loss:
\begin{equation}
    \max _{\theta, \phi} \min _\lambda L_{\theta, \phi, \lambda}=E_{q_\phi(z, x)}\left[\log \left(p_\theta(x \mid z)\right)\right]-
E_{p_\theta(z)}\left[\log \left(1-D_\lambda(z)\right)\right]+E_{q_\phi(z)}\left[\log \left(D_\lambda(z)\right)\right],
\end{equation}
which is the ELBO with the GAN loss instead of the KL term.
This model is called Adversarial Autoencoder \cite{aae}. The associated generative model is
\begin{equation}
    p_{\theta}(x|z)=\mathcal{N}(Dec_{\theta}(z),\sigma), \quad\text{with}\quad z\sim\mathcal{N}(0, I_R)
\end{equation}

\paragraph{Boundary Equilibrium GAN}
\label{par:began}

Another solution to the GAN problems is the Boundary Equilibrium GAN \cite{began}. In this particular model, the discriminator is an autoencoder
\begin{equation}
    D_{\lambda}(x)=Dec_{\lambda}(Enc_{\lambda}(x))
\end{equation}
and the objective is to learn the autoencoder loss function
\begin{equation}
    f_{\lambda}=||x-D_{\lambda}(x)||_{2}
\end{equation}
while maintaining a balance between the data loss and the generator loss to prevent instabilities caused by the min-max training. Let $G_{\theta}$ the generator, this is achieved in three steps using control theory:
\begin{itemize}
    \item Solve $\min_{\lambda} E_{p_{X}(x)}[f_{\lambda}(x)]-k_{t} E_{p_{X}(x)}[f_{\lambda}(G_{\theta}(Enc_{\lambda}(x)))]$ over $\lambda$.
    \item Solve $\max_{\theta} E_{p_{Z}(z)}[f_{\lambda}(G_{\theta}(z)].$
    \item Update $k_{t}=k_{t-1}+a(\gamma E_{p_{X}(x)}[f_{\lambda}(x)]- E_{p_{Z}(z)}[f_{\lambda}(G_{\theta}(z)]).$
\end{itemize}

The associated generative model is
\begin{equation}
    p_{\theta}(x|z)=G_{\theta}(z), \quad\text{with}\quad z\sim\mathcal{N}(0, I_R).
\end{equation}

\subsection{Generative modelling with multilinear constraints}
\label{subsec:cgm}
In this section, we introduce our new framework to impose linear or multilinear constraints on GMs. The idea is to rely on the generalization properties of NNs. The cGMs layers' objective is to approximate the training distributions while correcting the unwanted deformations that do not satisfy exactly the geometrical constraints. multilinear constraints are enforced subsequently in the same way as linear ones. We will show that the volume and the barycenter position of 3d objects can be preserved in our numerical experiments~\cref{sec:results}.

The best architectures that adapt to datasets structured on meshes are graph neural networks~\cite{bronstein2017geometric}. Unfortunately, compared to convolutional neural networks for data structured on Cartesian grids, they are quite heavy to train for large problems supported on computational meshes. Thus, we leave to further studies the implementation of our methodology with GNNs and we focus on generative models that employ principal component analysis (PCA) to perform dimension reduction as a preprocessing step. After having collected a training and test dataset of 3d point clouds,
\begin{equation}
    \mathcal{X}_{\text{train}}=\begin{pmatrix}
    \mid & \mid & \mid & \mid\\
    \mathbf{x}_1 & \mathbf{x}_2 & \dots & \mathbf{x}_{n_{M\times\text{train}}} \\
    \mid & \mid & \mid & \mid
  \end{pmatrix}^T\in\mathbb{R}^{n_{\text{train}}\times M},\qquad
  \mathcal{X}_{\text{test}}=\begin{pmatrix}
    \mid & \mid & \mid & \mid\\
    \mathbf{x}_1 & \mathbf{x}_2 & \dots & \mathbf{x}_{n_{M\times\text{test}}} \\
    \mid & \mid & \mid & \mid
  \end{pmatrix}^T\in\mathbb{R}^{n_{\text{test}}\times M},
\end{equation}
PCA is applied to obtain a set of $r$ modes of variation $U_{\text{PCA}}\in\mathbb{R}^{M\times r}$ with $M \gg r>0$, such that fixed a tolerance $1\gg\epsilon>0$ for the reconstruction error in Frobenious norm $\lVert\cdot\rVert_F$, we have:
\begin{equation}
    \lVert (I_{M}-U_{\text{PCA}}U_{\text{PCA}}^T)\mathcal{X}_{\text{train}}\rVert_{F}\leq \epsilon.
\end{equation}
The notation $M$ for the dimension of a single point cloud includes the $3$ (x,y,z) components of each point stacked on each other on the same vector. Without considering the enforcement of the geometrical constraints, the GMs could now be trained to approximate the lower dimensional distributions $U_{\text{PCA}}^T\mathcal{X}_{\text{train}}\in\mathbb{R}^{r\times n_{\text{train}}}$ with a great saving in terms of computational cost since $M\gg r$. However, this reasoning may not work with more complex distributions that cannot be accurately approximated as combinations of PCA modes. A direct consequence would be that the generalization error on the test dataset would perform significantly worse than the training error. Fortunately, this is not the case of our numerical studies: $r=140$ and $r=30$ PCA modes are sufficient for the Stanford's bunny and DTCHull test cases, respectively.

We denote with $Y=U_{\text{PCA}}X$ the random variable associated to the dataset at hand projected on the PCA's modes. With $\mathbf{y}\in\mathbb{R}^r$ we represent its realizations. To reconstruct the point clouds on the full space, a matrix multiplication $U_{\text{PCA}}\mathbf{y}=\Tilde{\mathbf{x}}\in\mathbb{R}^M$ is needed. We denote with $\Tilde{X}:(A, \mathcal{A}, P)\rightarrow\mathbb{R}^M$ the random variable associated to the GMs, that is used to approximate $X$ without the enforcement of geometrical constraints. 

Our methodology to impose linear or multilinear constraints on GMs, affects both the training and predictive stages. It consists in the addition of a final constraint enforcing layer $l_{\text{enforcing}}:\mathbb{R}^M\rightarrow\mathbb{R}^M$ that acts on the reconstructed outputs $\Tilde{\mathbf{x}}\in\mathbb{R}^M$. The final layer $l_{\text{enforcing}}$ perturbs the outputs $\Tilde{\mathbf{x}}+\delta\Tilde{\mathbf{x}}=\Tilde{\Tilde{\mathbf{x}}}\in\mathbb{R}^M$ such that the geometrical constraints with fixed values represented by the vector $\mathbf{c}\in\mathbb{R}^{n_c}$ are exactly satisfied through the solution of the least-squares problem
\begin{equation}
    \delta\Tilde{\mathbf{x}}=\argmin_{\delta\mathbf{x}\in\mathbb{R}^M} \lVert\delta\mathbf{x}\rVert_2,\qquad\text{s.t.}\qquad A_c(\Tilde{\mathbf{x}}+\delta\mathbf{x})=\mathbf{c},
\end{equation}
where $A_c\in\mathbb{R}^{n_c\times M}$ is the matrix representing the linear geometrical constraints and $\delta\Tilde{\mathbf{x}}\in\mathbb{R}^M$ are the point clouds perturbations. multilinear constraints, like the volume, are imposed subsequently component after component. The final outputs $\Tilde{\Tilde{\mathbf{x}}}\in\mathbb{R}^M$ that satisfy exactly the geometrical constraints are then forwarded to the loss function during the training or used directly in the predictive phase.

\begin{figure}[htpb!]
    \centering
    \includegraphics[width=0.8\textwidth]{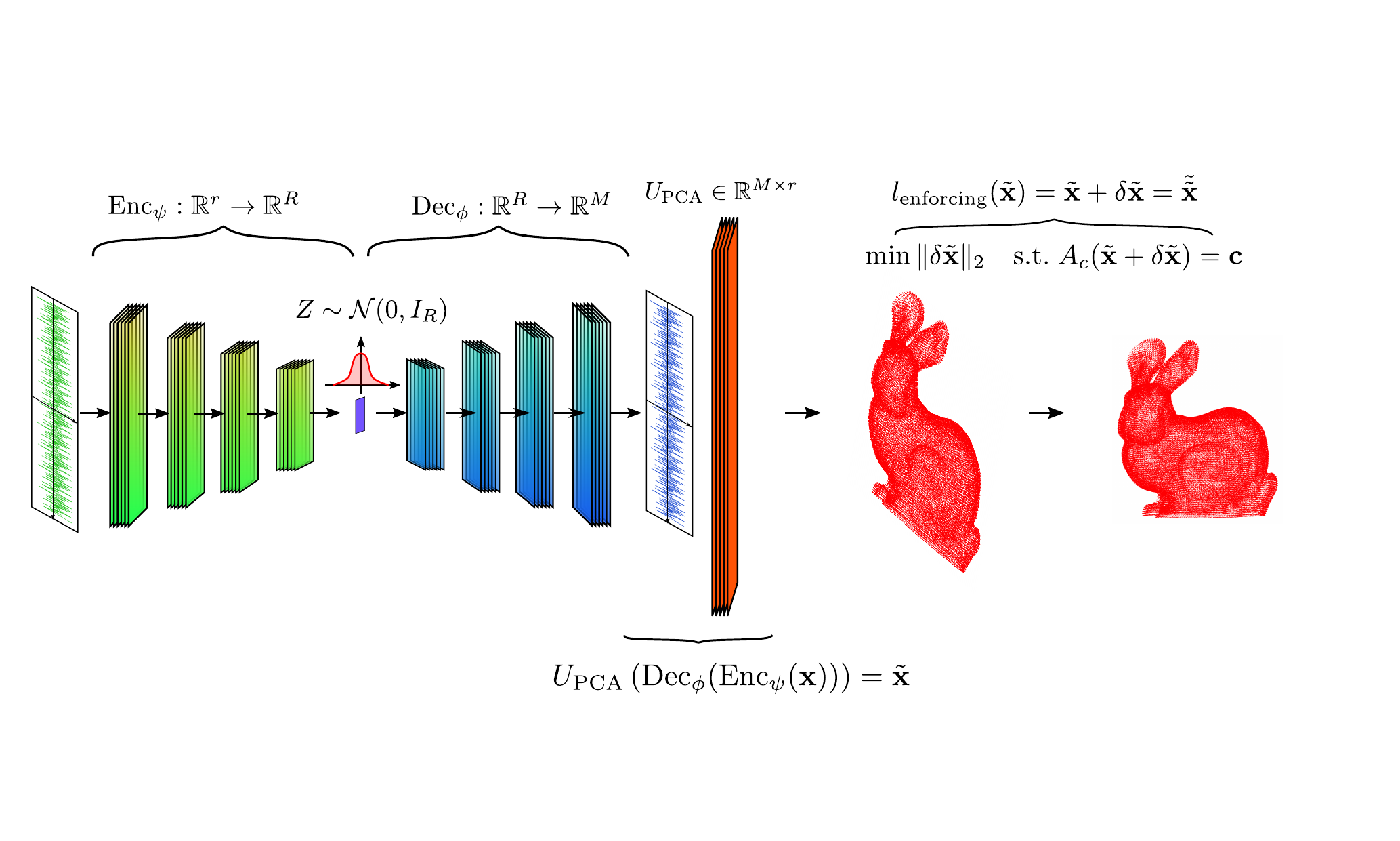}
    \caption{Training of the constrained autoencoder generative model composed by a parametrized encoder $\text{Enc}_{\psi}$, a decoder $\text{Dec}_{\phi}$, a reconstruction layer that employs PCA modes $U_{\text{PCA}}\in\mathbb{R}^{M\times r}$ and a final multilinear constraint enforcing layer $l_{\text{enforcing}}$. The generalization to other generative models of is straight-forward.}
    \label{fig:enforcing_l}
\end{figure}

It is crucial to remark that the linear constraints are imposed exactly on the outputs of the GMs, after solving a least-squares problem similar to equation~\ref{eq:cFFD}. The main difference is that the perturbations affect directly the cloud of points coordinates and not the displacements of the control points as in the constrained FFD presented in section~\ref{subsec:cffd}. Figure~\ref{fig:enforcing_l} shows a sketch of the methodology, including the application of PCA and of the constraints enforcing layer, for the simple autoencoder GM. The procedure is easily adapted for the other GMs of the previous section~\ref{subsec:gm}.

Despite this quite arbitrary correction $\Tilde{\mathbf{x}}+\delta\Tilde{\mathbf{x}}=\Tilde{\Tilde{\mathbf{x}}}\in\mathbb{R}^M$ during the training of the GMs, the NNs' layers manage to correct the outputs $y\in\mathbb{R}^r$ to balance the PCA reconstruction and the perturbation $\delta\Tilde{\mathbf{x}}\in\mathbb{R}^M$. The simplicity of this strategy together with the generalization capabilities of GMs make our methodology effective and relatively easy to implement. We remark that the methodology is not strictly linked with FFD and cFFD: arbitrary techniques, like~\cite{vonfunck, eisemberg}, can be employed to generate the training and test datasets as long as the geometrical constraints are linear or multilinear.

Our method cannot be applied to nonlinear constraints. A possible way to implement them is to use substitute the linear constraints enforcing layer with nonlinear optimization layers~\cite{cvxpylayers2019}. The disadvantage is the higher computational cost and the fact that the nonlinear constraints would not be satisfied exactly.

%% file: sections/mor.tex
In our numerical experiments in section~\ref{sec:results}, we validate the distributions of the trained cGMs with some geometrical and physical metrics of interest. To evaluate these metrics we need to compute some physical fields with numerical simulations on computational domains affected by the newly sampled 3d objects from the cGMs. A Poisson problem is considered in section~\ref{subsec:bunny} and the Navier-Stokes equations in section~\ref{subsec:bulb}.

The first test case in section~\ref{subsec:bunny} considers as domain the deformed Stanford bunny~\cite{turk} with fixed barycenter location. The second test case in section~\ref{subsec:bulb} considers as domain a parallelepiped with two separated phases (water and air) in which is embedded a ship hull whose bulb is deformed with cGMs while keeping its volume fixed.

We show how reduced order modelling~\cite{hesthaven2016certified, rozza2022advanced} can benefit from the employment of cGMs. The main advantage is dimension reduction of the parameter space. In fact, the space of parameters that define the geometrical deformations changes from a possibly high-dimensional space associated with the cFFD method to the usually smaller latent space of the cGMs. The consequence of this dimension reduction is the increased efficiency of non-intrusive ROMs based on interpolation methods from the space of parameters to the coefficients of proper orthogonal decomposition (POD) modes.

The outputs of cGMs are 3d point clouds organized into STL files. To obtain a mesh from each STL we use interpolation with radial basis functions (RBF)~\cite{de2007mesh} of a reference computational mesh: the 3d point clouds generated by the cGM are used to define an interpolation map that will deform the reference computational mesh. We proceed in this way because it is simpler to design ROMs on computational meshes with the same number of degrees of freedom. In fact, with RBF interpolation, each mesh generated from a STL file has the same number of cells of the reference mesh. A possible solution is to evaluate projection and extrapolation maps from a fine mesh common to all the others associated with different sampled 3d objects with possibly different numbers of dofs. Since MOR is not our main focus, we employ RBF interpolation, knowing that in this way, the geometrical constraints are not imposed exactly anymore but are subject to the level of discretization of the mesh and to the accuracy of the RBF interpolation. 

In this section, we briefly summarize the employment of RBF interpolation and present the non-intrusive ROMs we use in the results' section~\ref{sec:results}: proper orthogonal decomposition with interpolation (PODI) performed by Gaussian process regression (GPR)~\cite{rasmussen2003gaussian}, radial basis functions (RBF) or feed-forward neural networks (NNs)~\cite{goodfellow2016deep}. The active subspaces method (AS)~\cite{constantine2015active} is also introduced as reference dimension reduction method for the space of parameters and coupled with PODI model order reduction~\cite{tezzele2018combined}.

\begin{figure}[H]
  \centering
  \includegraphics[width=0.9\textwidth]{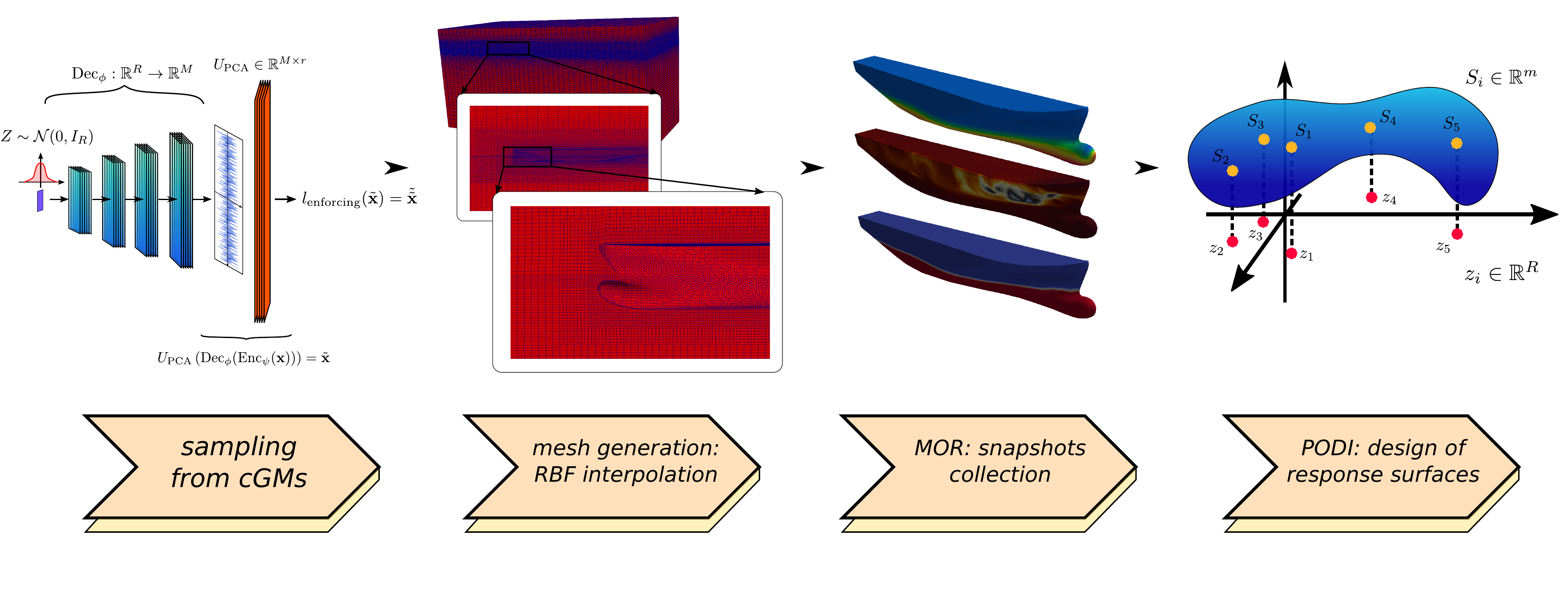}
 \caption{Pipeline for non-intrusive reduced order modelling: from 3d object generation with constrained generative models (cGMs) to mesh deformation with radial basis functions (RBF) to surrogate modelling with proper orthogonal decomposition with interpolation (PODI).}
\end{figure}

\subsection{Mesh interpolation with radial basis functions}
\begin{figure}[H]
  \centering
  \includegraphics[width=0.9\textwidth]{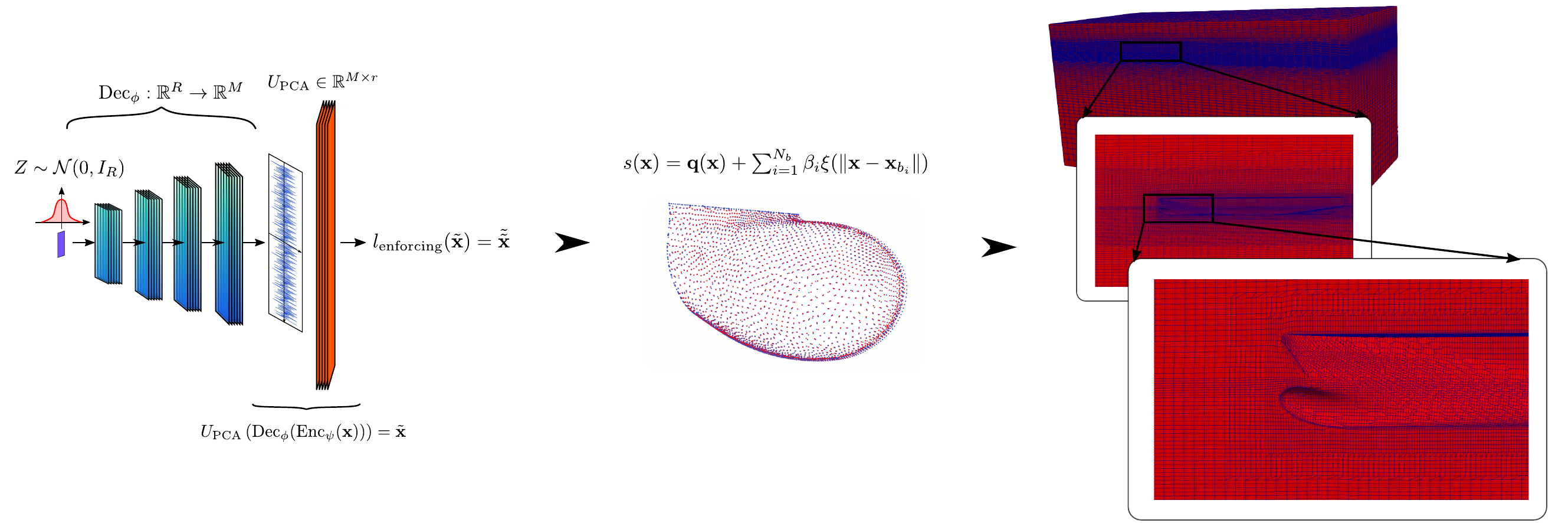}
 \caption{Pipeline for the deformation of the reference computational mesh for the Duisburg test case \textbf{HB}. The numerical results are reported in section~\ref{subsec:bulb}. In red the 3d point cloud of the reference STL, in blue the 3d point cloud of a deformed bulb generated by a constrained autoencoder used as cGM on the left. From these two points cloud $x^2$ and $x^1$ a RBF interpolation map is defined from equation~\eqref{eq:rbf} and used to deform the reference mesh on the right. The bulb location with respect to the whole computational domain is shown.}
 \label{fig:rbf_int}
\end{figure}

After having obtained a new 3d object as STL file from the trained cGMs represented by the 3d points cloud $x^1=\{\mathbf{x}^1_i\}_{i=1}^{\Tilde{M}}\in\mathbb{R}^{3}$, we define a RBF interpolation map $s:\mathbb{R}^3\rightarrow\mathbb{R}^3$,
\begin{equation}
  s(\mathbf{x})=\mathbf{q}(\mathbf{x})+\sum_{i=1}^{N_b} \beta_i \xi(\lVert\mathbf{x}-\mathbf{x}_{b_{i}}\rVert),\qquad\forall\mathbf{x}\in\mathbb{R}^3,
\end{equation}
where $\{\mathbf{x}_{b_{i}}\}_{i=1}^{N_b}\subset\mathbb{R}^3$ are the RBF control points, $N_b$ is the number of RBF control points, $\xi:\mathbb{R}^3\rightarrow\mathbb{R}^3$ is a radial basis function and $\mathbf{q}\in[\mathbb{P}^1(\mathbb{R}^3)]^3$ a vector-valued polynomial of degree $1$ with coefficients $\delta=\{\delta^1_0, \delta^1_1, \delta^2_0, \delta^2_1, \delta^3_0, \delta^3_1\}$ to be defined. The coefficients of the RBF interpolation map $s$ are $\beta=\{\beta_{i}\}_{i=1}^{N_b}\subset\mathbb{R}^3$. The parameters of the RBF interpolation map to be defined are $\beta$ and $\delta$. Given a reference mesh, the STL file from which it was obtained is employed as reference STL, with associated 3d point cloud $x^2=\{\mathbf{x}^2_i\}_{i=1}^{\Tilde{M}}\in\mathbb{R}^{3}$. Since the newly sampled STL files and the reference one have the same number of points $\Tilde{M}$, the following interpolation problem gives $\beta$ and $\delta$:
\begin{subequations}
  \label{eq:rbf}
  \begin{align}
    s(\mathbf{x}_j^2)&=q(\mathbf{x}_j^1)+\sum_{i=1}^{N_b} \beta_i \xi(\lVert\mathbf{x}_j^1-\mathbf{x}_{b_{i}}\rVert),\qquad \forall j\in\{1,\dots,\Tilde{M}\}\\
    0&=\sum_{i=1}^{N_b} \beta_i q(\mathbf{x}_{b_{i}}).
  \end{align}
\end{subequations}
In practice, to each point cloud, $x^1$ and $x^2$ is usually added a set of points to keep fixed such that the reference mesh is only deformed in a limited region of the computational domain. This is the reason why we employed the notation $\Tilde{M}>M$ to differentiate between the dimension $M$ of the 3d point clouds that are the output of the cGMs and the dimension $\Tilde{M}$ of the RBF control points that are enriched with additional points to be kept fixed. An example of reference (in red) and deformed $M=5000$ (in blue) control points from the constrained autoencoder generative model is shown in Figure~\ref{fig:rbf_int}.

We remark that due to the RBF interpolation, the multilinear constraints are not imposed exactly anymore, but only approximately. To avoid this problem, the same discretization of the 3d object embedded in the computational mesh must be used and converted as point cloud also for the training of the cGMs. This is effectively done for the Stanford's bunny test case, but not for the DTCHull's one. In fact, we decided to reduce the resolution of the hull for the numerical simulations in order to lower the computational costs for these preliminary studies.

\subsection{Proper orthogonal decomposition with interpolation}
After having obtained a deformed mesh with RBF interpolation, as for the DTCHull test case, or directly without RBF interpolation, as for the Stanford's bunny test case, numerical simulations can be performed. For each one of the newly sampled geometries from the cGMs, some physical fields of interest are collected and organized in a snapshots matrix $S\in\mathbb{R}^{m\times n_{\text{ROM}, \text{train}}}$, where $m$ is the number of degrees of freedom while the rest $n_{\text{ROM}, \text{test}}=n_{test}-n_{\text{ROM}, \text{train}}$ is used for testing the accuracy of the reduced order models (ROMs). We remind that $n_{\text{test}}$ is the number of newly sampled geometries from the cGMs, divided in $n_{\text{ROM}, \text{train}}$ training and $n_{\text{ROM}, \text{test}}$ test datasets to validate our model order reduction procedure.

We employ the finite volumes method (FVM)~\cite{versteeg2007introduction} implemented in the open-source software library OpenFoam~\cite{weller1998tensorial}, for both test cases in section~\ref{sec:results}. For the Stanford's bunny we set a simple Poisson problem, while for the DTCHull we start from the DTCHull multiphase tutorial~\cite{white2019numerical} for the Reynolds Averaged Navier-Stokes equations (RAS) with associated solver \textit{interFoam}.

Proper orthogonal decomposition (POD), already introduced as PCA previously, is employed to compute a set of reduced basis $U_{\text{ROM}}\in\mathbb{R}^{m\times r_{\text{ROM}}}$. For our studies, we will employ only non-intrusive ROMs: new values of the physical fields are obtained with interpolations or regressions of the expansion with respect to the POD modes $U_{\text{ROM}}$, this technique is generally called POD with interpolation (PODI). The input-output dataset of the interpolation or regression are the latent coordinates of the cGMs $\{z_i\}_{i=1}^{n_{\text{ROM}, \text{train}}}\in\mathbb{R}^R$ or the displacements $\delta P$ of cFFD for the inputs and the coefficients $U_{\text{ROM}}^T S\in\mathbb{R}^{r_{\text{ROM}}\times n_{\text{ROM}, \text{train}}}$ for the outputs:
\begin{subequations}
  \label{eq:inpout}
  \begin{align}
    \{(z_l, S_l)\}_{p=1}^{n_{\text{ROM}, \text{train}}}\subset\mathbb{R}^R\times\mathbb{R}^m,\qquad \text{(inputs-outputs of PODI for the cGMs)}\\
    \{(\delta P^l, S_l)\}_{p=1}^{n_{\text{ROM}, \text{train}}}\subset\mathbb{R}^p\times\mathbb{R}^m,\qquad \text{(inputs-outputs of PODI for the cFFD)}
  \end{align}
\end{subequations}
where $p$ is the number of cFFD displacements $\{\delta P^l\}_{l=1}^p=\{\{\delta \mathbf{P}^l_{i,j,k}\}_{i,j,k=0}^{\Tilde{m}, \Tilde{n}, \Tilde{o}}\}_{l=1}^p$ different from $0$, $\Tilde{m}, \Tilde{n}, \Tilde{o}$ count the non-zero displacements, possibly less than the $m, n, o$ control points of the lattice of FFD. The notation $S_j\in\mathbb{R}^m$ refers to the rows of the snapshots matrix $S\in\mathbb{R}^{m\times n_{\text{ROM}, \text{train}}}$.

Once the interpolation or regression maps are defined from the training input-output datasets, the new $n_{\text{ROM}, \text{test}}$ physical fields associated with the geometries of the test dataset, are efficiently evaluated through these interpolations or regressions, without the need for full-order numerical simulations.

The method we employ to perform the interpolation is RBF interpolation, while Gaussian process regression (GPR) and feed-forward neural networks (NNs) are employed to design a regression of the POD coefficients. These techniques are compared in section~\ref{sec:results}. Also, the Active subspaces method (AS) will be used to build response surfaces from the latent coordinates for the cGMs or the cFFD displacements while performing also a further reduction in the space of parameters.

\subsection{Active subspaces method}
We briefly sketch the Active Subspaces method (AS). It will be used as reference dimension reduction method in the space of parameters and as regression method to perform model order reduction through the design of response surfaces. We will define response surfaces for the ROMs built on top of both cFFD and cGMs deformed meshes: we will see that not only the initial parameter space dimension $p$ associated to the cFFD's non-zero displacements (see equation~\eqref{eq:inpout}) can be substantially reduced for our test cases, but also the latent spaces' dimensions of the cGMs.

Given a function $f:\chi\subset\mathbb{R}^R\rightarrow\mathbb{R}$ from the space of parameters to a scalar output of interest $f:\boldsymbol{\mu}\mapsto f(\boldsymbol{\mu})$, the active subspaces are the leading eigenspaces of the uncentered covariance matrix of the gradients:
\begin{equation}
  \label{eq:covSA}
  \Sigma = E[\nabla_{\boldsymbol{\mu}} (f\nabla_{\boldsymbol{\mu}} f)^T] = \int_{\chi} (\nabla_{\boldsymbol{\mu}} f)(\nabla_{\boldsymbol{\mu}} f)^T\rho\,d\mathbf{\boldsymbol{\mu}},
\end{equation}
where $\rho:\chi\subset\mathbb{R}^R\rightarrow\mathbb{R}$ is the probability density function of the distribution of the inputs $\boldsymbol{\mu}$ considered as a vector-valued random variable in the probability space $(\chi, \mathcal{A}, P)$. The gradients of $f$ are usually approximated with regression methods if they cannot be evaluated directly. The uncentered covariance matrix $\Sigma\in\mathbb{R}^{R\times R}$ is computed approximately with the simple Monte Carlo method and then the eigenvalue decomposition
\begin{equation}
  \Sigma = W\Lambda W^T,\qquad W=[W_1, W_2],\qquad \Lambda=\text{diag}(\lambda_1,\dots,\lambda_R),\qquad W_1\in\mathbb{R}^{R\times r_{\text{AS}}},\ W_2\in\mathbb{R}^{R\times (R-r_{\text{AS}})},
\end{equation}
highlights the active $W_1\in\mathbb{R}^{R\times r_{\text{AS}}}$ and the inactive $W_2\in\mathbb{R}^{R\times (R-r_{\text{AS}})}$ subspaces, corresponding to the first $r_{\text{AS}}$ eigenvalues $\{\lambda_1,\dots,\lambda_{r_{\text{AS}}}\}$ and last $R-r_{\text{AS}}$ eigenvalues $\{\lambda_{r_{\text{AS}}},\dots,\lambda_{R}\}$.

A response surface can be obtained with the approximation
\begin{equation}
  f(\boldsymbol{\mu})\approx g(W_1^T\boldsymbol{\mu})=g(\boldsymbol{\mu_1}),
\end{equation}
where $g:W_1^T(\chi)\subset\mathbb{R}^{r_{\text{AS}}}\rightarrow\mathbb{R}$ is a surrogate model for $f$ from the reduced space of active variables $\boldsymbol{\mu}_1=W_1^T\boldsymbol{\mu}$, instead of the full parameter space. We will employ Gaussian process regression to evaluate $g$.

%% file: sections/results.tex
We present two numerical studies to validate our novel methodology for constrained generative modelling previously introduced in section~\ref{sec:constrained generative models}. The test cases we consider are the Stanford Bunny~\cite{turk} (\textbf{SB}) and the bulb of the hull of the Duisburg test case~\cite{white2019numerical} (\textbf{HB}). After generating new samples from the constrained generative models (cGMs), the crucial task of validating the results must be carried out. In fact, a natural metric that evaluates the quality of the generated distribution of the cGMs is not available: for each problem at hand, we have to decide which criteria, summary statistics, geometrical and physical properties are most useful to validate the generated distribution. For this reason, we decide that the \textbf{SB} test case's generated 3d objects will be employed to solve a Poisson problem with fixed barycenter position and the \textbf{HB} test cases' generated bulbs will be embedded on a larger computational domain to solve the multiphase Navier-Stokes equations with fixed volume of the hull.

\begin{figure}[!htpb]
    \centering
    \includegraphics[width=0.25\textwidth, trim={0 0 700 0}, clip]{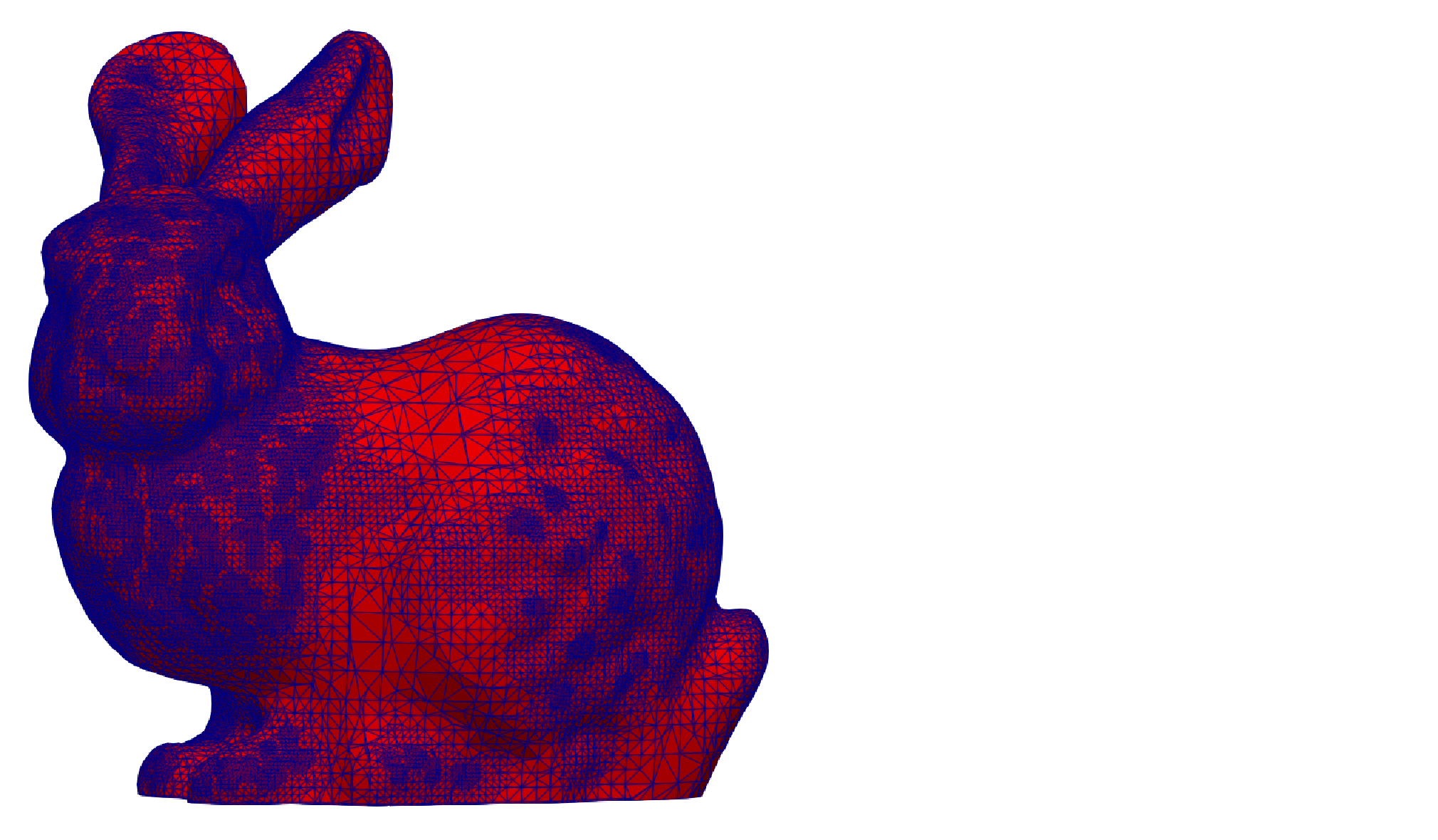}
    \includegraphics[width=0.48\textwidth, trim={0 100 0 0}, clip]{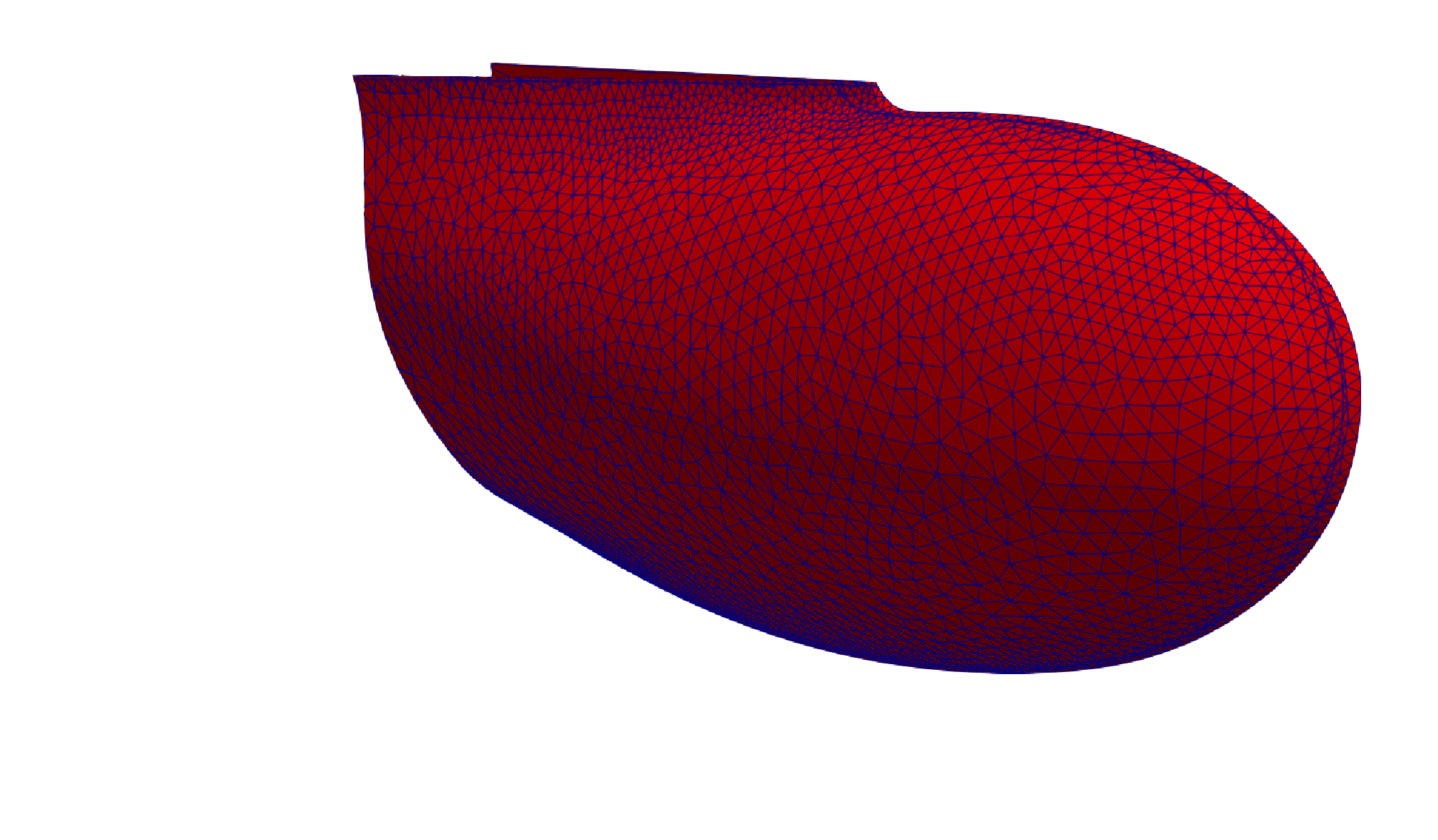}\\
    \includegraphics[width=0.25\textwidth, trim={0 0 700 0}, clip]{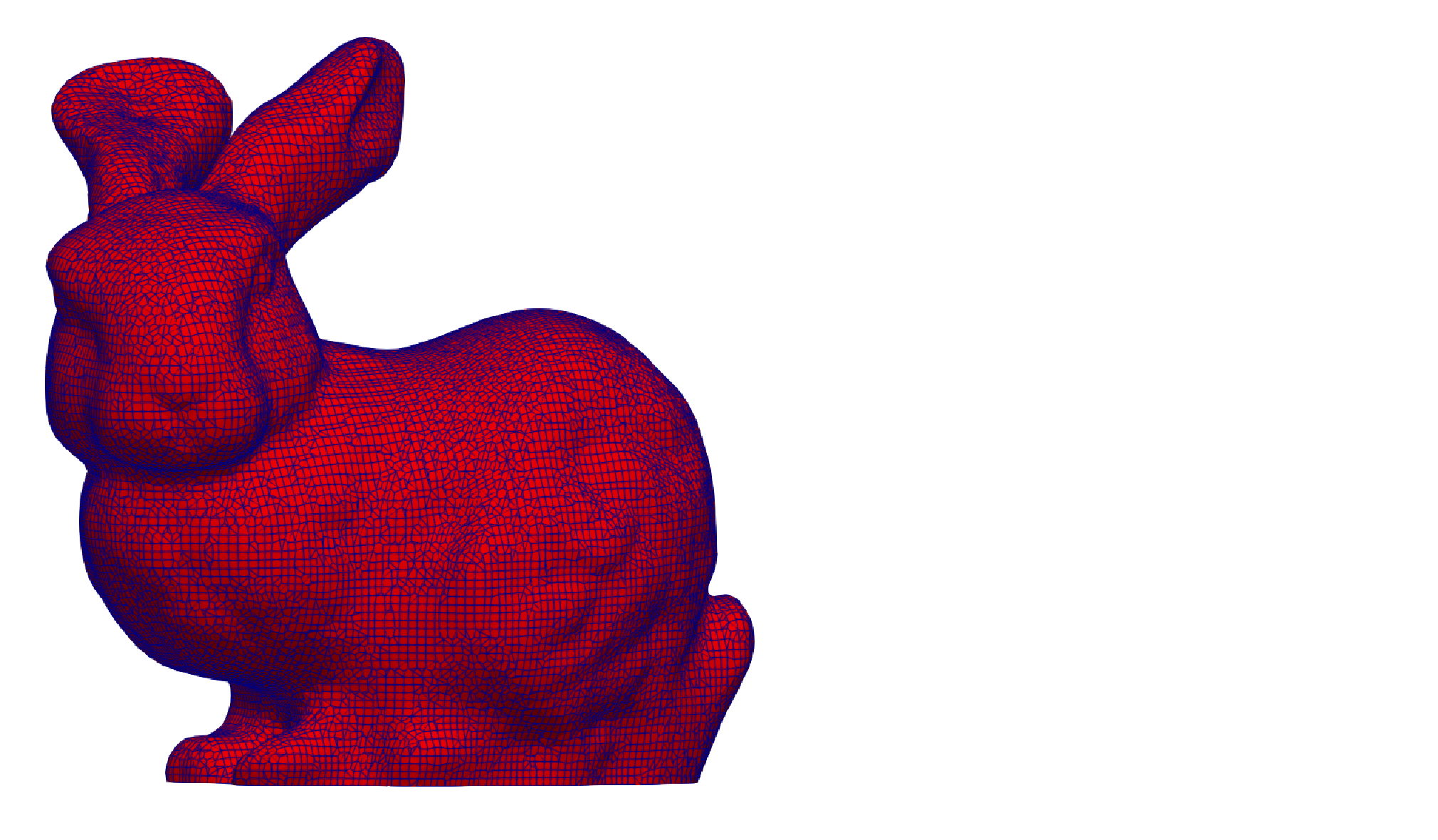}
    \includegraphics[width=0.48\textwidth, trim={400 0 0 300}, clip]{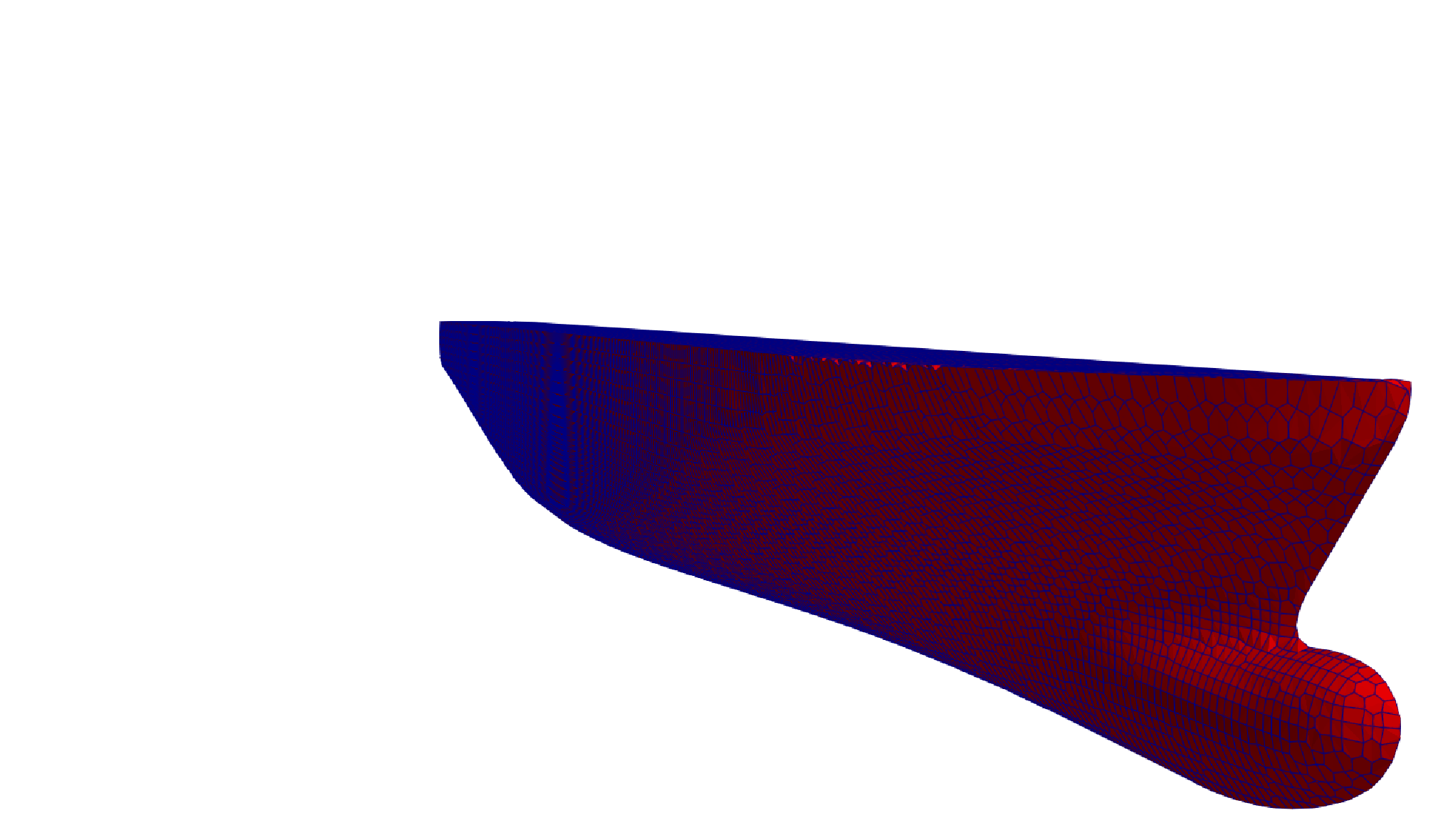}
    \caption{\textbf{Top left:} STL of \textbf{SB} with \textbf{145821} points. \textbf{Top right:} STL of \textbf{HB} with \textbf{5000} points. \textbf{Bottom left:} computational mesh of \textbf{SB} with \textbf{M=114354} cells and dofs. \textbf{Bottom right:} computational mesh of \textbf{HB} with \textbf{M=1335256} cells: only the boundary of the computational domain that intersects the hull is shown.}
    \label{fig:meshes}
\end{figure}

The computational meshes employed and the STL files used for the 3d object generation are shown in Figure~\ref{fig:meshes}. The STL files are considered as point clouds and will be the inputs and outputs of the cGMs: the number of points is constant for each training and generated geometry. The STL files of test cases \textbf{SB} and \textbf{HB} have \textbf{145821} and \textbf{5000} points, respectively. We remark that only the deformations of the bulb are generated by our cGMs, that is only the \textbf{5000} points of the bulb over the \textbf{33866} points representing the whole hull are deformed. The computational meshes have sizes \textbf{114354} and \textbf{1335256} for the \textbf{SB} and \textbf{HB} test cases, respectively.

We check the Jensen-Shannon Distance (JSD)~\cite{murphy2012machine} between the two quantity of interests $X$ and $Y$
\begin{equation}
    JSD(X,Y)=\sqrt{\frac{1}{2}KL(p_{X}||0.5\cdot p_{X}+0.5\cdot p_{Y})+0.5\cdot KL(p_{Y}||0.5\cdot p_{X}+0.5\cdot p_{Y})}
\end{equation}
where $p_{X}$ and $p_{Y}$ are the p.d.f. of $X$ and $Y$ respective (and are estimated from the samples using kernel density estimation) and \begin{equation}
KL(p||q)=\int_{\mathbb{R}^{n}}p(x)log\left(\frac{p(x)}{q(x)}\right)dx
\end{equation}
We choose JSD because: it is a distance on the probability space, it is bounded between $0$ and $1$ and it is invariant under affine transformation. These properties permit us to compare the model performances on different quantities.

For every architecture of each test case, we will also evaluate the sum of the variance (\textit{Var}) of each point of the generated 3d point clouds. It will be a useful metric to determine which cGM produces the richest distribution in terms of variability of the sampled 3d point clouds.

In the Appendix~\ref{sec:appendix} are reported the architectures' details and training specifics.

\subsection{Stanford Bunny (SB)}
\label{subsec:bunny}

The 3d object we employ is the Stanford bunny~\cite{turk}. Our objective is to preserve the barycenter $\mathbf{x}_B\in\Omega_{\text{bunny}}\subset\mathbb{R}^3$, so we have the following set of constraints:
\begin{equation}
    c_x=\frac{1}{n}\sum_{i=1}^M x_{i},\qquad
    c_y=\frac{1}{n}\sum_{i=1}^M y_{i},\qquad
    c_z=\frac{1}{n}\sum_{i=1}^M z_{i}.
\end{equation}

The numerical model we are going to use to validate the results of cGMs is a mixed Poisson problem:
\begin{subequations}
    \begin{align}
        \Delta u(\mathbf{x}) &= f(\mathbf{x}),\qquad  \mathbf{x}\in\Omega_{\text{bunny}},\\
        u(\mathbf{x})&=0,\qquad \mathbf{x}\in\partial\Omega_{\text{bunny}}\cap\{y=0\}=\Gamma_{D},\\
        \mathbf{n}\cdot\nabla u(\mathbf{x})&=0,\qquad \mathbf{x}\in\partial\Omega_{\text{bunny}}\cap\{y=0\}^c=\Gamma_{N},
    \end{align}
    where the source term $f:\mathbb{R}^3\rightarrow\mathbb{R}$ is :
    \begin{equation}
        f(\mathbf{x})=\begin{cases}
            f(\mathbf{x})=e^{\frac{1}{100-\lVert \mathbf{x}-\mathbf{x_B}\rVert_2^2}}, \qquad \lVert \mathbf{x}-\mathbf{x_B}\rVert_2<10\\
            0,\qquad\text{otherwise}
        \end{cases}.
    \end{equation}
\end{subequations}
The geometrical properties we are going to compare are the moments of inertia with uniform density equal to $1$ ($I_{xx}$, $I_{xy}$, $I_{xz}$, $I_{yy}$, $I_{yz}$, $I_{zz}$) and the integral of the heat on the boundary with homogeneous Neumann conditions for the \textbf{SB} test case:
\begin{subequations}
    \label{eq:SB_met}
    \begin{align}
        I_{xx}=\int_{\Omega_{\text{bunny}}} r_{X}^2(\mathbf{x})d\mathbf{x},\ I_{yy}=\int_{\Omega_{\text{bunny}}} r_{Y}^2(\mathbf{x})d\mathbf{x},\ I_{zz}=\int_{\Omega_{\text{bunny}}} r_{Z}^2(\mathbf{x})d\mathbf{x},\qquad\text{(principal moments of inertia)}\label{eq:pmomentsBunny}\\
        I_{xy}=\int_{\Omega_{\text{bunny}}} r_{X}(\mathbf{x})r_{Y}(\mathbf{x})d\mathbf{x},\ I_{xz}=\int_{\Omega_{\text{bunny}}} r_{X}(\mathbf{x})r_{Z}(\mathbf{x})d\mathbf{x},\ I_{yz}=\int_{\Omega_{\text{bunny}}} r_{Y}(\mathbf{x})r_{Z}(\mathbf{x})d\mathbf{x},\qquad\text{(moments of inertia)}\label{eq:momentsBunny}\\
        I_{u} = \int_{\partial \Gamma_N} u(\mathbf{x})\ d\mathbf{x}\qquad\text{(integral over the boundary }\partial\Gamma_N\text{)}\label{eq:intergralBunny}
    \end{align}
\end{subequations}
where $r_{X}, r_{Y},r_{Z}:\Omega_{\text{bunny}}\rightarrow\mathbb{R}_{+}$ are the distances from the $x$-, $y$- and $z$-axes, respectively. These quantities are evaluated on the discrete STL point cloud for the moments of inertia $I_{xx},I_{yy},I_{zz},I_{xy},I_{xz},I_{yz}$ and on the computational mesh for the integral of the solution on the Neumann boundary $I_u$.

\begin{figure}[htpb!]
    \centering
    \includegraphics[width=0.195\textwidth, trim={0 0 700 0}, clip]{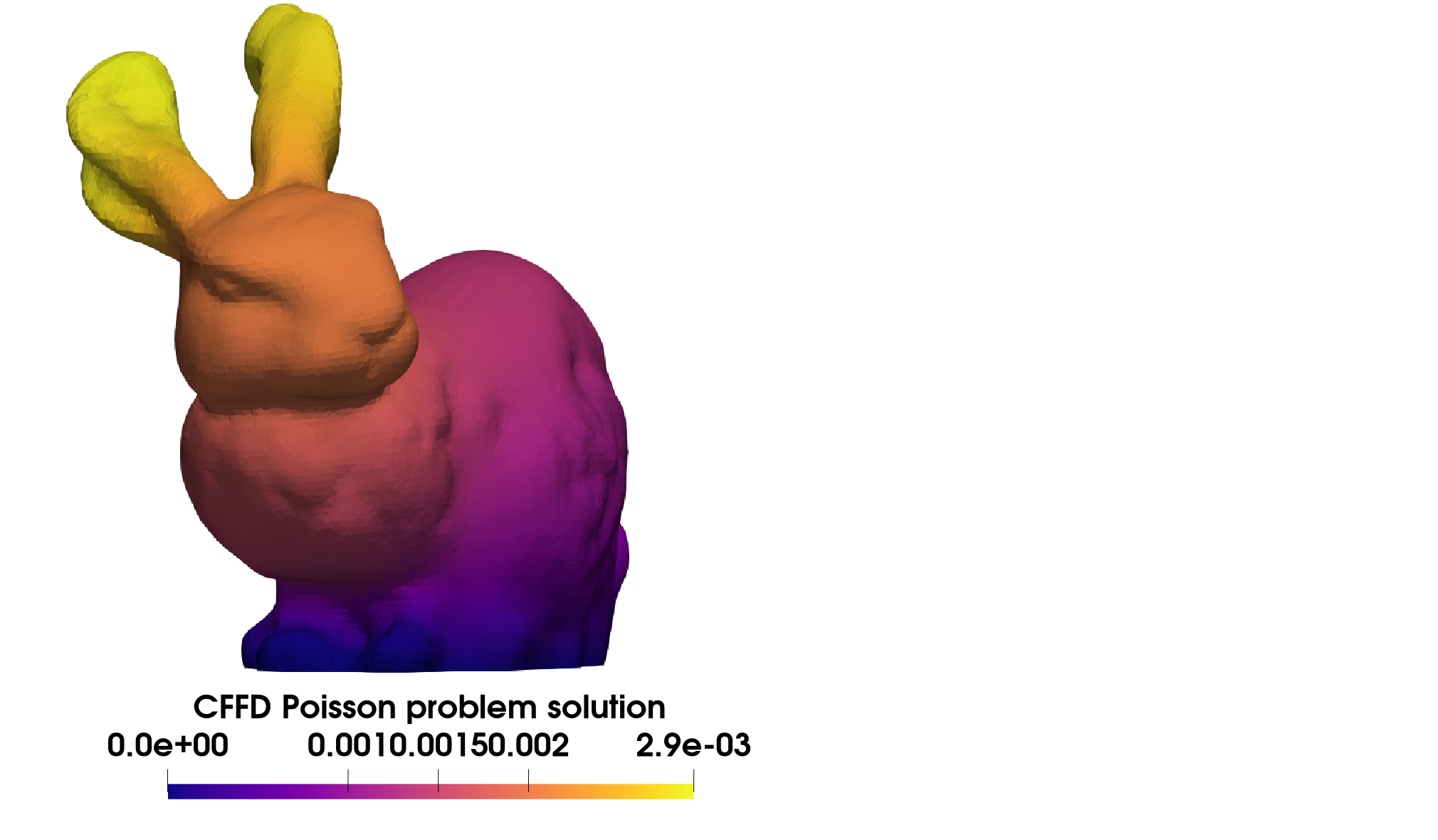}
    \includegraphics[width=0.195\textwidth, trim={0 0 700 0}, clip]{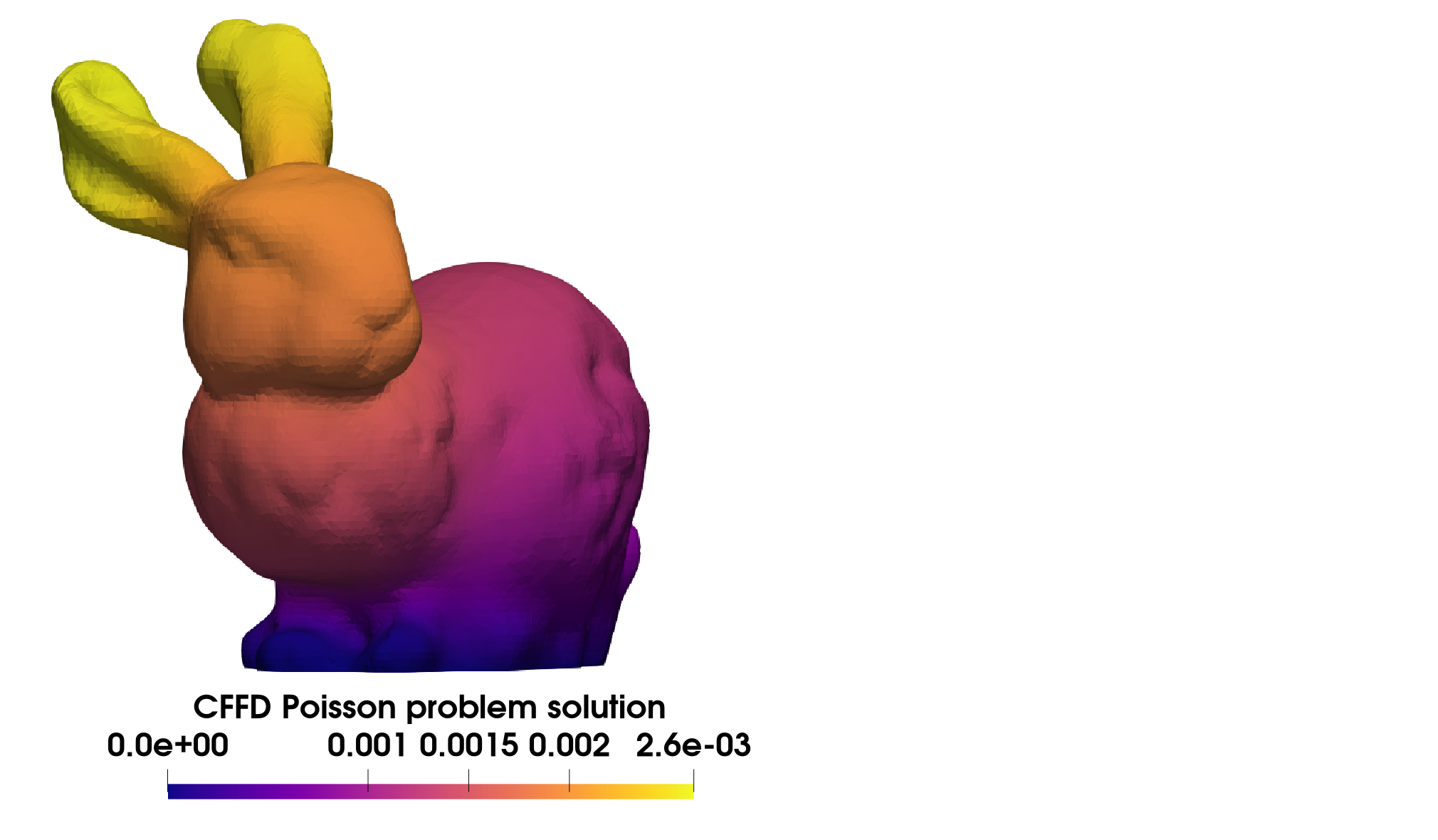}
    \includegraphics[width=0.195\textwidth, trim={0 0 700 0}, clip]{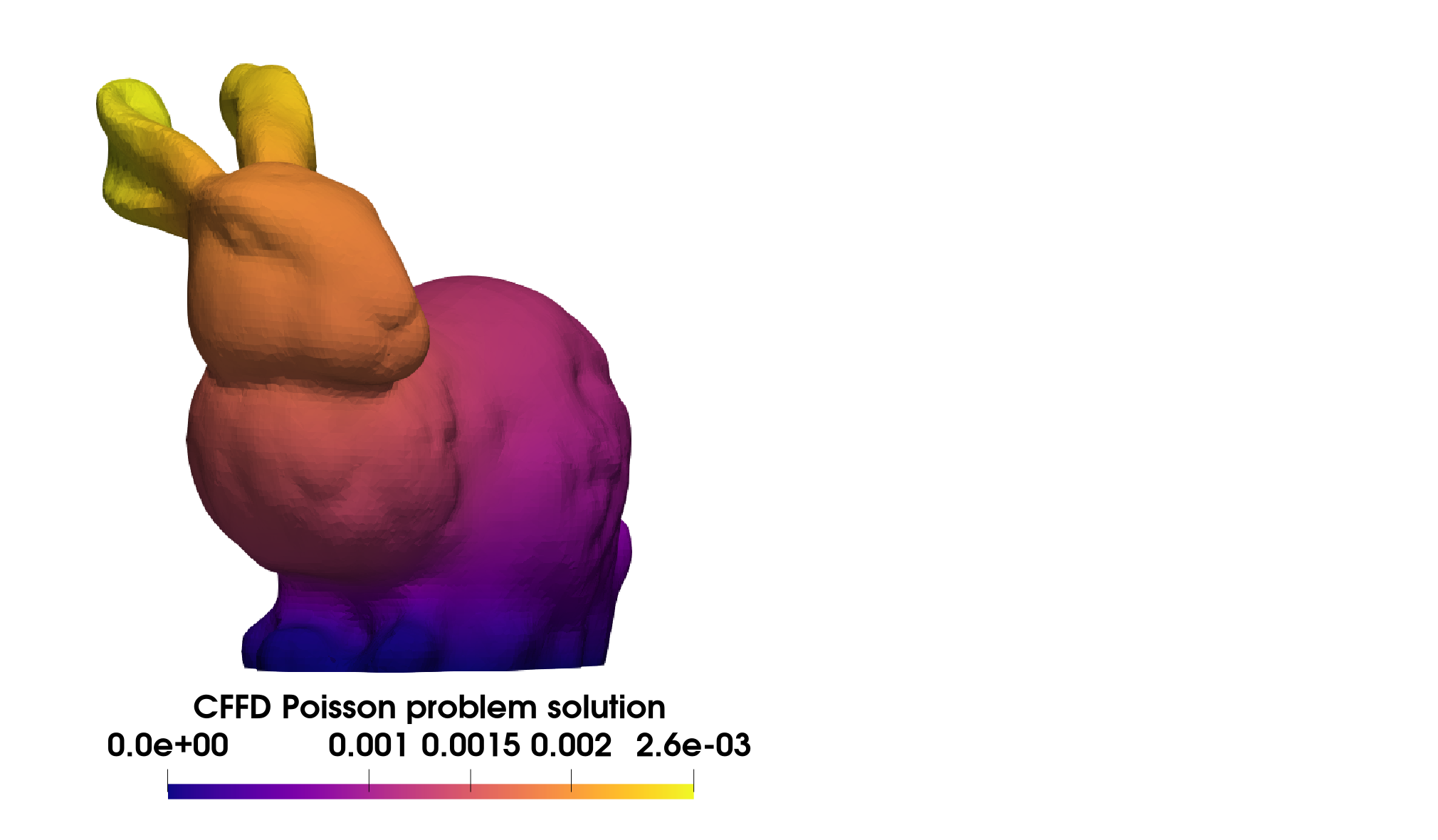}
    \includegraphics[width=0.195\textwidth, trim={0 0 700 0}, clip]{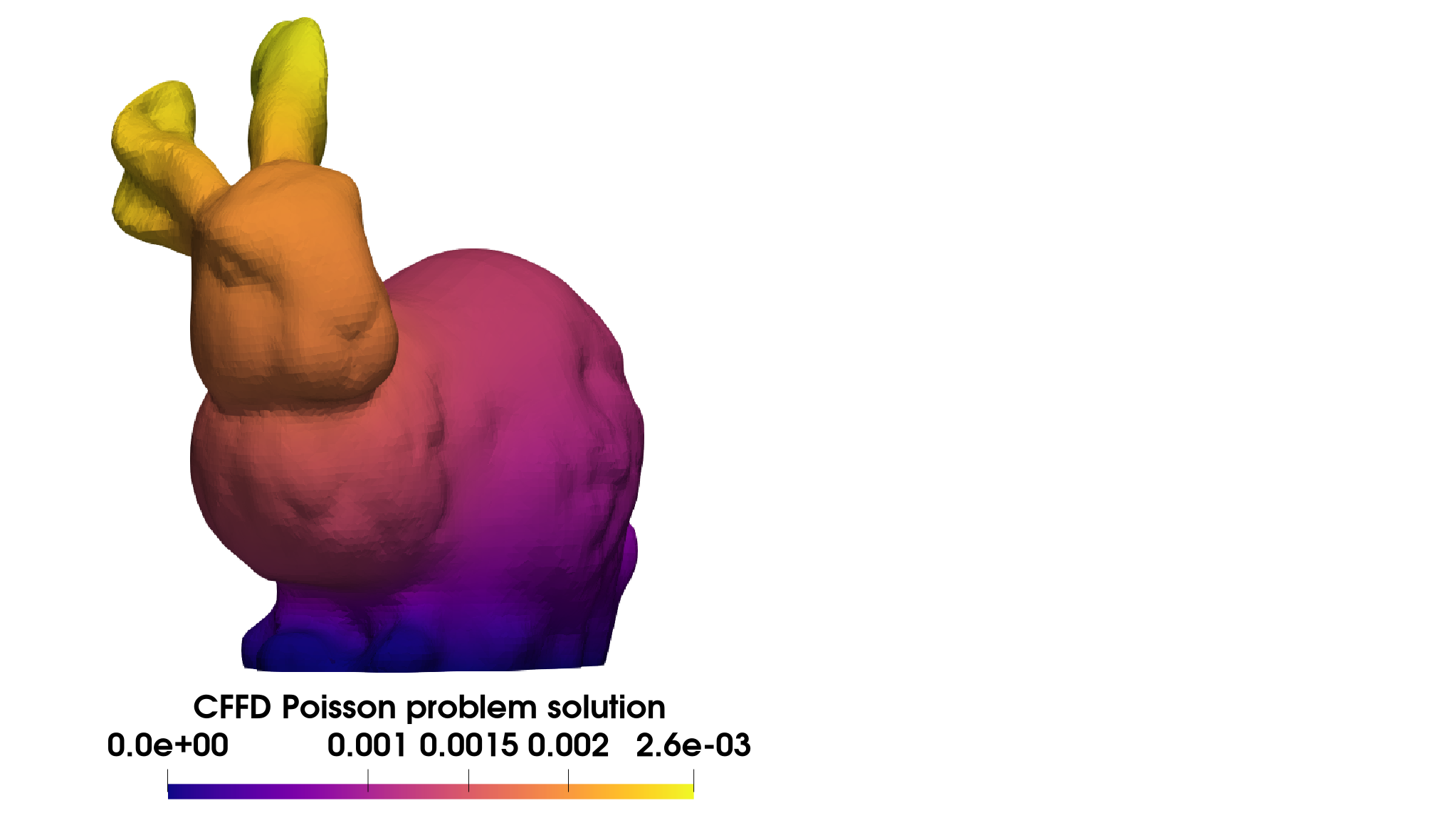}
    \includegraphics[width=0.195\textwidth, trim={0 0 700 0}, clip]{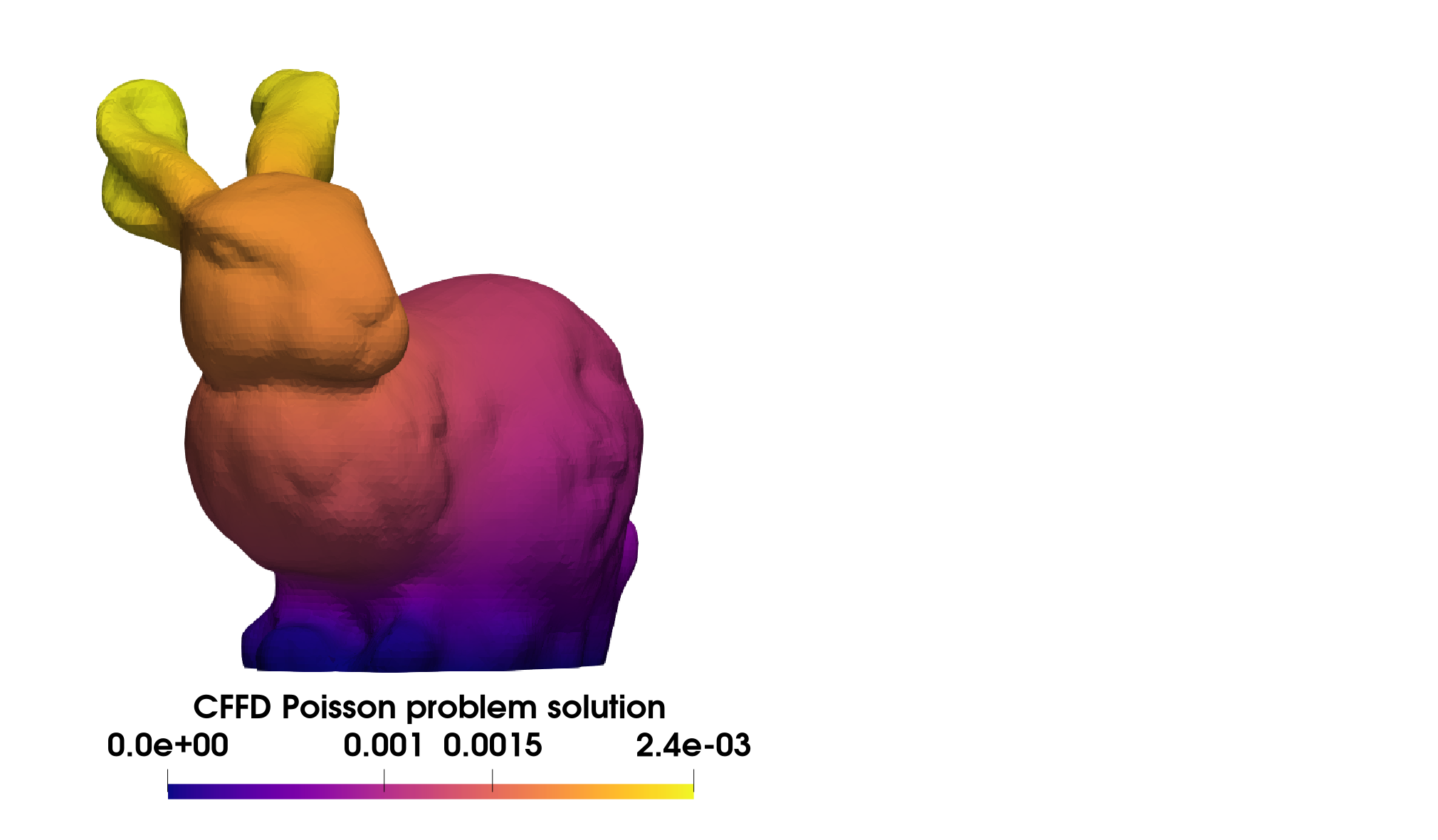}\\
    \includegraphics[width=0.195\textwidth, trim={0 0 700 0}, clip]{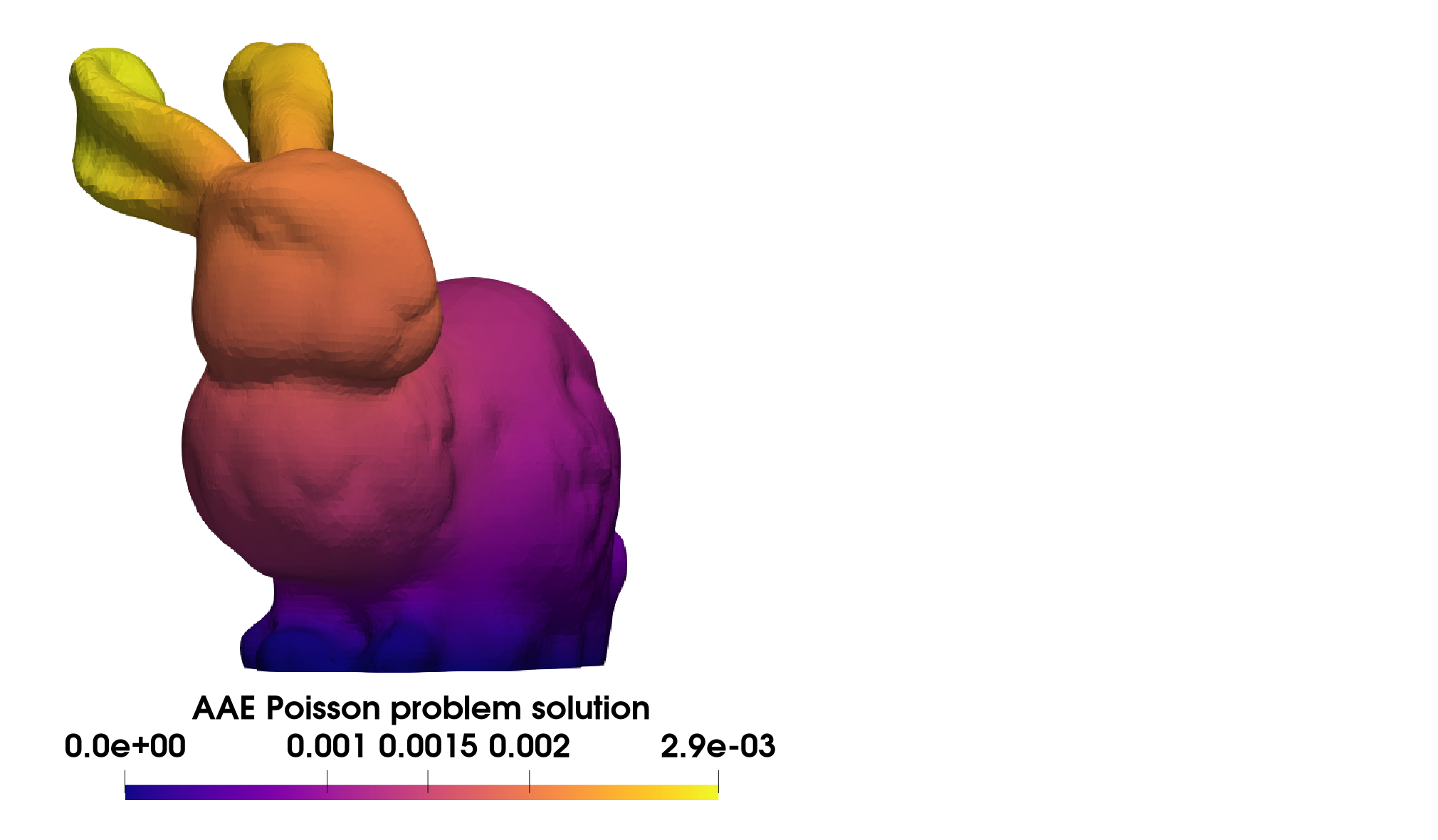}
    \includegraphics[width=0.195\textwidth, trim={0 0 700 0}, clip]{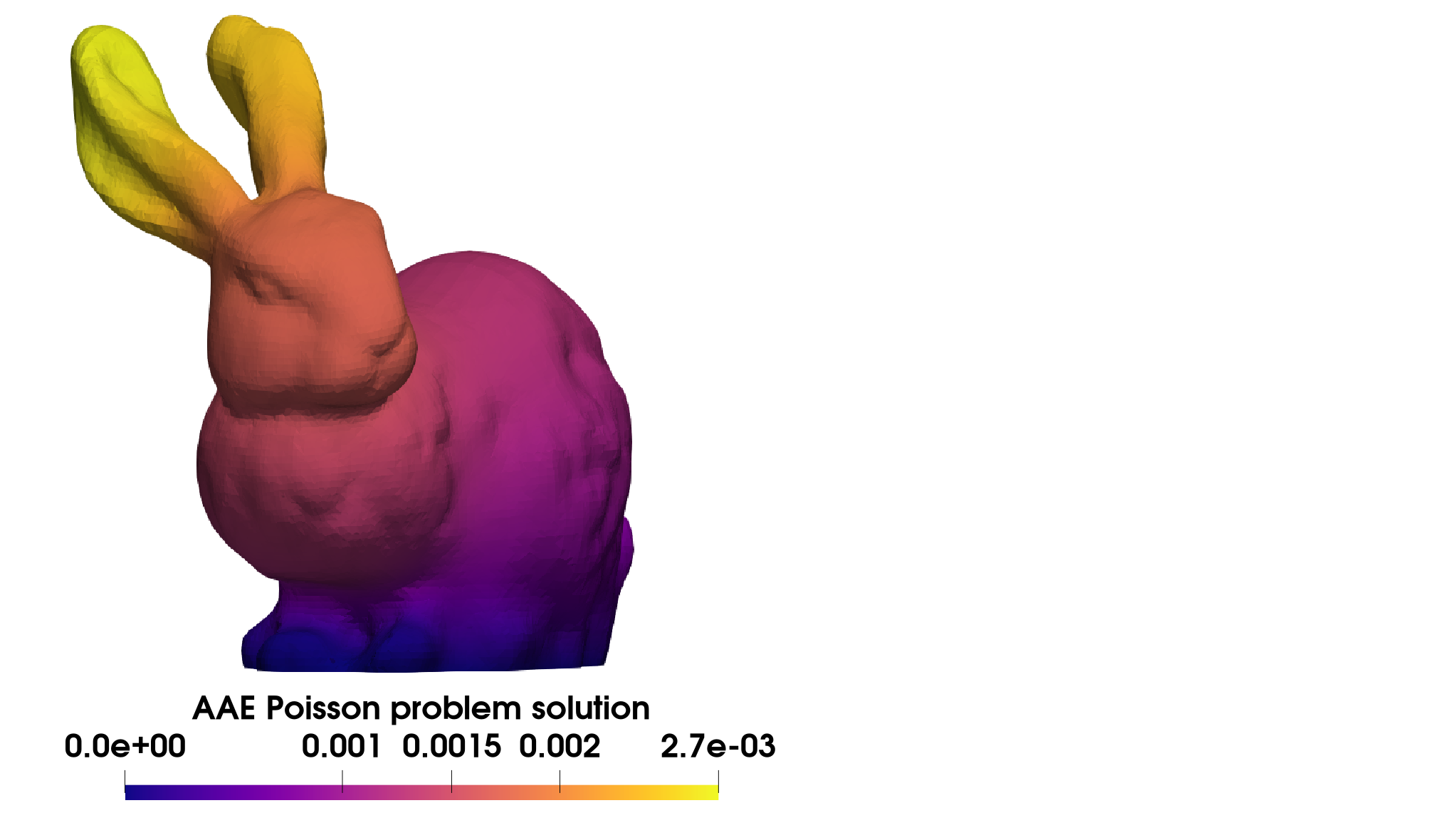}
    \includegraphics[width=0.195\textwidth, trim={0 0 700 0}, clip]{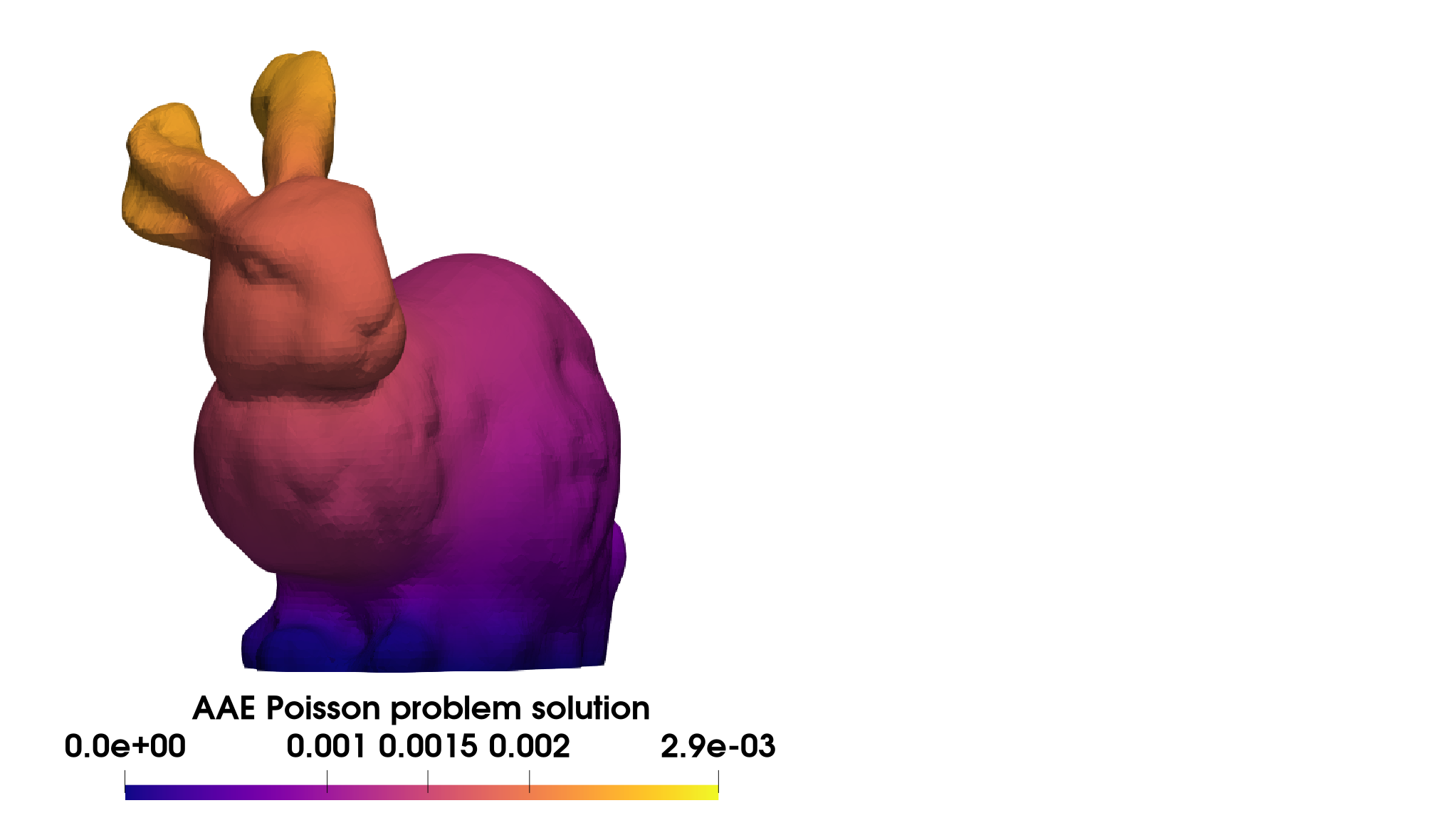}
    \includegraphics[width=0.195\textwidth, trim={0 0 700 0}, clip]{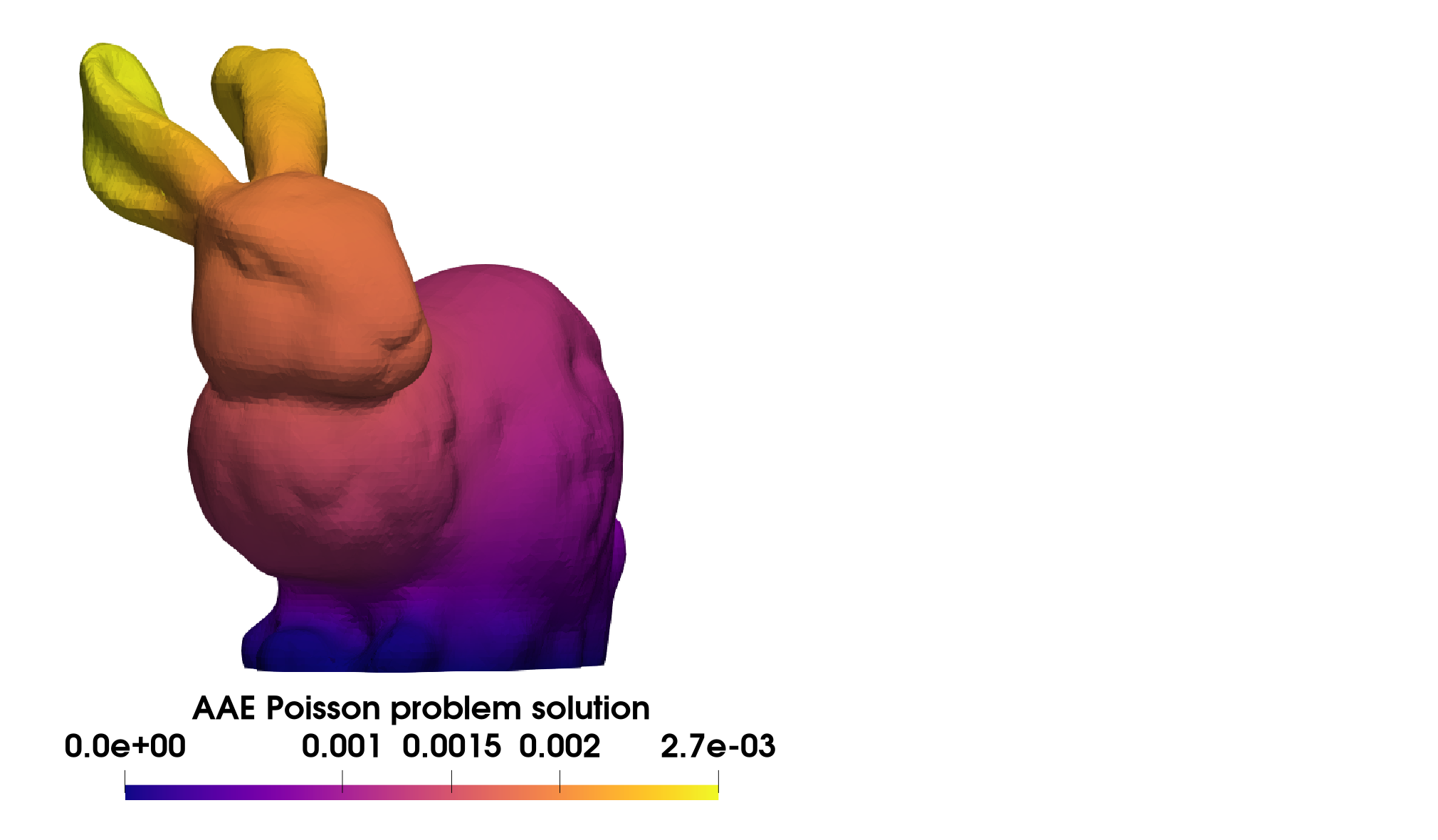}
    \includegraphics[width=0.195\textwidth, trim={0 0 700 0}, clip]{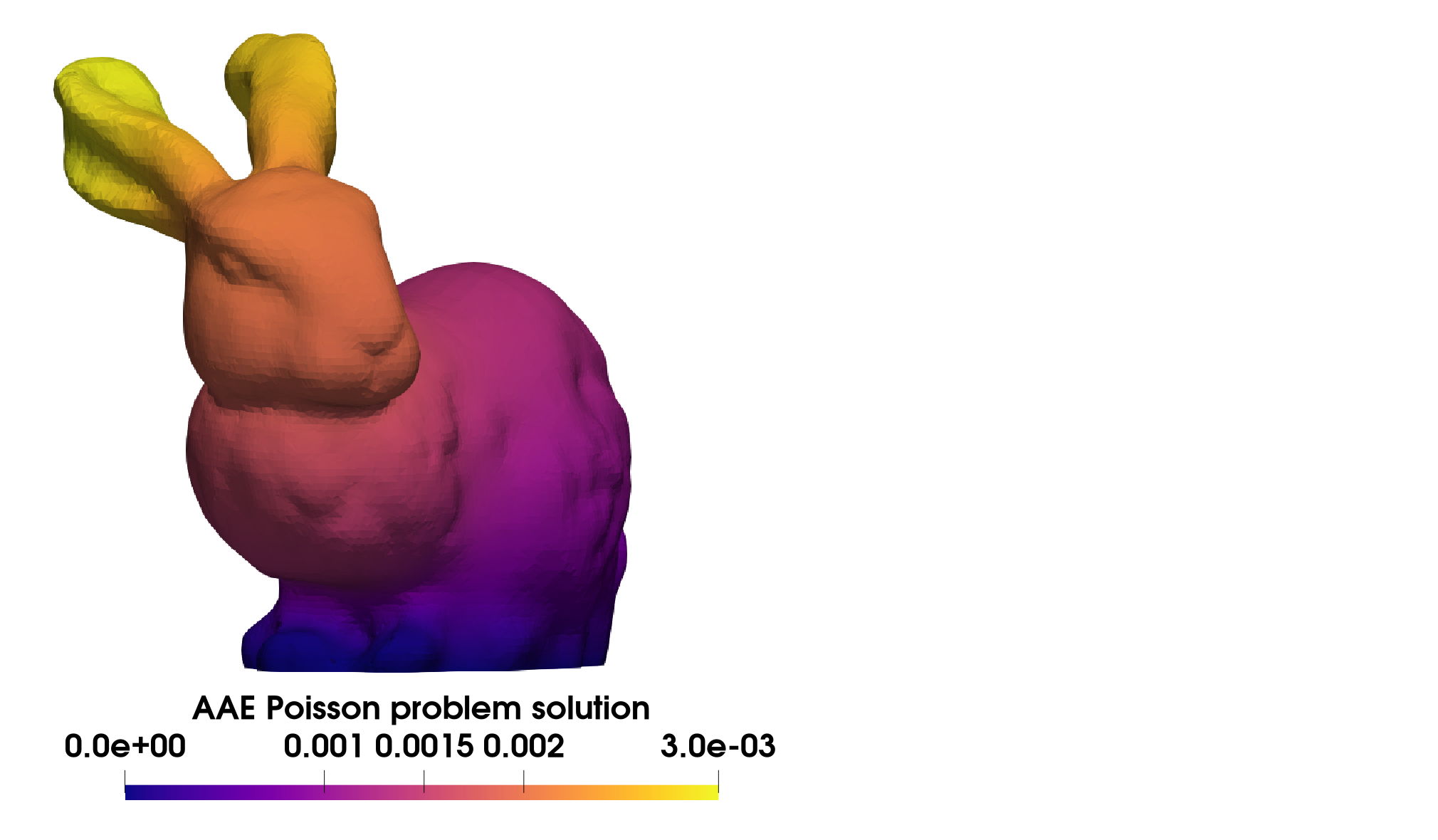}\\
    \includegraphics[width=0.195\textwidth, trim={0 0 700 0}, clip]{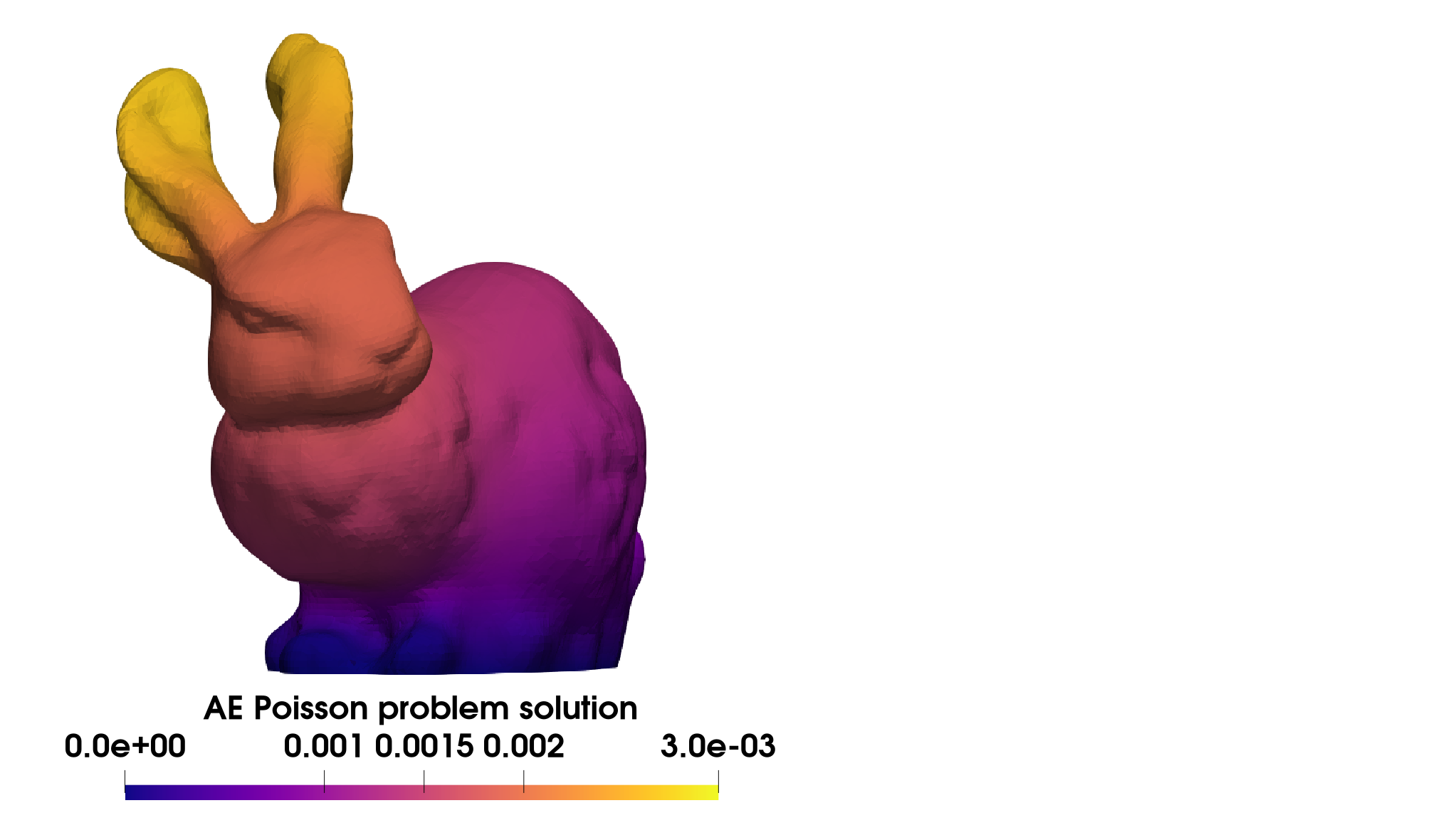}
    \includegraphics[width=0.195\textwidth, trim={0 0 700 0}, clip]{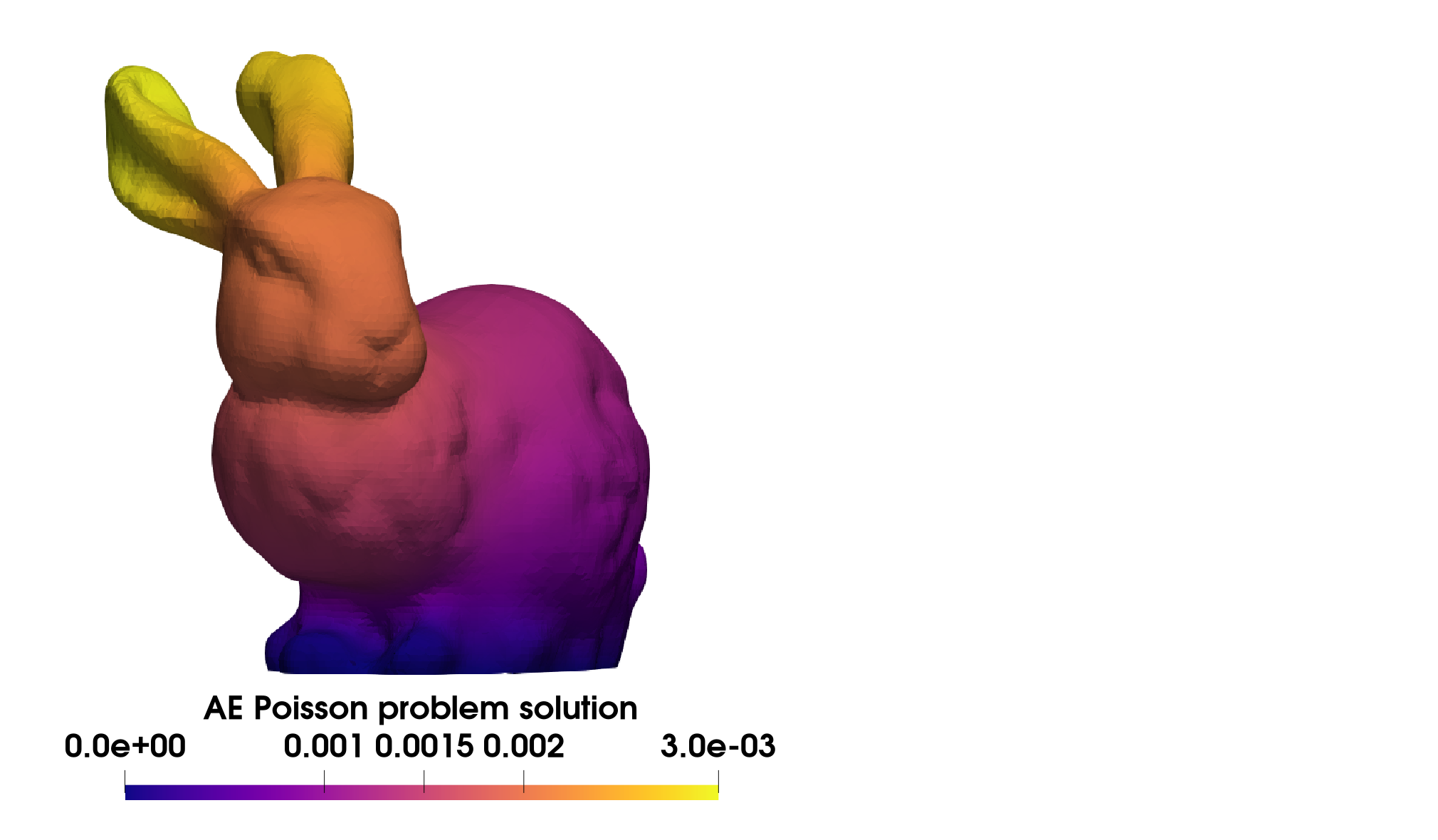}
    \includegraphics[width=0.195\textwidth, trim={0 0 700 0}, clip]{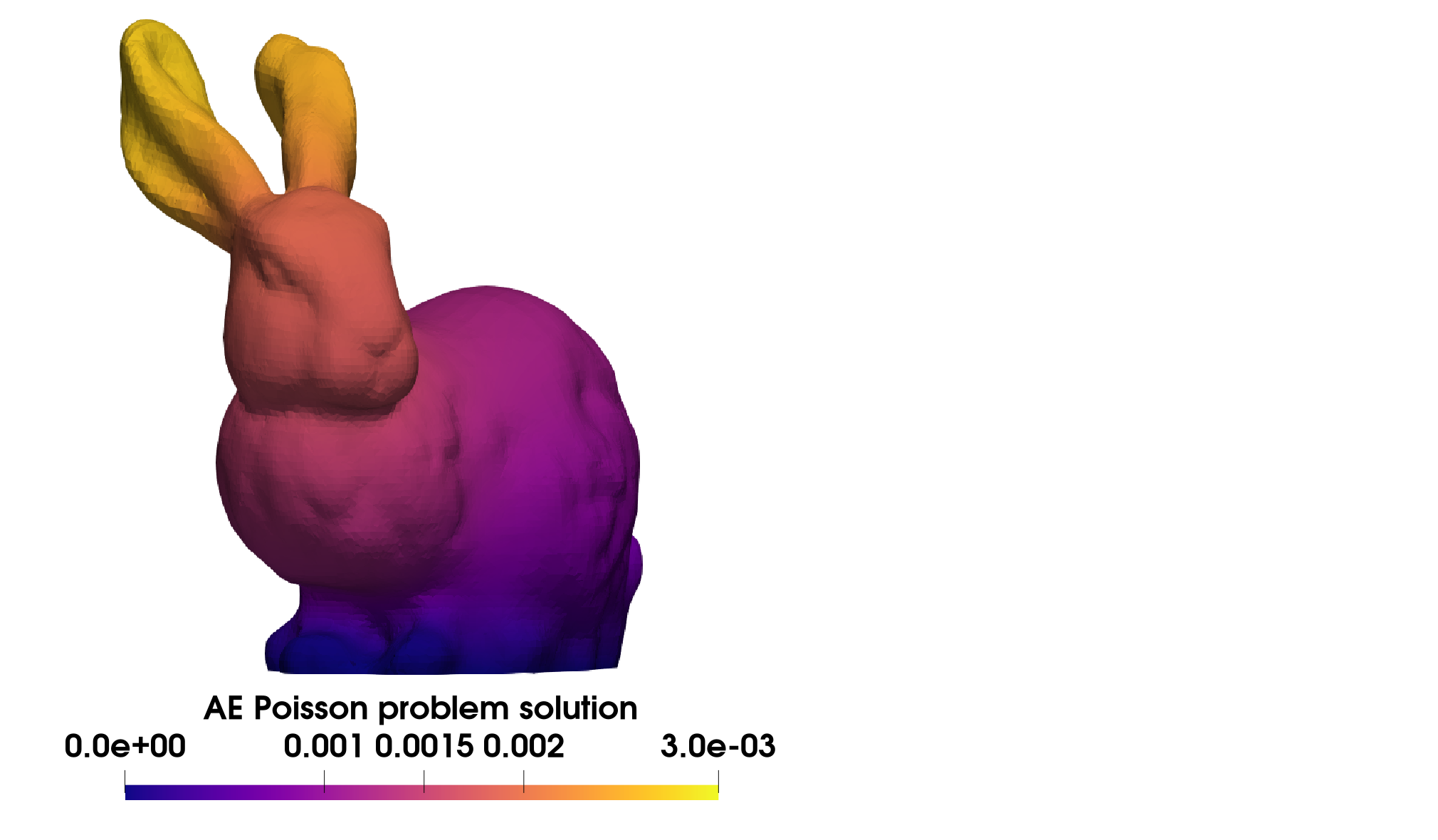}
    \includegraphics[width=0.195\textwidth, trim={0 0 700 0}, clip]{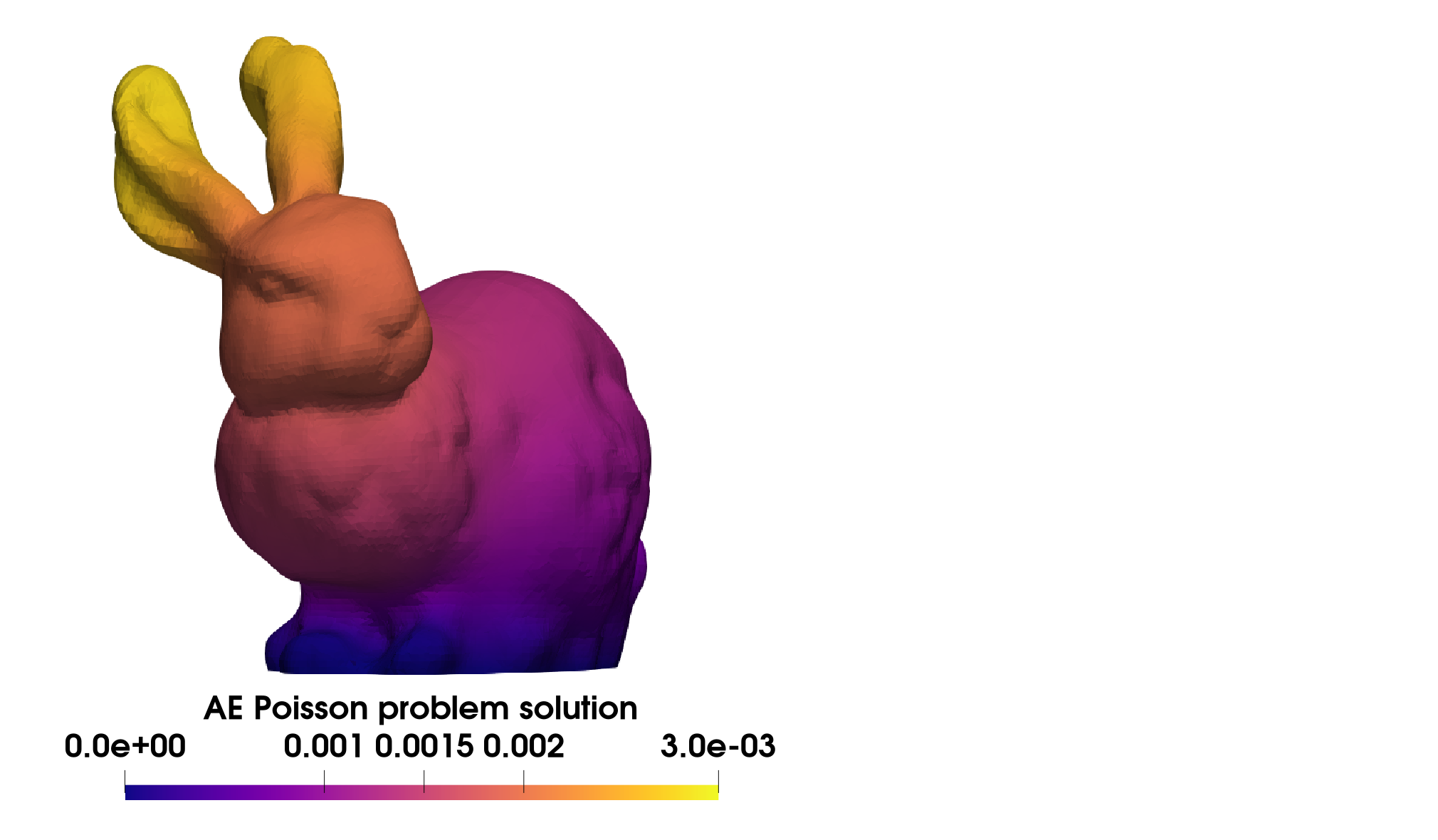}
    \includegraphics[width=0.195\textwidth, trim={0 0 700 0}, clip]{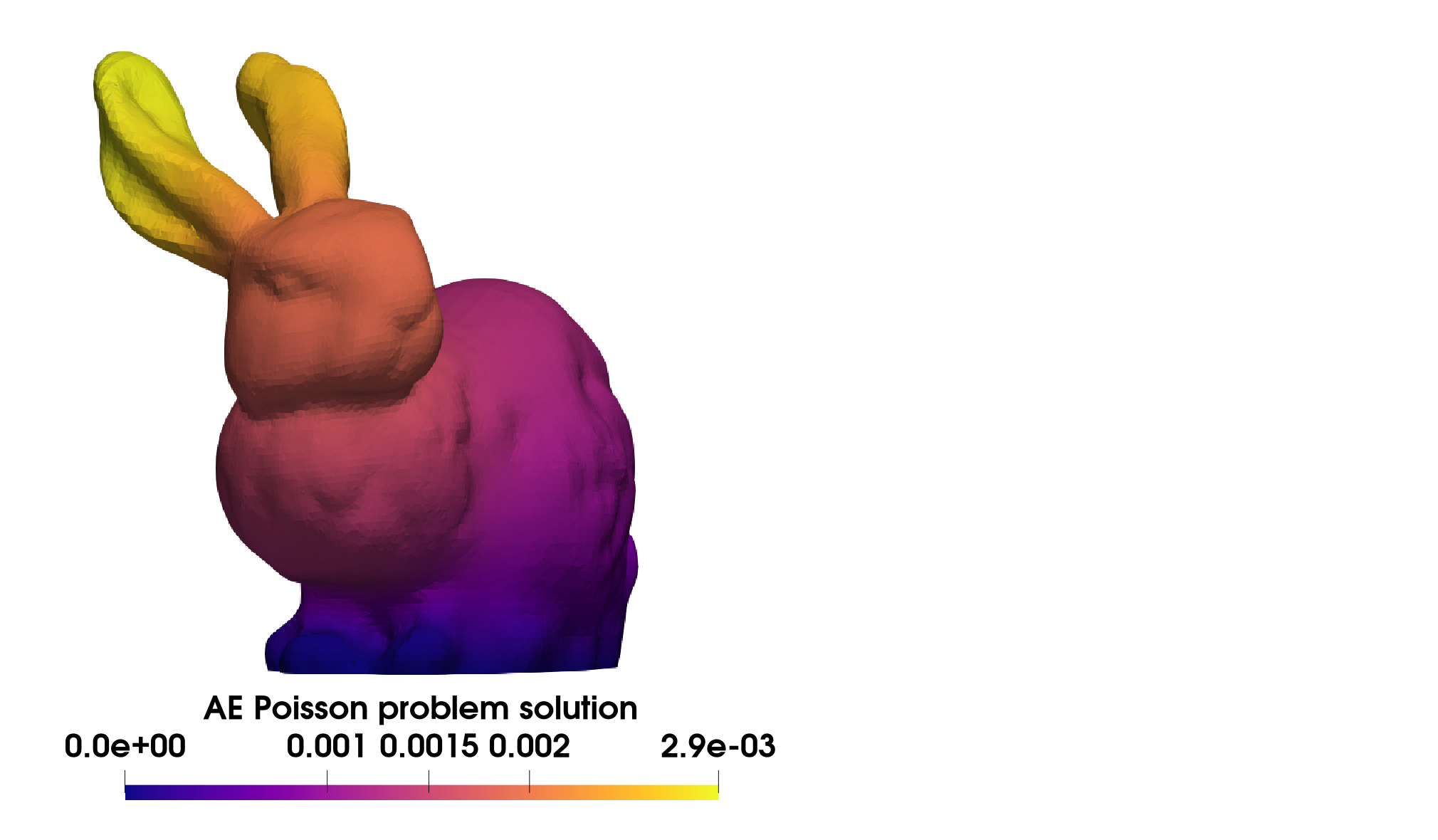}\\
    \includegraphics[width=0.195\textwidth, trim={0 0 700 0}, clip]{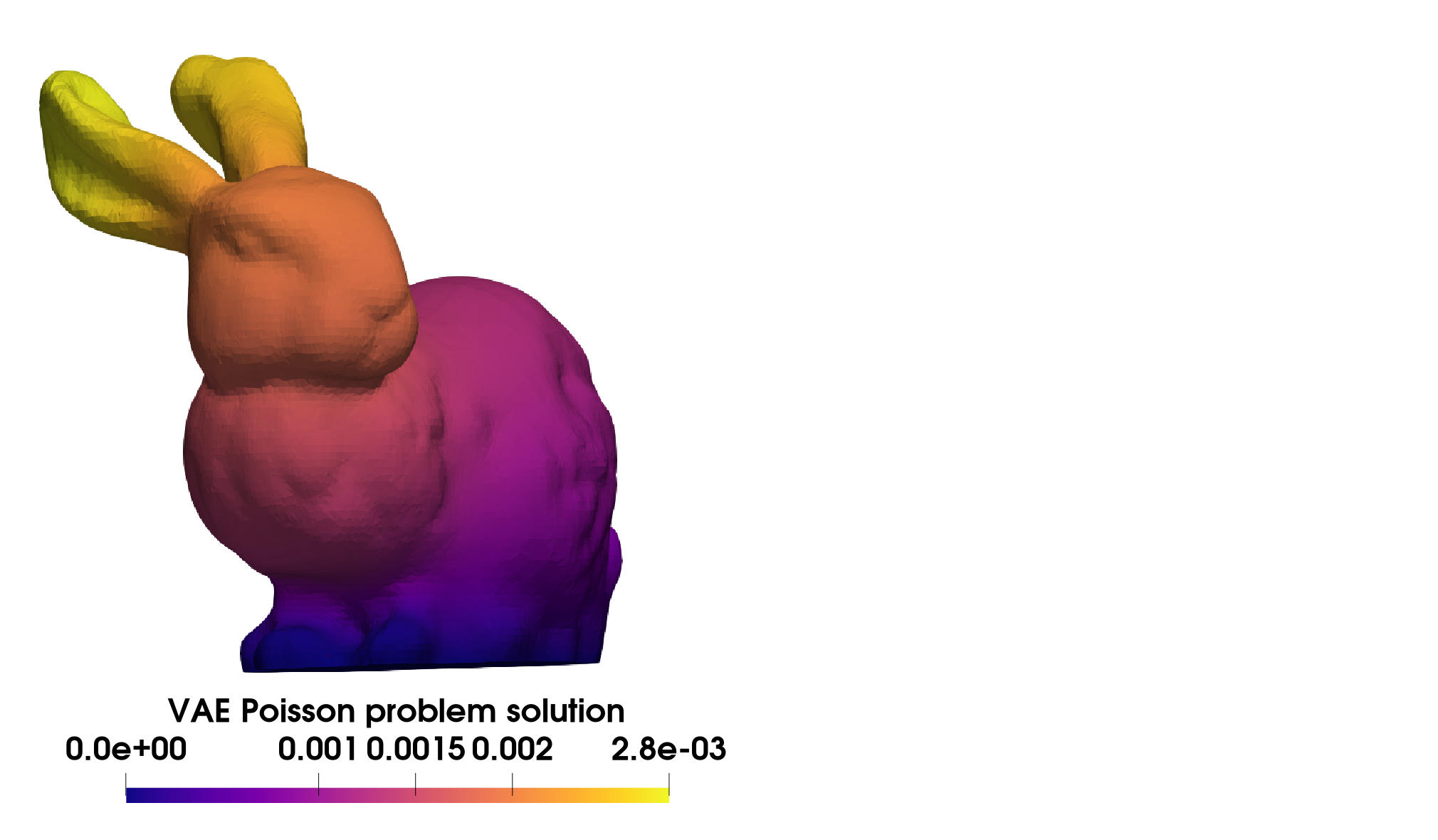}
    \includegraphics[width=0.195\textwidth, trim={0 0 700 0}, clip]{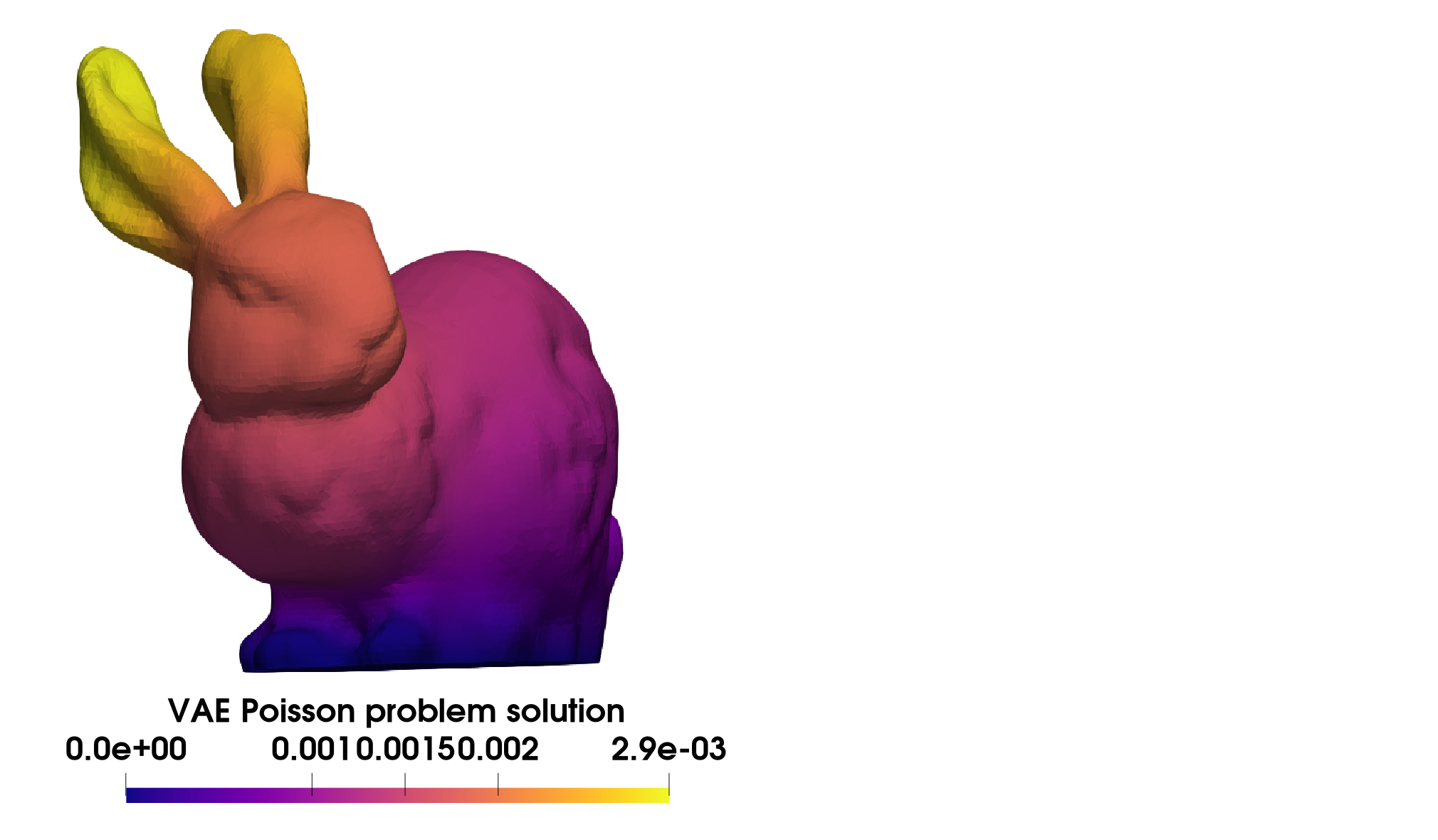}
    \includegraphics[width=0.195\textwidth, trim={0 0 700 0}, clip]{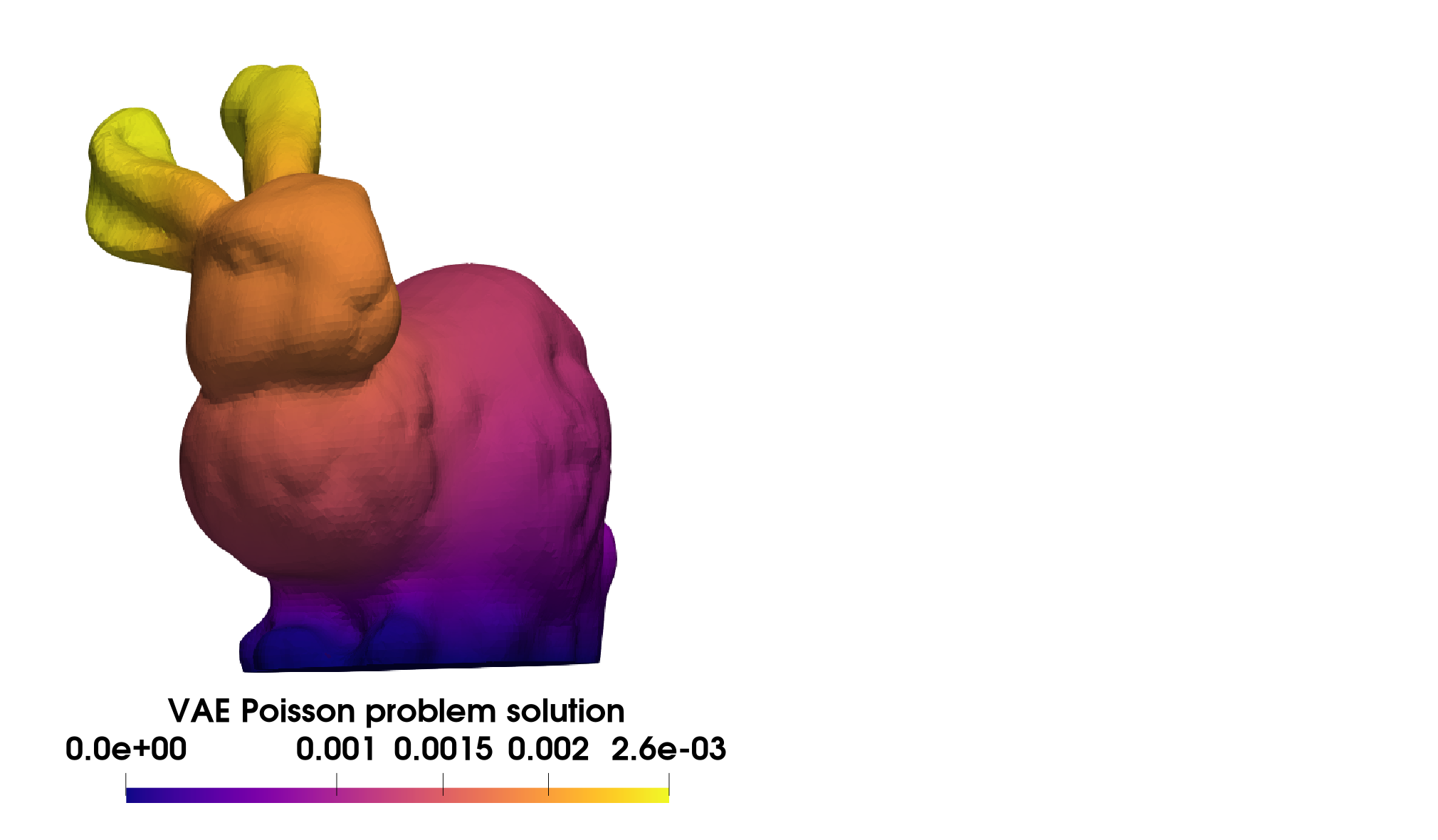}
    \includegraphics[width=0.195\textwidth, trim={0 0 700 0}, clip]{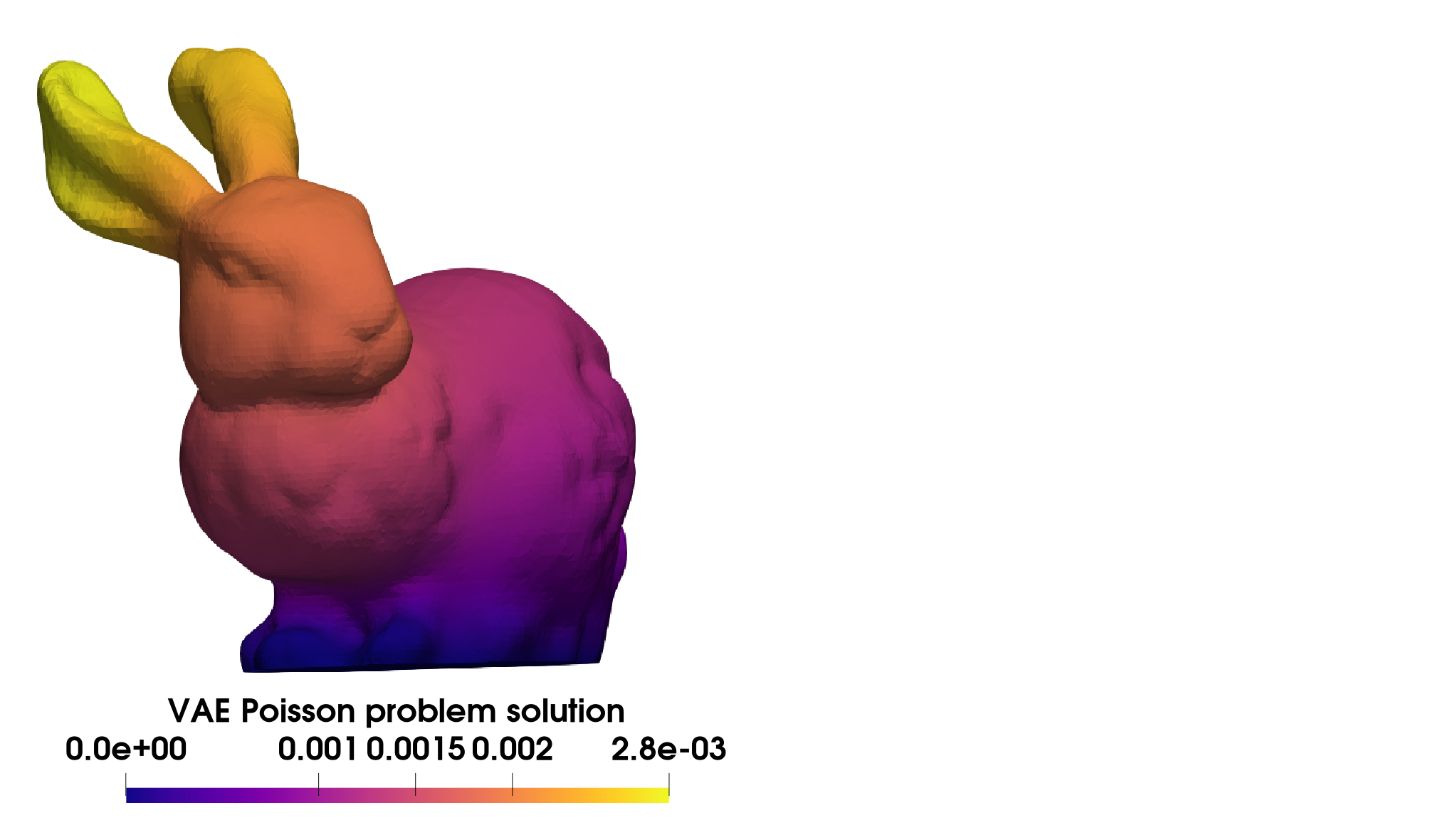}
    \includegraphics[width=0.195\textwidth, trim={0 0 700 0}, clip]{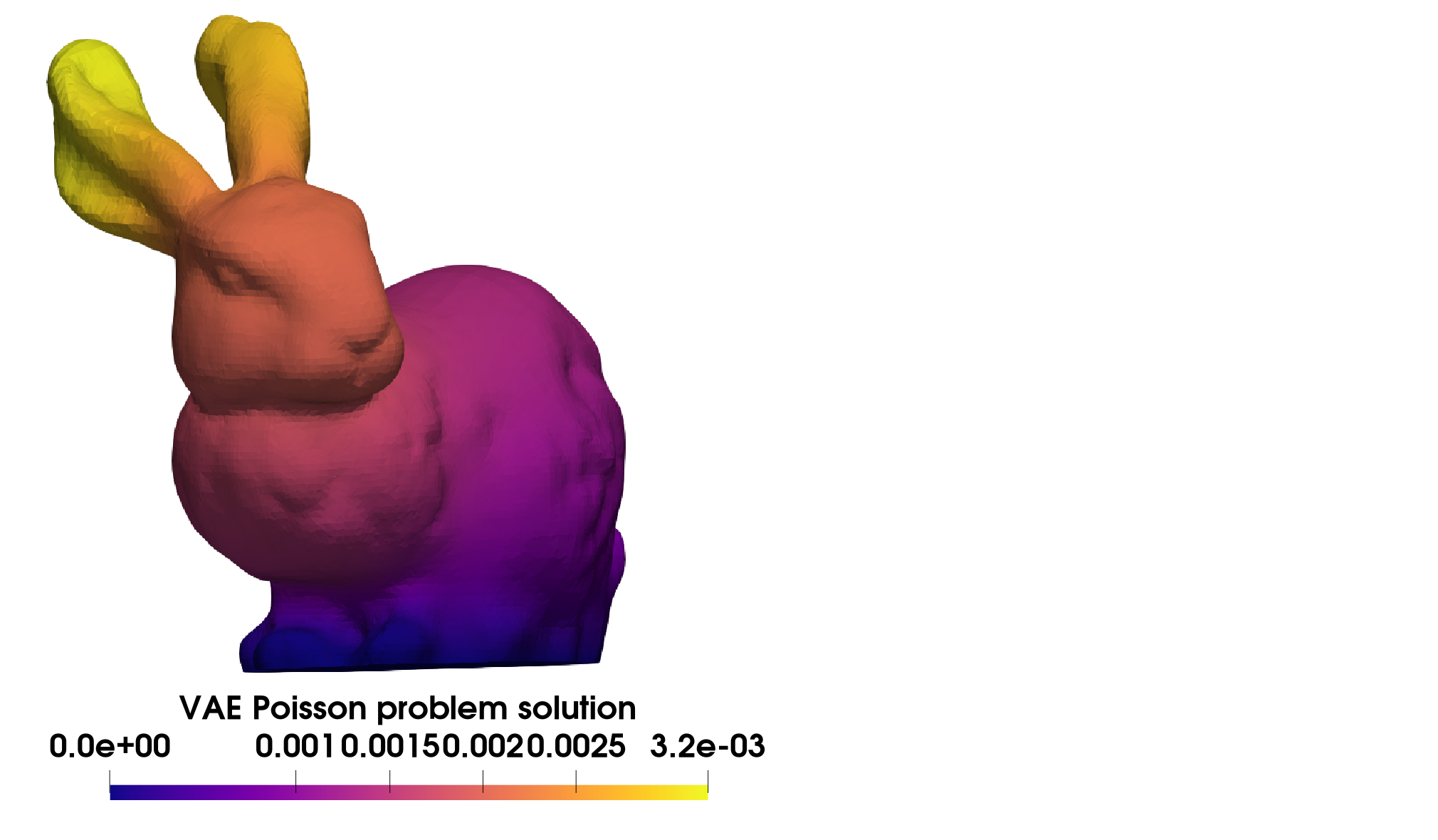}\\
    \includegraphics[width=0.195\textwidth, trim={0 0 700 0}, clip]{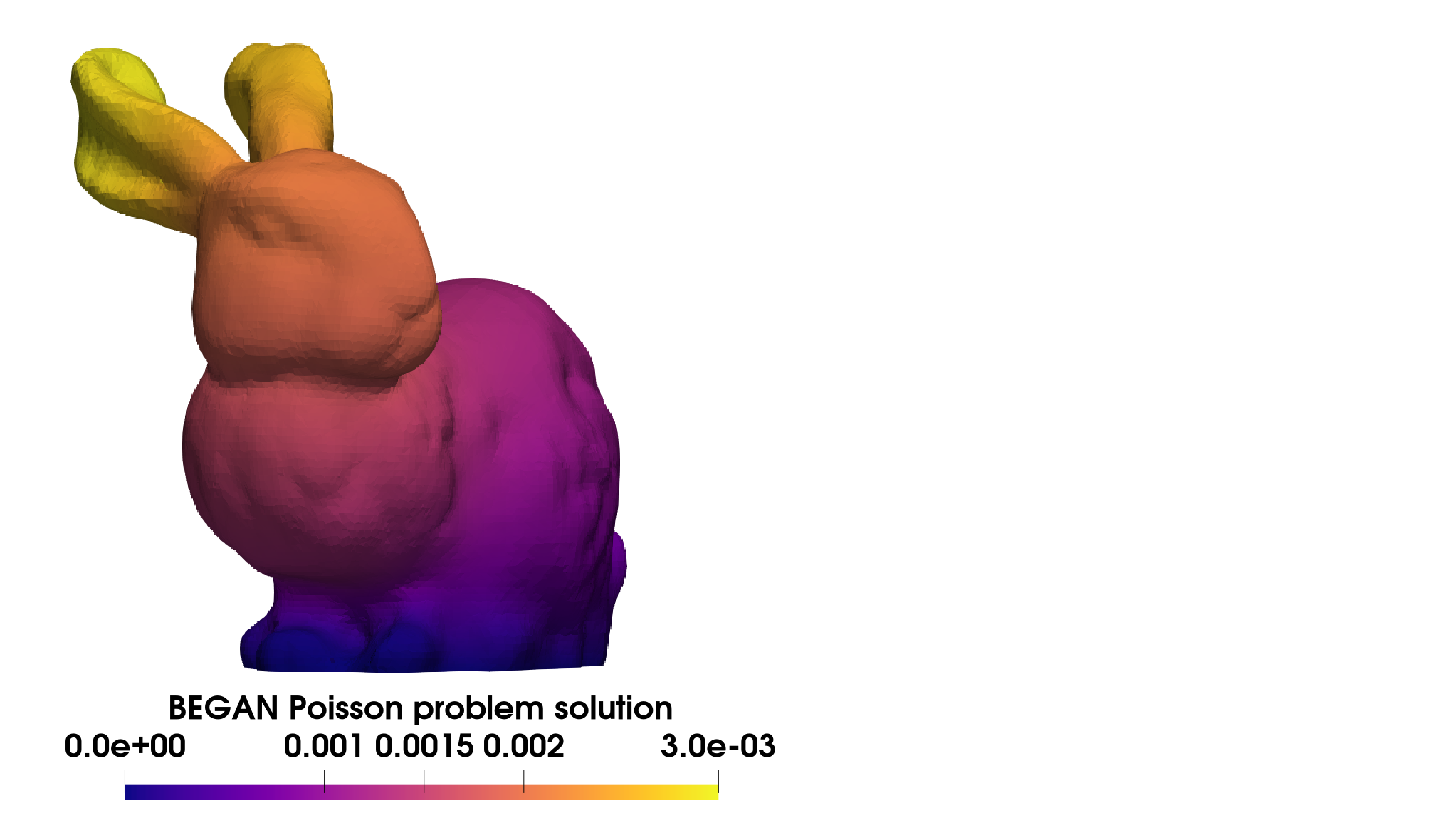}
    \includegraphics[width=0.195\textwidth, trim={0 0 700 0}, clip]{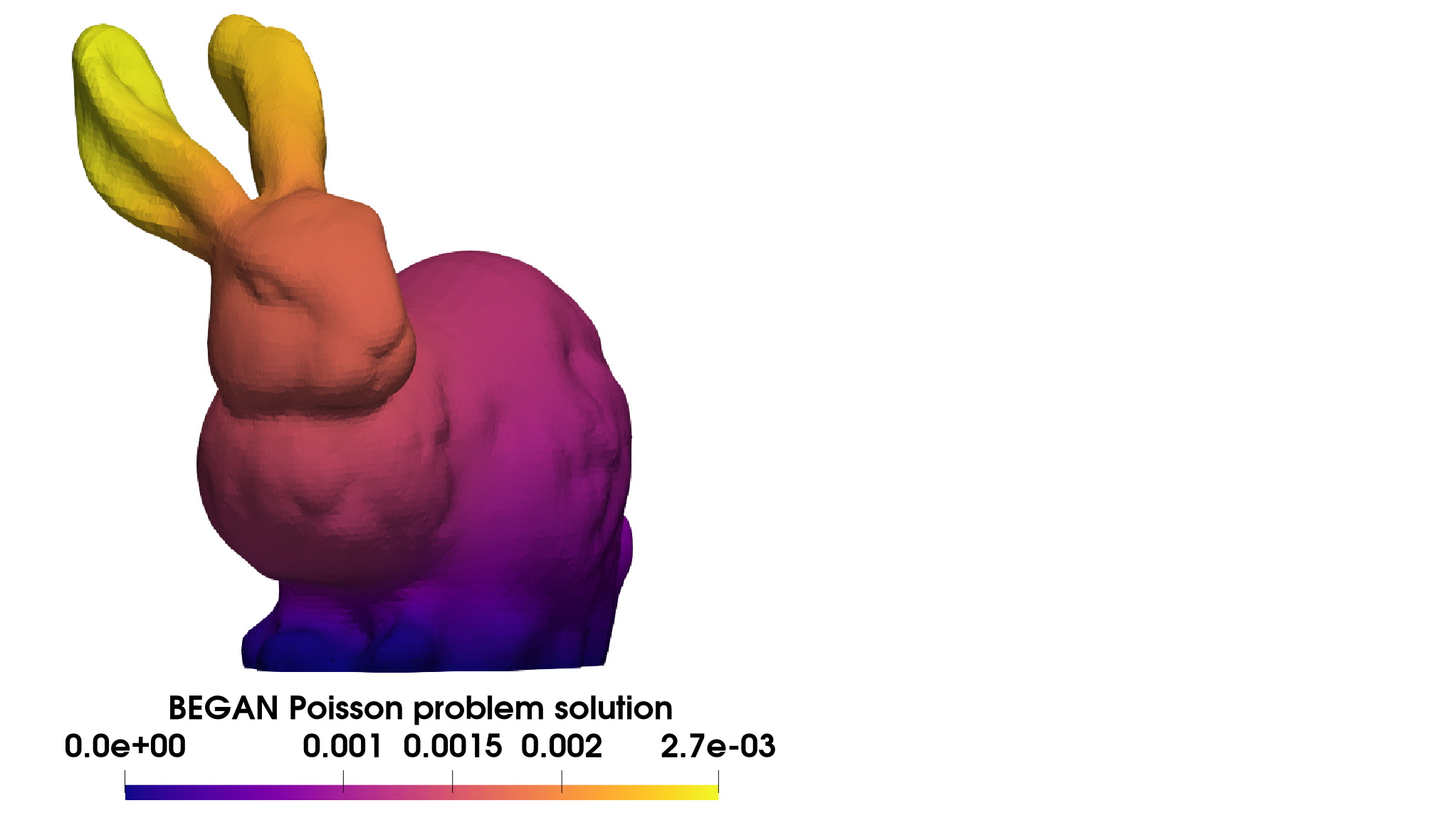}
    \includegraphics[width=0.195\textwidth, trim={0 0 700 0}, clip]{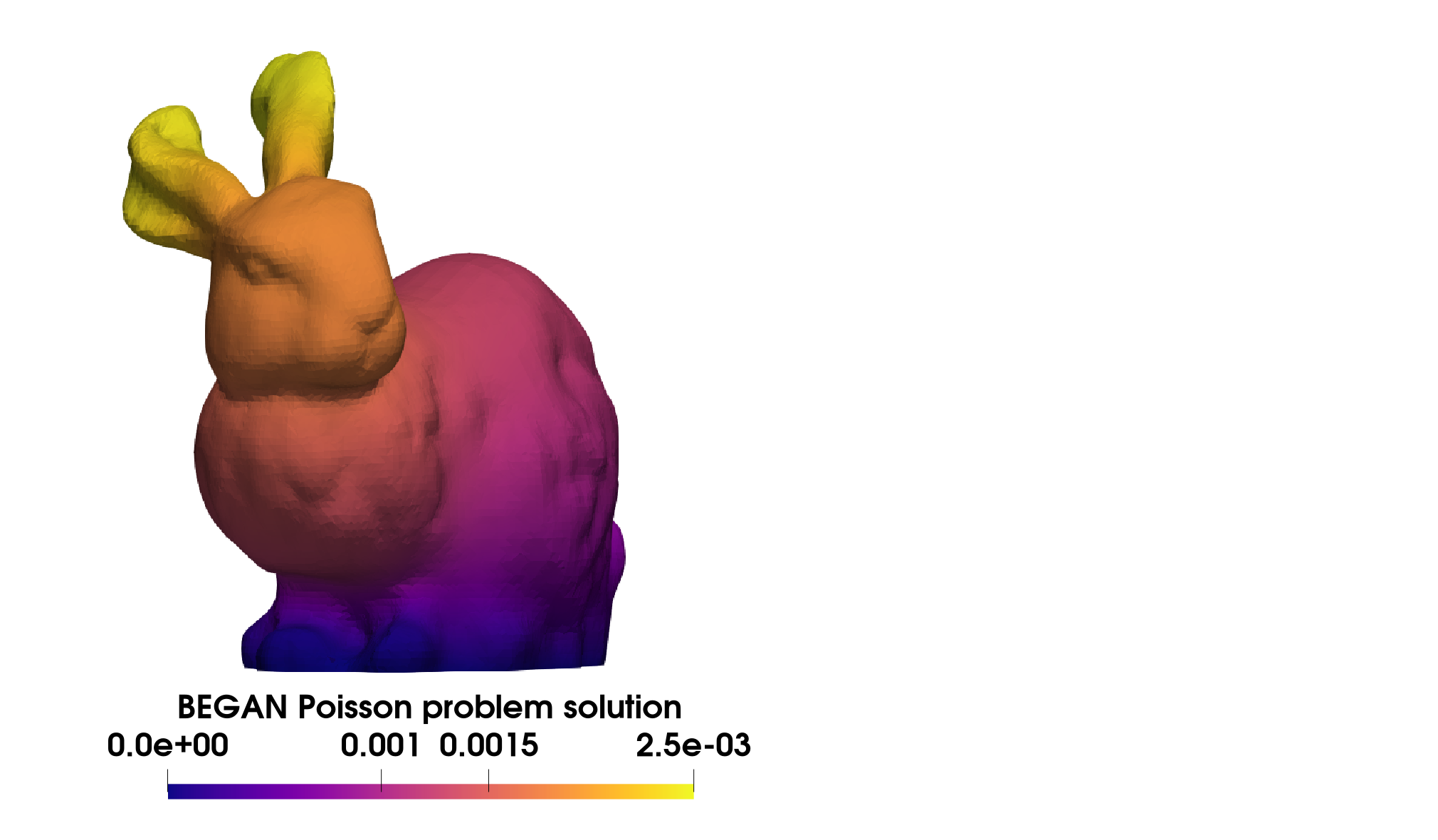}
    \includegraphics[width=0.195\textwidth, trim={0 0 700 0}, clip]{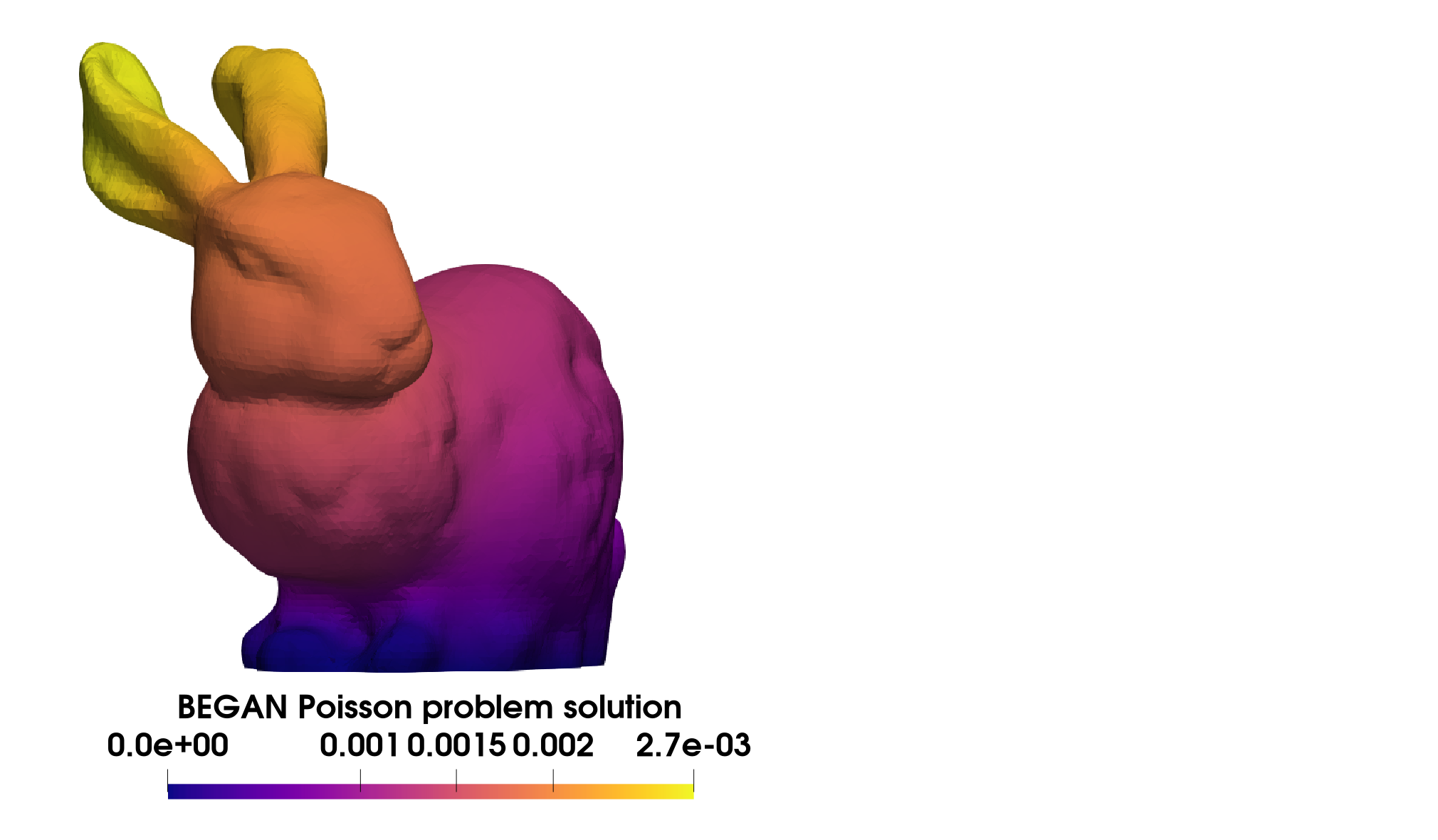}
    \includegraphics[width=0.195\textwidth, trim={0 0 700 0}, clip]{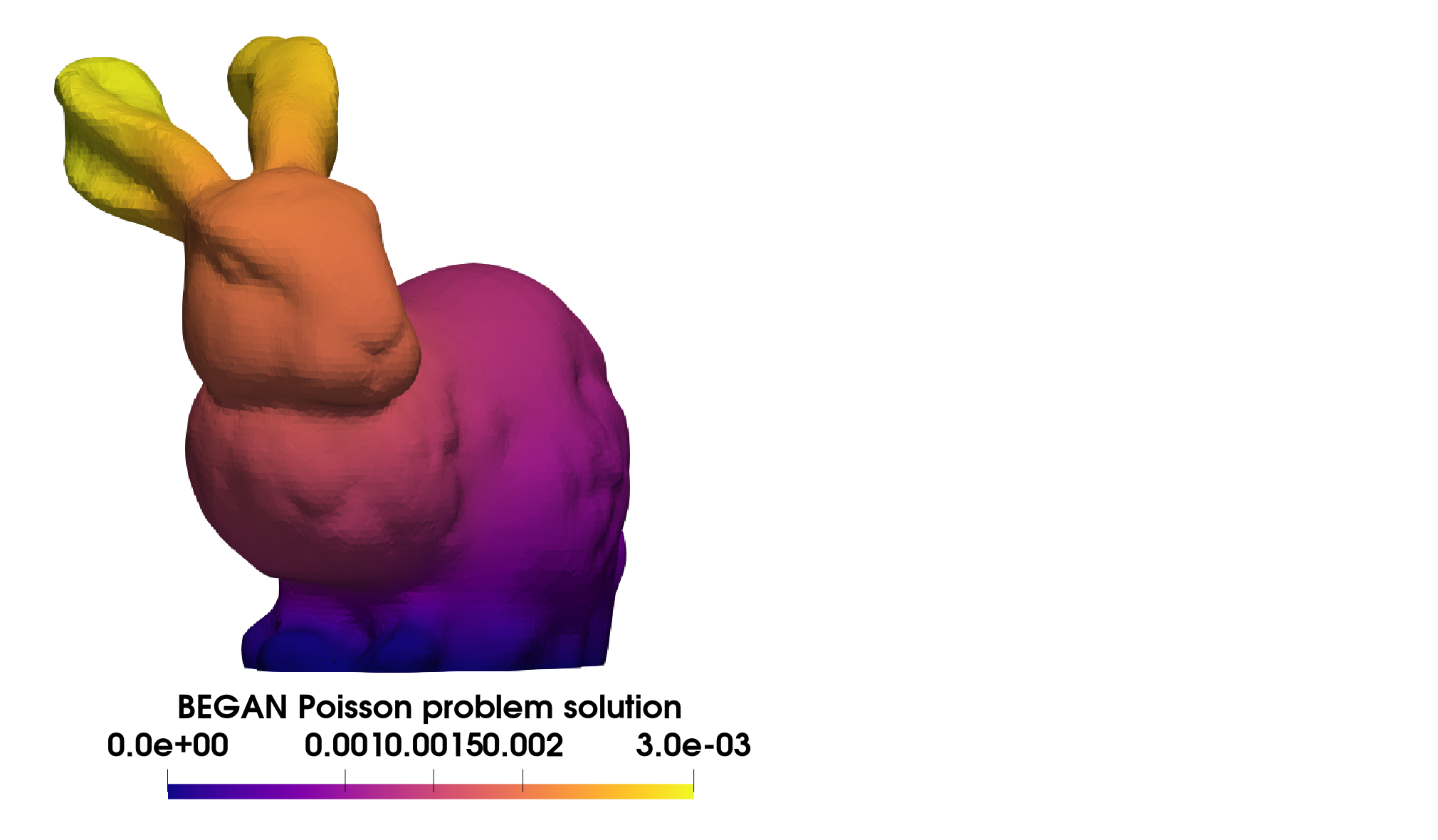}
    \caption{\textbf{SB. }Some geometrical deformations of the Stanford bunny are shown: the first row refers to constrained cFFD of section~\ref{subsec:cffd}, the second to the adversarial autoencoder of paragraph~\ref{par:aae}, the third to the simple autoencoder of paragraph~\ref{par:ae}, the fourth to the beta variational autoencoder of paragraph~\ref{par:vae} and the last to the boundary equilibrium generative adversarial networks of paragraph~\ref{par:began}. All the generative models implement the linear constraints enforcing layer of section~\ref{sec:constrained generative models} to preserve the position of the barycenter.}
    \label{fig:def_SB}
\end{figure}

The number of training and test samples are $\mathbf{n_{\text{train}}=400}$ and $\mathbf{n_{\text{test}}=200}$, respectively. For model order reduction, the number of training and test components are instead $\mathbf{n_{_{\text{ROM}},\text{train}}=80}$ and $\mathbf{n_{_{\text{ROM}},\text{test}}=20}$, respectively. We use $\mathbf{r_{\text{PCA}}=30}$ PCA modes for preprocessing inside the cGMs and $\mathbf{r_{\text{POD}}=3}$ modes to perform model order reduction. The parameters' space dimension of cFFD is $\mathbf{p=54}$, while the latent space dimensions of the cGMs is $\mathbf{R=15}$. Some deformations from the cGMs are shown in Figure~\ref{fig:def_SB}: the field shown is the solution to the mixed Poisson problem on different deformed geometries.

\begin{figure}[htpb!]
    \centering
    \includegraphics[width=0.49\textwidth]{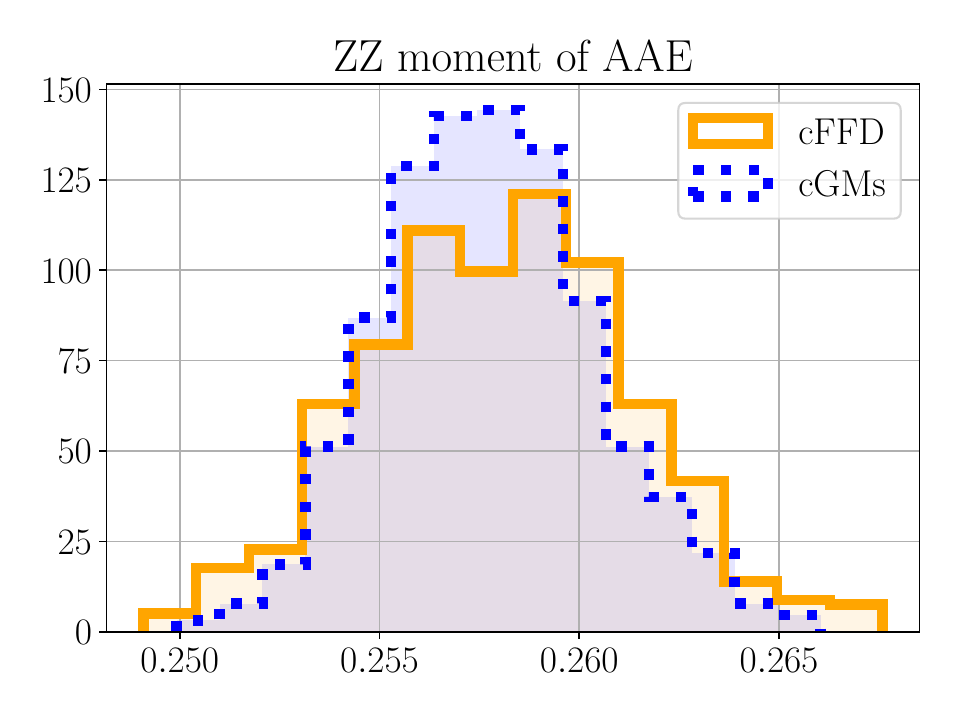}
    \includegraphics[width=0.49\textwidth]{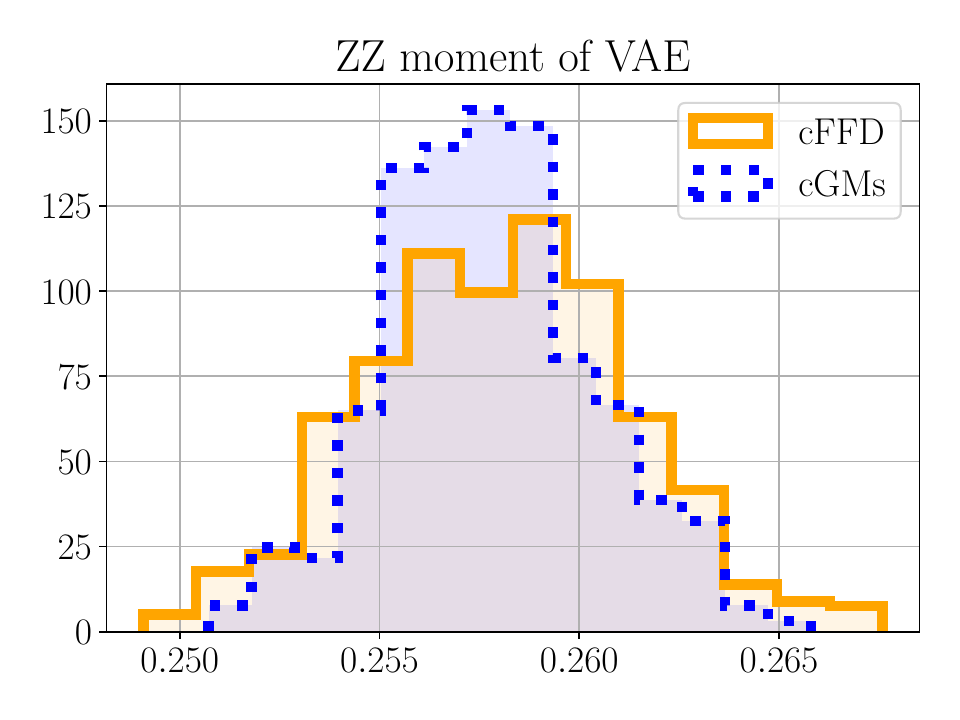}\\
    \includegraphics[width=0.49\textwidth]{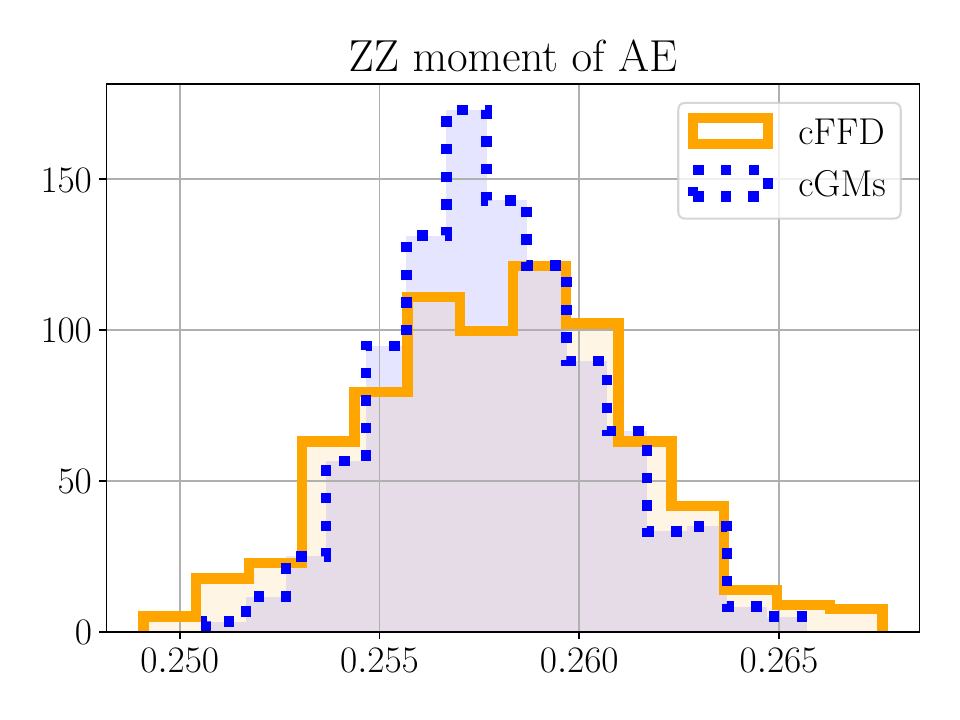}
    \includegraphics[width=0.49\textwidth]{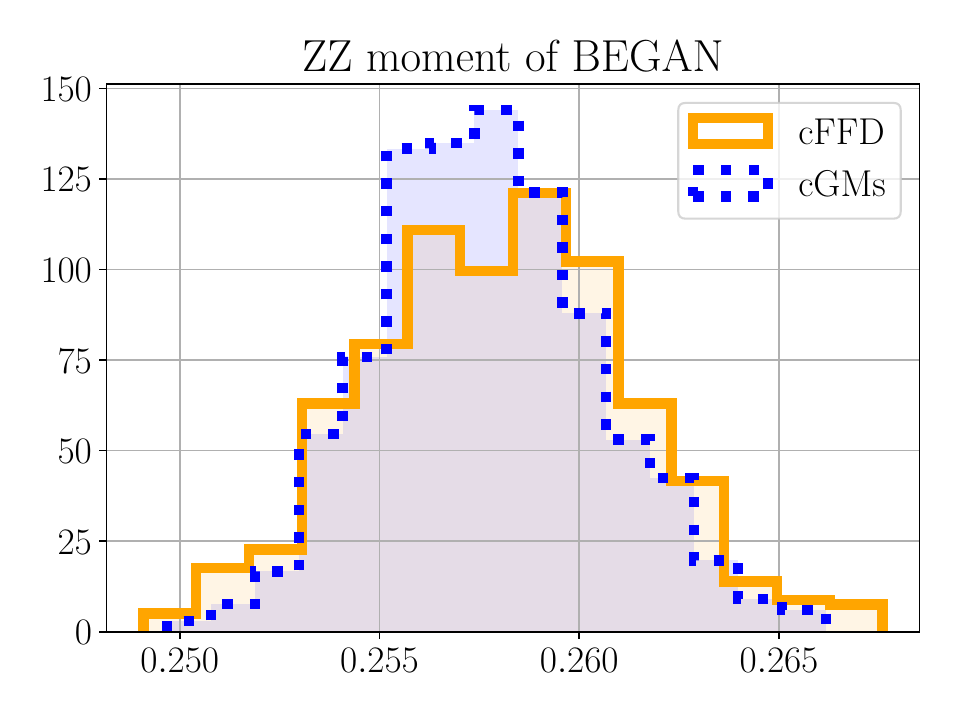}
    \caption{\textbf{SB. }In this figure the histograms of the ZZ component of the inertia tensor of the cFFD and the cGMs are shown. The histogram area intersection is qualitatively more than $50\%$. For a quantitative measure see Table~\ref{tab:SB}. The histograms are obtained from the $\mathbf{n_{\text{test}}=200}$ test samples.}
    \label{fig:hist_I}
\end{figure}

\begin{figure}[htpb!]
    \includegraphics[width=\textwidth]{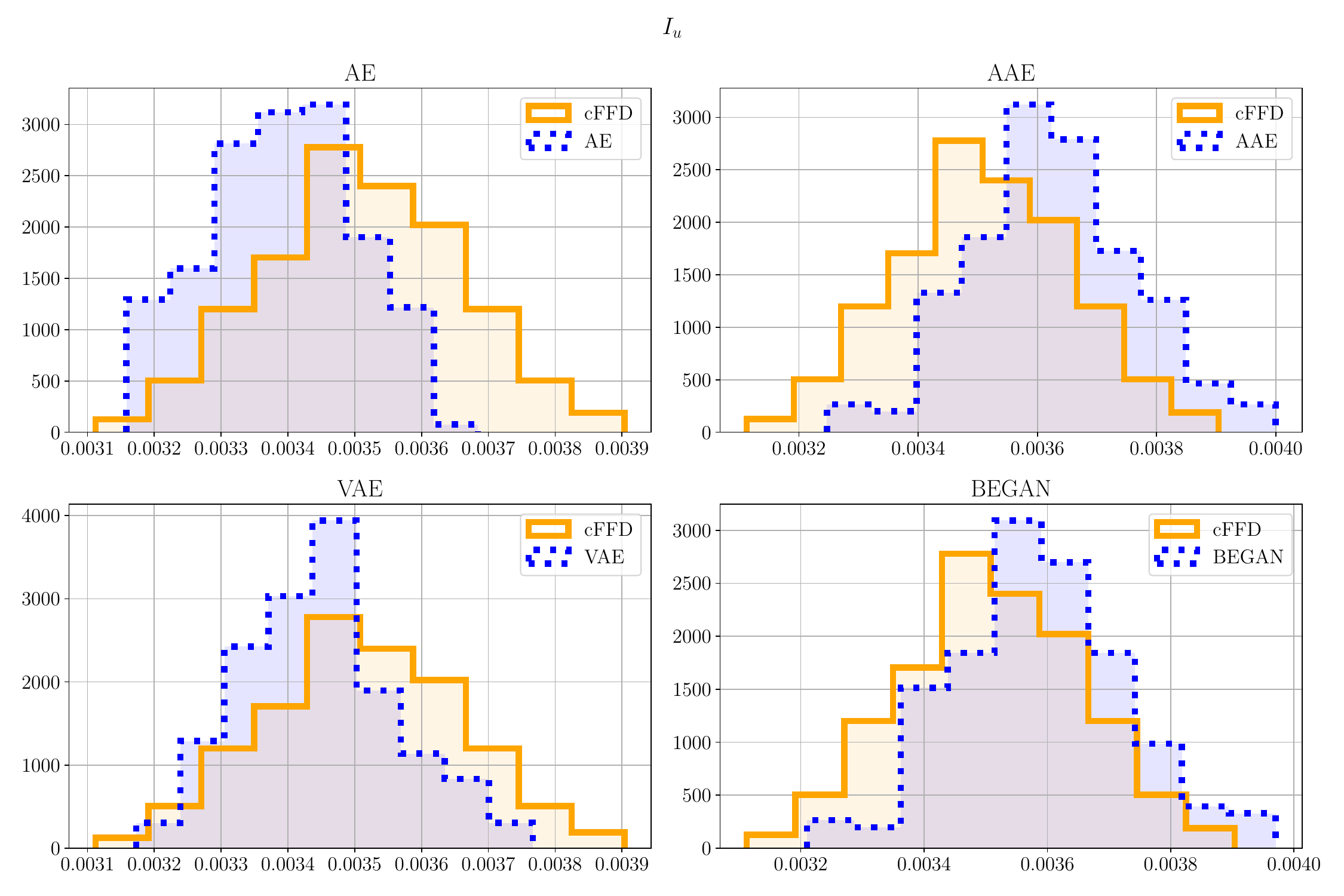}
        \caption{\textbf{SB. }In this figure the histograms of the integral of the solution on the Neumann boundary $\Gamma_N$ of the cFFD and the cGMs are shown. The histogram area intersection is more than $50\%$. For a quantitative measure see Table~\ref{tab:SB}. The histograms are obtained from the $\mathbf{n_{\text{ROM},\text{test}}=20}$ test samples.}
        \label{fig:hist_heat}
\end{figure}

The results with respect to the geometrical and physical metrics defined previously are shown in Table~\ref{tab:SB}. Qualitatively the histograms of each architecture are reported in Figure~\ref{fig:hist_I} for the $I_{zz}$ moment of inertia and in Figure~\ref{fig:hist_heat} for the $I_u$ integral of the solution of the mixed Poisson problem on the Neumann boundary $\Gamma_N$.

\begin{figure}[htpb!]
    \includegraphics[width=\textwidth]{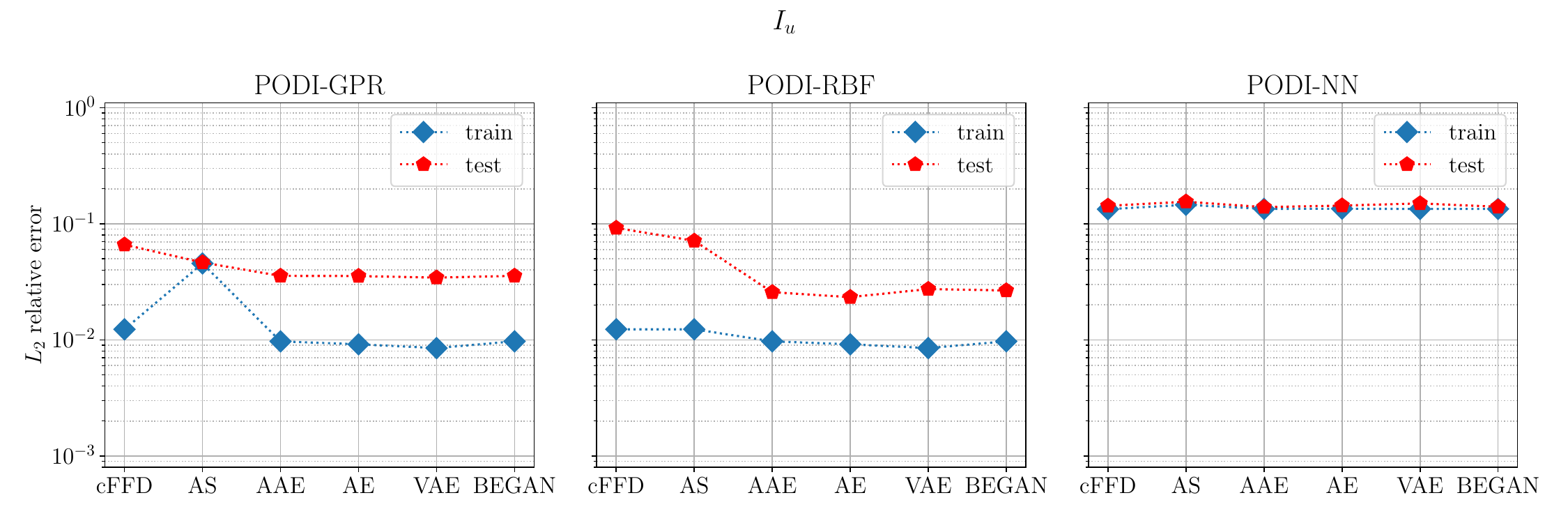}
    \caption{\textbf{SB. }Here we show the ROM performance over the training and test datasets, using different interpolation and regression techniques combined with proper orthogonal decomposition with interpolation (PODI): Gaussian process regression (GPR), radial basis functions interpolation (RBF) and feed-forward neural networks (NN). For every method there is at least a cGM that performs slightly better than the cFFD in terms of accuracy, while reducing the parameters' space dimension from $\mathbf{p=54}$ to $\mathbf{R=15}$. The active subspace (AS) dimension chosen is $\mathbf{r_{\text{AS}}=1}$. The training and test errors are evaluated on $\mathbf{n_{\text{ROM},\text{test}}=80}$ and $\mathbf{n_{\text{ROM},\text{test}}=20}$ training and test samples.}
    \label{fig:SB_rom}
\end{figure}

\begin{figure}[htpb!]
    \includegraphics[width=\textwidth]{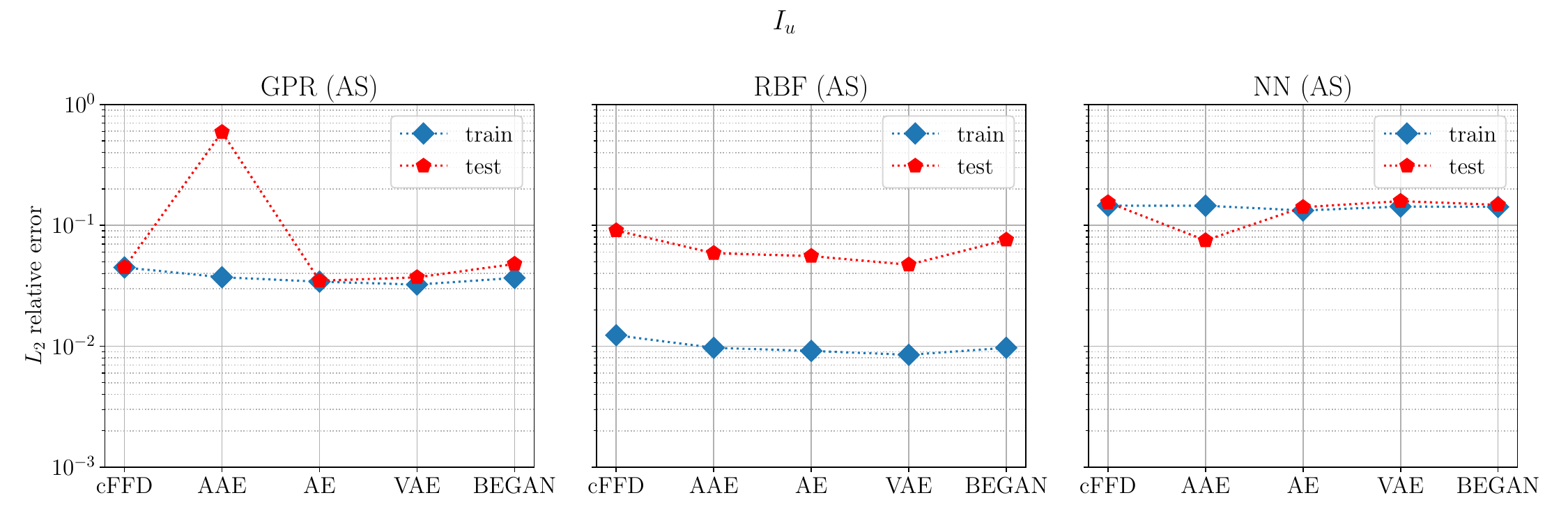}
    \caption{\textbf{SB. }In this figure we show the ROM-PODI performance on the training and test datasets coupled with AS dimension reduction in the space of parameters. It can be seen that the latent dimension can be reduced further with AS $\mathbf{r_{\text{AS}}=1}$ without compromising too much the accuracy. The AS response surface is built over the GPR, RBF, NN interpolations/regressions with an additional GPR from the active one-dimensional variables to the same outputs. The training and test errors are evaluated on $\mathbf{n_{\text{ROM},\text{test}}=80}$ and $\mathbf{n_{\text{ROM},\text{test}}=20}$ training and test samples.}
    \label{fig:SB_rom_AS}
\end{figure}

Non-intrusive model order reduction for the solutions of the mixed Poisson problem with GPR-PODI, RBF-PODI, NN-PODI and AS response surface design with $\mathbf{r_{\text{AS}}=1}$ for the cFFD data is shown in Figure~\ref{fig:SB_rom}. Accurate surrogates models are built even if the dimension of the space of parameters changes from $\mathbf{p=54}$ for the cFFD geometries to $\mathbf{R=15}$ for the cGMs newly generated geometries. We also apply an additional level of reduction in the space of parameters with response surface design with AS in Figure~\ref{fig:SB_rom_AS}. The inputs, in this case, are the one-dimensional active variables $\mathbf{r_{\text{AS}}=1}$ also for the cGMs: the parameters space's dimension changes from $\mathbf{R=15}$ to $\mathbf{r_{\text{AS}}=1}$. For simplicity, we show the AS response surface with dimension $\mathbf{r_{\text{AS}}=1}$, even if from the plot in Figure~\ref{fig:AS_evals} of the first $20$ eigenvalues of the uncentered covariance matrix from equation~\eqref{eq:covSA}, the spectral gap~\cite{constantine2015active} suggests $\mathbf{r_{\text{AS}}=9}$.

The speedup for the generation of the geometries employing the cGMs instead of cFFD is around $60$. The speedup of PODI non-intrusive model order reduction with respect to the full-order simulations is around $126000$.

\begin{table}[htpb!]
    \centering
    \caption{\textbf{SB. }In this table we show the evaluation metrics of the Stanford bunny defined in equations~\eqref{eq:SB_met}: $I_{xx},I_{xy},I_{xz},I_{yy},I_{yz},I_{zz}$ are the components of the inertia tensor, $I_u$ is the integral of the solution of the mixed Poisson problem on the Neumann boundary $\Gamma_N$. The distribution of the components of the inertia tensor is obtained from $\mathbf{n_{\text{test}}=200}$ test samples, while the physical metric $I_u$ is based on  $n_{\text{ROM},\text{test}}=20$ test samples. The BEGAN is overall the best model we managed to train for this test case. The model with the highest output variance is the BEGAN, even though it does not reach the total variance of the training dataset from cFFD that is \textbf{98}.}
    \begin{tabular}{l|l|l|l|l}
        \hline
        \hline
        &   AE &  AAE &  VAE &  BEGAN   \\
        \hline
        \hline
        $JSD(I_{xx})$ & {8.6e-03} & 7.5e-03 & 2.1e-02 & \textbf{6.3e-03} \\ \hline
        $JSD(I_{xy})$ & 1.6e-02  & \textbf{1.2e-02} &2.4e-02 & \textbf{1.2e-02}\\ \hline
        $JSD(I_{xz})$ & \textbf{2.1e-02}	& \textbf{2.1e-02} & 2.6e-02 & \textbf{2.1e-02} \\ \hline
        $JSD(I_{yy})$ & 2.0e-02 & 1.2e-02 & {2.7e-02} & \textbf{8e-03}\\ \hline
        $JSD(I_{yz})$ & 1.9e-02 & \textbf{1.5e-02} & {2.3e-02}  & 1.6e-02\\ \hline
        $JSD(I_{zz})$ & {3.1e-02} & 2.2e-02 & 2.9e-02 & \textbf{1.8e-02}\\ \hline
        $JSD(I_u)$ & 4.0e-01 & 3.5e-01 & \textbf{2.5e-01} & 2.6e-01 \\ \hline
   		$Var$ & 62 & 74 & 58 & \textbf{76} \\
    \end{tabular}
    \label{tab:SB}
\end{table}


\newpage
\subsection{DTCHull bulb (HB)}
\label{subsec:bulb}

The Duisburg test case~\cite{white2019numerical} models a two-phase water-air turbulent incompressible flow over a naval hull. We start from OpenFoam's~\cite{weller1998tensorial} tutorial \textit{DTCHull} related to the \textit{interFoam} multiphase solver for the Reynolds Averaged Navier-Stokes equations (RAS) using the volume of fluid modelling. Water and air are considered as isothermal immiscible fluids. The system of partial equations to be solved is the following
\begin{subequations}
    \begin{align}
        \partial_t (\rho \mathbf{u}) + \nabla\cdot (\rho\mathbf{u}\otimes\mathbf{u})+\nabla p -\rho g-\nabla\cdot\nu\nabla\mathbf{u}-\nabla\cdot R&=0,\qquad(\mathbf{x}, t)\in\Omega\times[0, T]\\
        \nabla\cdot\mathbf{u}&=0,\qquad(\mathbf{x}, t)\in\Omega\times[0, T]\\
        \partial_t \alpha +\nabla\cdot (\mathbf{u}\alpha) &= 0,\qquad(\mathbf{x}, t)\in\Omega\times[0, T]\\
        \alpha\rho_W+(1-\alpha)\rho_A &=\rho,\qquad (\mathbf{x}, t)\in\Omega\times[0, T]\\
        \alpha\nu_W+(1-\alpha)\nu_A &=\nu,\qquad (\mathbf{x}, t)\in\Omega\times[0, T]
    \end{align}
where $\mathbf{u}$ is the velocity field, $p$ is the pressure field, $\rho_W,\rho_A$ are the densities of water and air, $\nu_W, \nu_A$ are the dynamic viscosities of water and air, $R$ is the Reynolds' stress tensor, $g$ is the acceleration of gravity and $\alpha$ represents the interphase between the fluids with values from $0$ (inside the air phase) to $1$ (inside the water phase). The turbulence is modelled with the $\kappa-\omega$ Shear Stress Transport (SST) model~\cite{menter1993zonal}. The initial conditions are:
\begin{align}
    \mathbf{u}(\mathbf{x})=(U_0, 0, 0),\qquad\mathbf{x}\in\Omega\\
    p(\mathbf{x})=0,\qquad\mathbf{x}\in\Omega\\
    \alpha(\mathbf{x})=\alpha_0(\mathbf{x}),\qquad\mathbf{x}\in\Omega
\end{align}
\end{subequations}
where $U_0\in\mathbb{R}$ is the initial velocity, $\alpha_0$ is the initial water level and the boundary conditions are specified in the OpenFoam DTCHull tutorial of the \textit{interFoam} solver. We search for steady-state solutions with a final pseudo time instant $T=4000$. The average physical quantities we will evaluate such as the drag and the angular momentum along the z-axis are obtained as the mean over the last $500$ pseudo-time instants.

\begin{figure}
    \centering
    \includegraphics[width=0.49\textwidth]{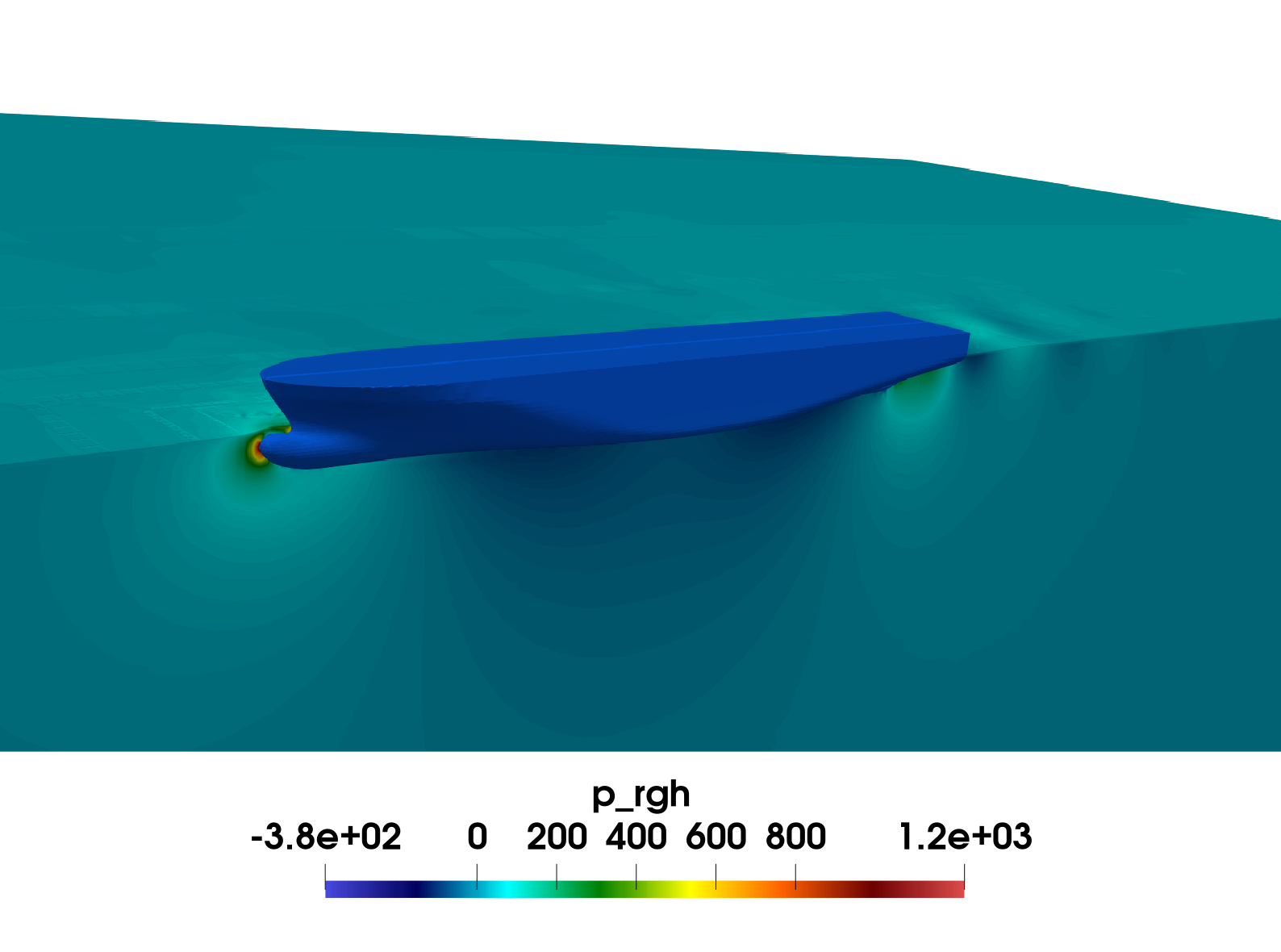}
    \includegraphics[width=0.49\textwidth]{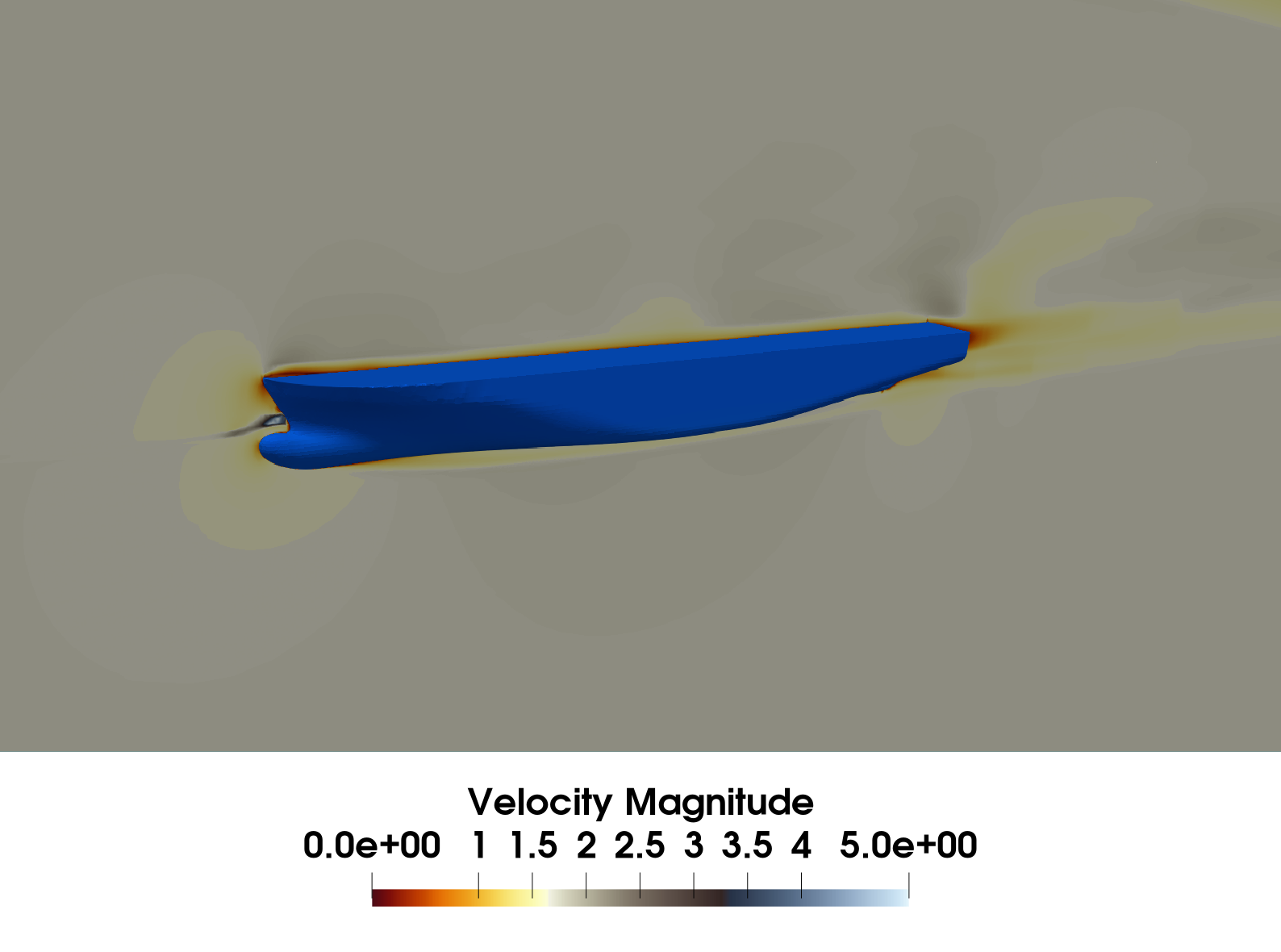}
    \caption{\textbf{HB. }\textbf{Left:} pressure field minus the hydrostatic pressure contribution, it is computed on the last pseudo-time instant $T=4000$. \textbf{Right:} velocity magnitude computed on the last pseudo-time instant $T=4000$.}
    \label{fig:p_v}
\end{figure}

\begin{figure}
    \centering
    \includegraphics[width=0.32\textwidth]{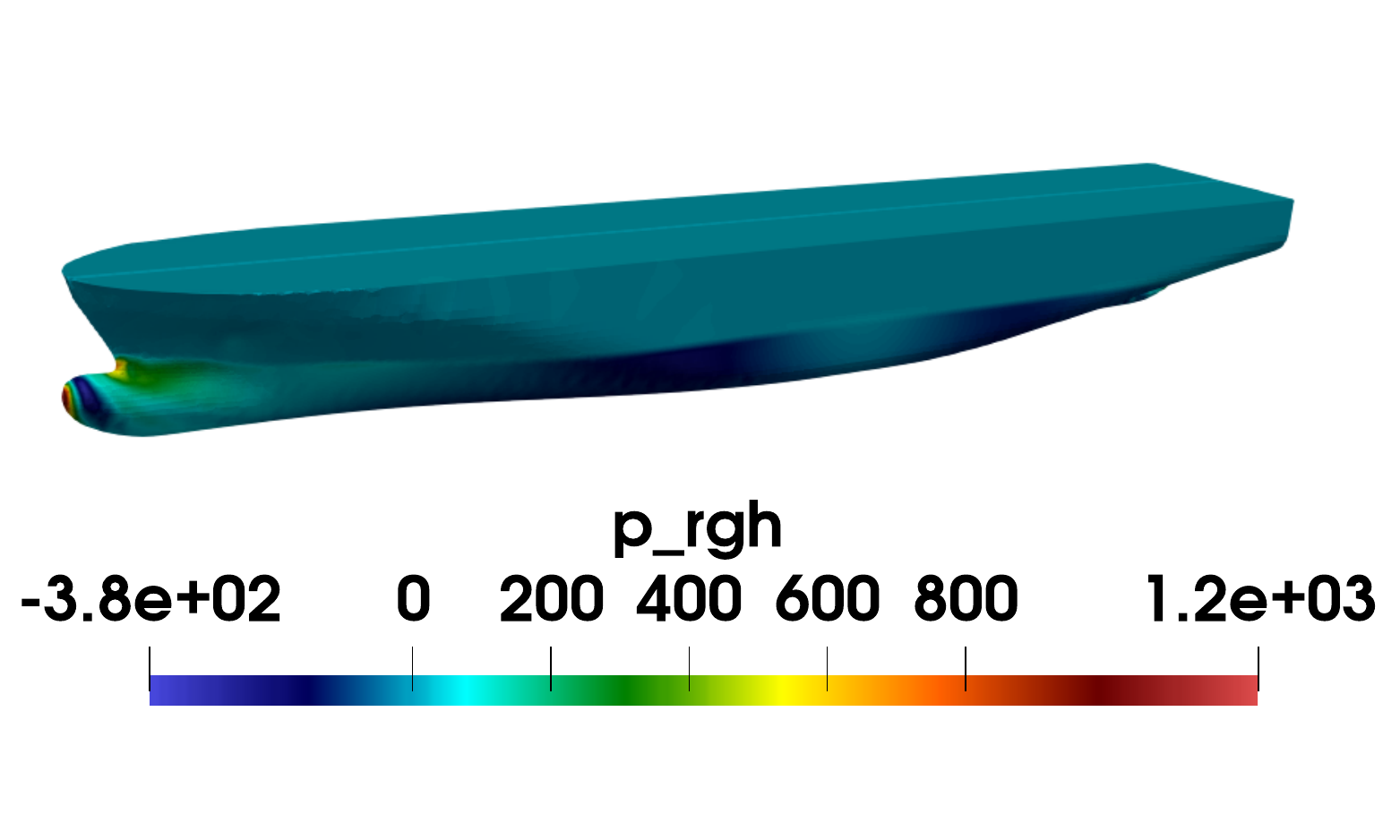}
    \includegraphics[width=0.32\textwidth]{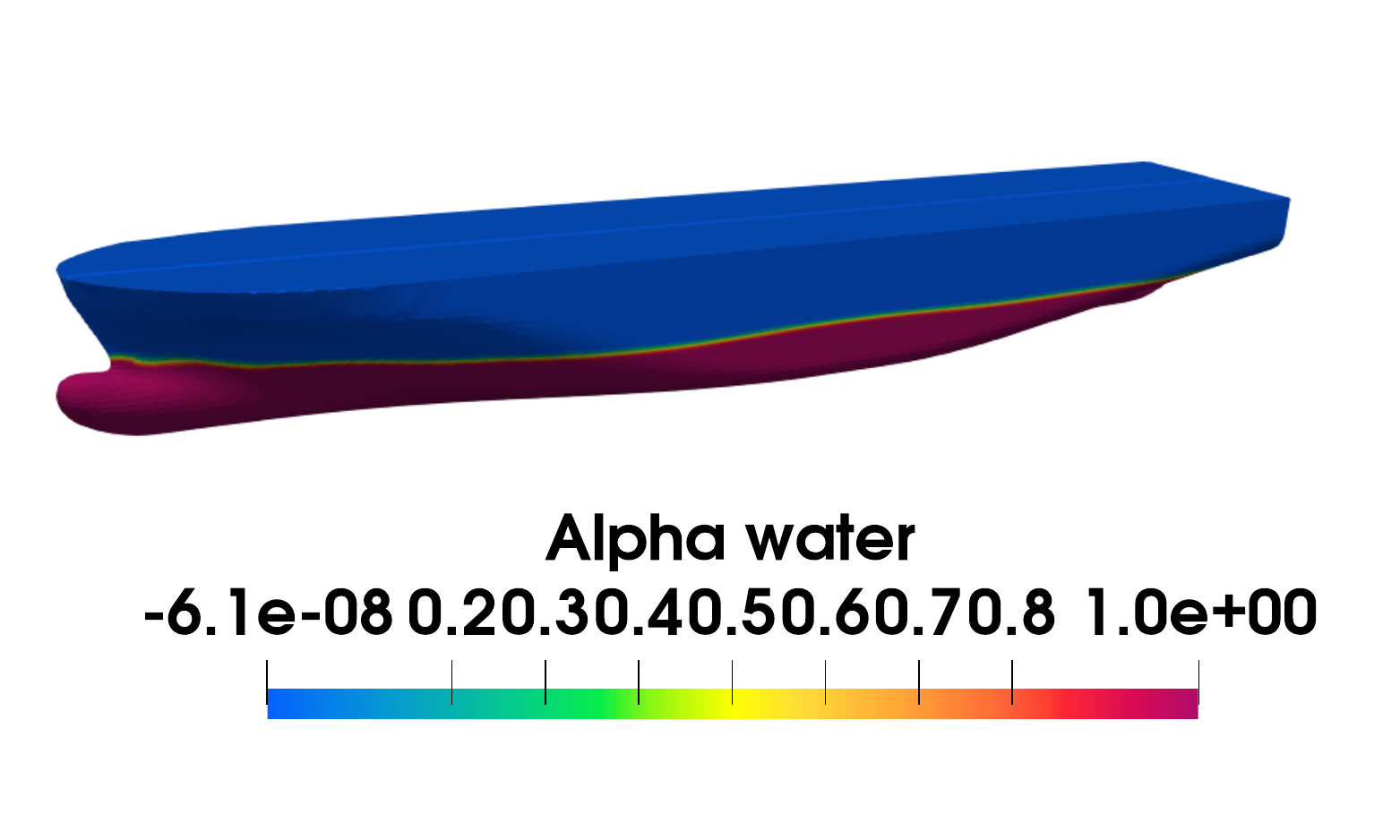}
    \caption{\textbf{HB. }\textbf{Left:} pressure field minus the hydrostatic pressure contribution computed on the last pseudo-time instant $T=4000$ on the hull's surface. \textbf{Right:} interphase field $\alpha$ on the hull at the final pseudo-time instant $T=4000$: the value of $1$ correspond to the water phase and the value of $0$ to the air phase.}
    \label{fig:hull_p_a}
\end{figure}

\begin{figure}
    \centering
    \includegraphics[width=0.49\textwidth]{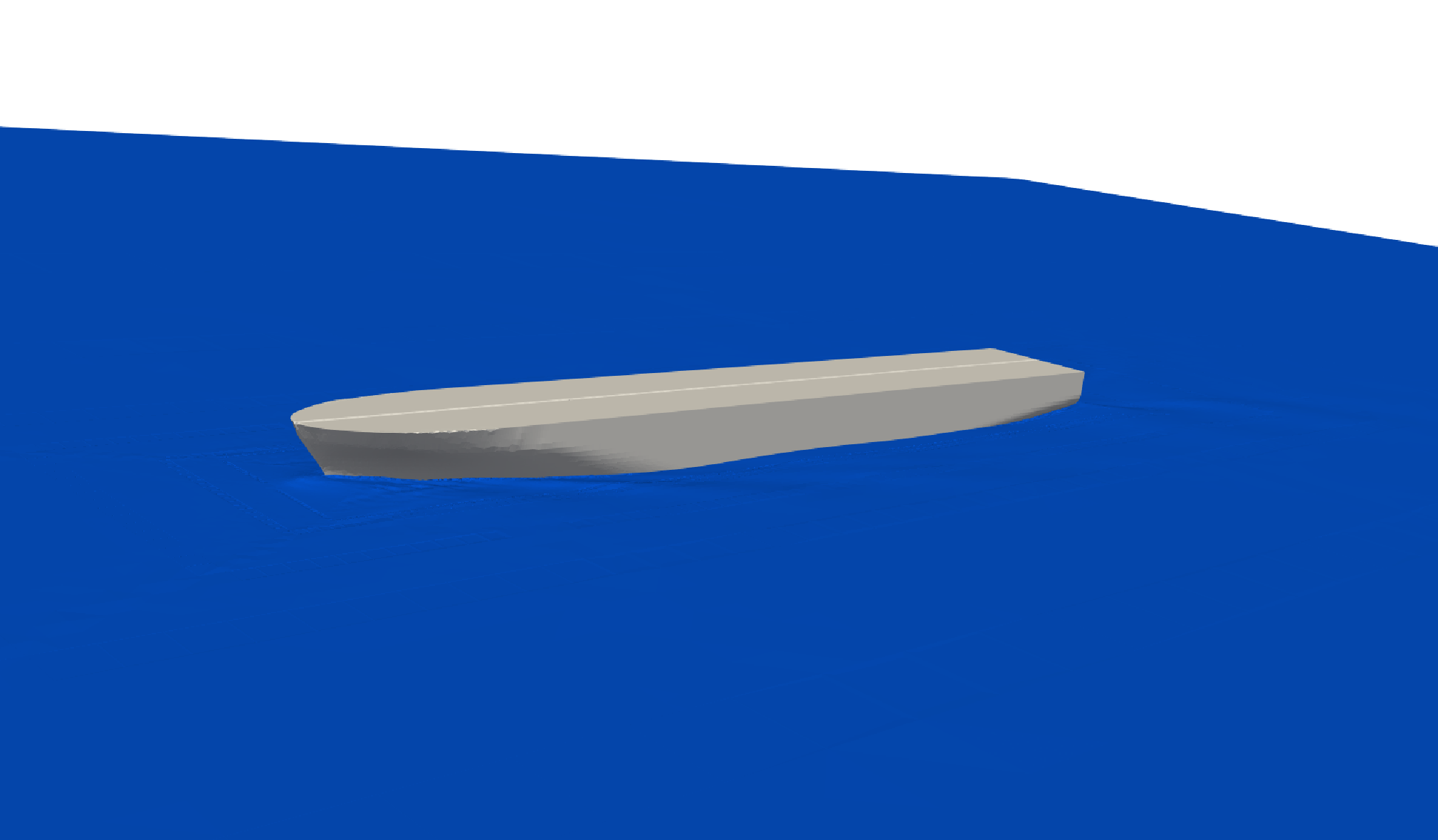}
    \includegraphics[width=0.49\textwidth]{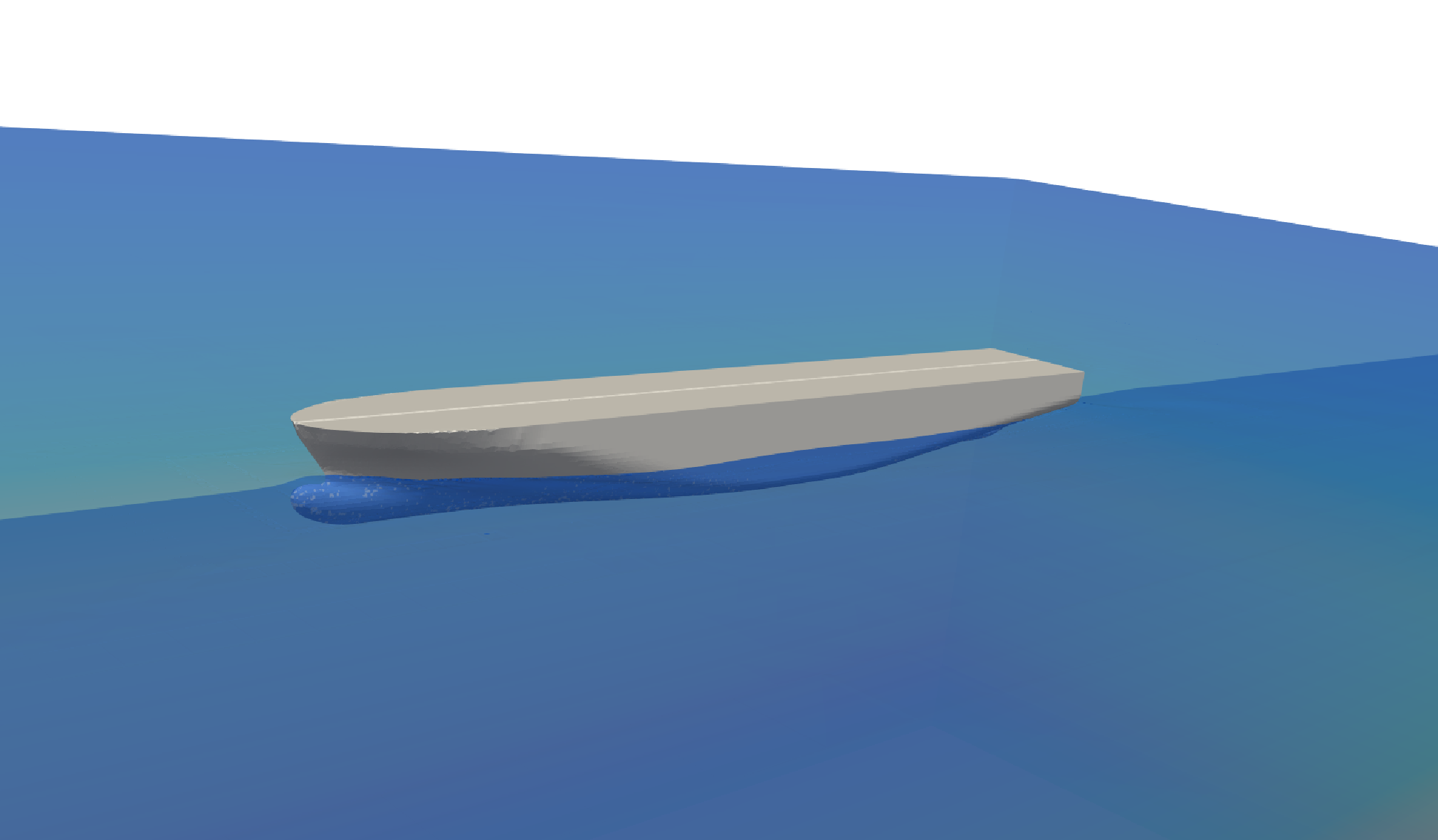}
    \caption{\textbf{HB. }Interphase field $\alpha$ on the whole domain at the final pseudo-time instant $T=4000$: only the values less than $0.5$ of the field $\alpha$ are shown, representing the water level.}
    \label{fig:alpha}
\end{figure}

Some physical fields of interest at the final pseudo-time instant are shown in Figure~\ref{fig:p_v} for the effective pressure and the velocity magnitude, in Figure~\ref{fig:hull_p_a} for the effective pressure on the hull surface and the interphase field $\alpha$ and in Figure~\ref{fig:alpha} for the interphase field $\alpha$ on the whole domain. Only half of the hull is employed for the numerical simulations.

For the DTCHull test case \textbf{HB} we consider the moments of inertia of the hull, the angular momentum along the z-axis of the hull, the surface area of the hull's bulb, and the drag on the hull:
\begin{subequations}
    \label{eq:HB_met}
    \begin{align}
        I_{xx}=\int_{\Omega_{\text{hull}}} r_{X}^2(\mathbf{x})d\mathbf{x},\ I_{yy}=\int_{\Omega_{\text{hull}}} r_{Y}^2(\mathbf{x})d\mathbf{x},\ I_{zz}=\int_{\Omega_{\text{hull}}} r_{Z}^2(\mathbf{x})d\mathbf{x},\qquad\text{(principal moments of inertia)}\label{eq:pmomentsDTCH}\\
        I_{xy}=\int_{\Omega_{\text{hull}}} r_{X}(\mathbf{x})r_{Y}(\mathbf{x})d\mathbf{x},\ I_{xz}=\int_{\Omega_{\text{hull}}} r_{X}(\mathbf{x})r_{Z}(\mathbf{x})d\mathbf{x},\ I_{yz}=\int_{\Omega_{\text{hull}}} r_{Y}(\mathbf{x})r_{Z}(\mathbf{x})d\mathbf{x},\qquad\text{(moments of inertia)}\label{eq:momentsDTCH}\\
        M_{z}= \int_{\Omega_{\text{hull}}} r_{z}^2(\mathbf{x})\times\mathbf{u}\ d\mathbf{x},\qquad\text{(angular momentum along z-axis)}\label{eq:zmomentDTCH}\\
        A_{\text{bulb}}=\int_{\partial\Omega_{bulb}} d\sigma,\qquad\text{(surface area of the bulb)}\label{eq:ADTCH}\\
        c_{\mathrm{d}}=\frac{1}{A_{hull}(\mathbf{u}\cdot\mathbf{e}_x) ^2} \left(\varointclockwise_{\delta\Omega_{\text{hull}}}p\mathbf{n}-[\nu_{\alpha}(\nabla \mathbf{u}+\nabla \mathbf{u}^{T})]\mathbf{n} ds\right)\cdot\mathbf{e}_x,\qquad\text{(drag on the hull)}\label{eq:dragDTCH}
    \end{align}
\end{subequations}
where the moments of inertia $I_{xx},I_{yy},I_{zz},I_{xy},I_{xz},I_{yz}$ and surface area $A_{\text{bulb}}$ are evaluated from the discrete STL file, and the angular momentum along the z-axis $M_z$ and the drag coefficient $c_d$ are evaluated from the computational mesh.

\begin{figure}
    \centering
    \includegraphics[width=0.12\textwidth, trim={200 100 500 100}, clip]{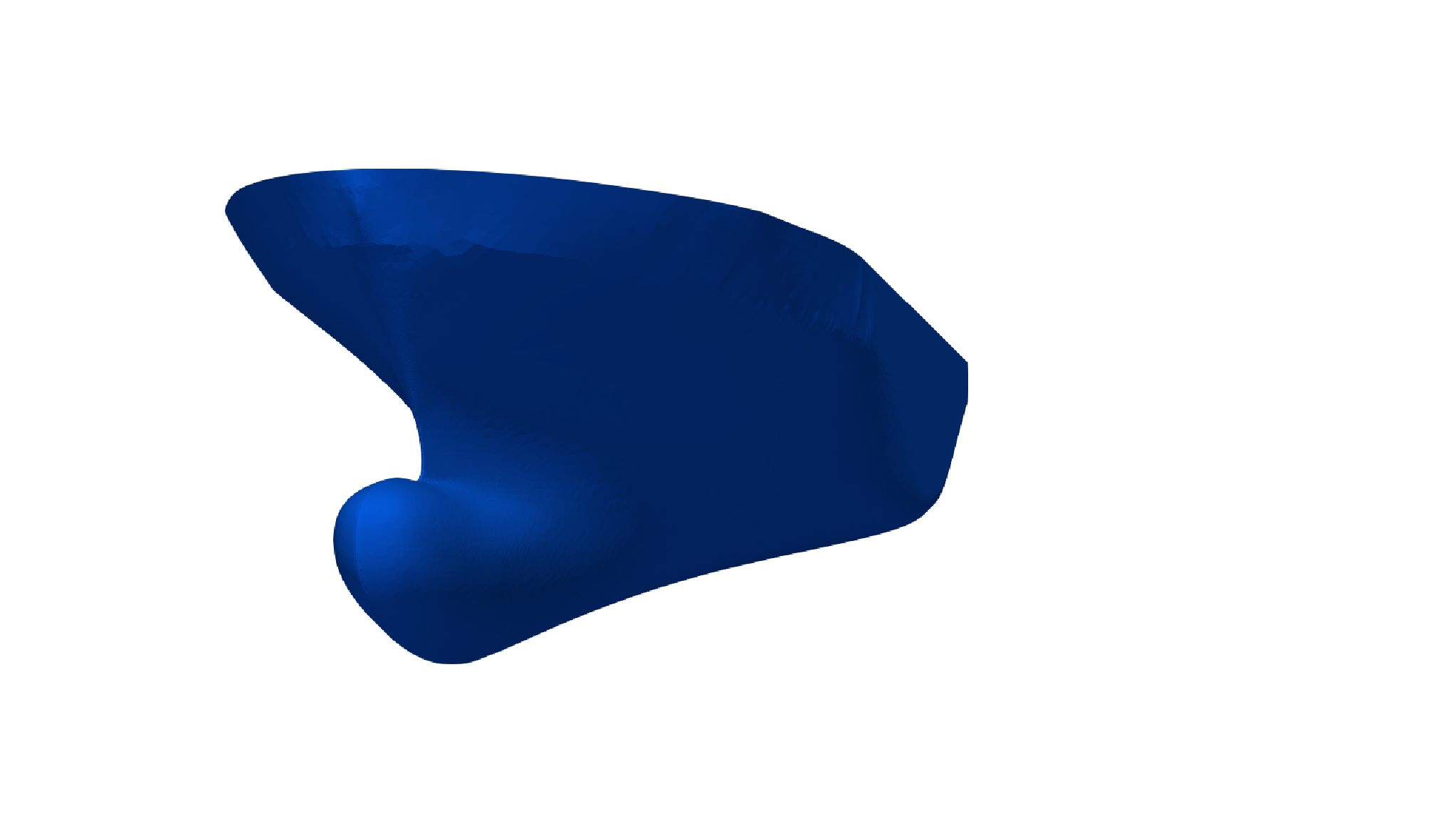}
    \includegraphics[width=0.12\textwidth, trim={200 100 500 100}, clip]{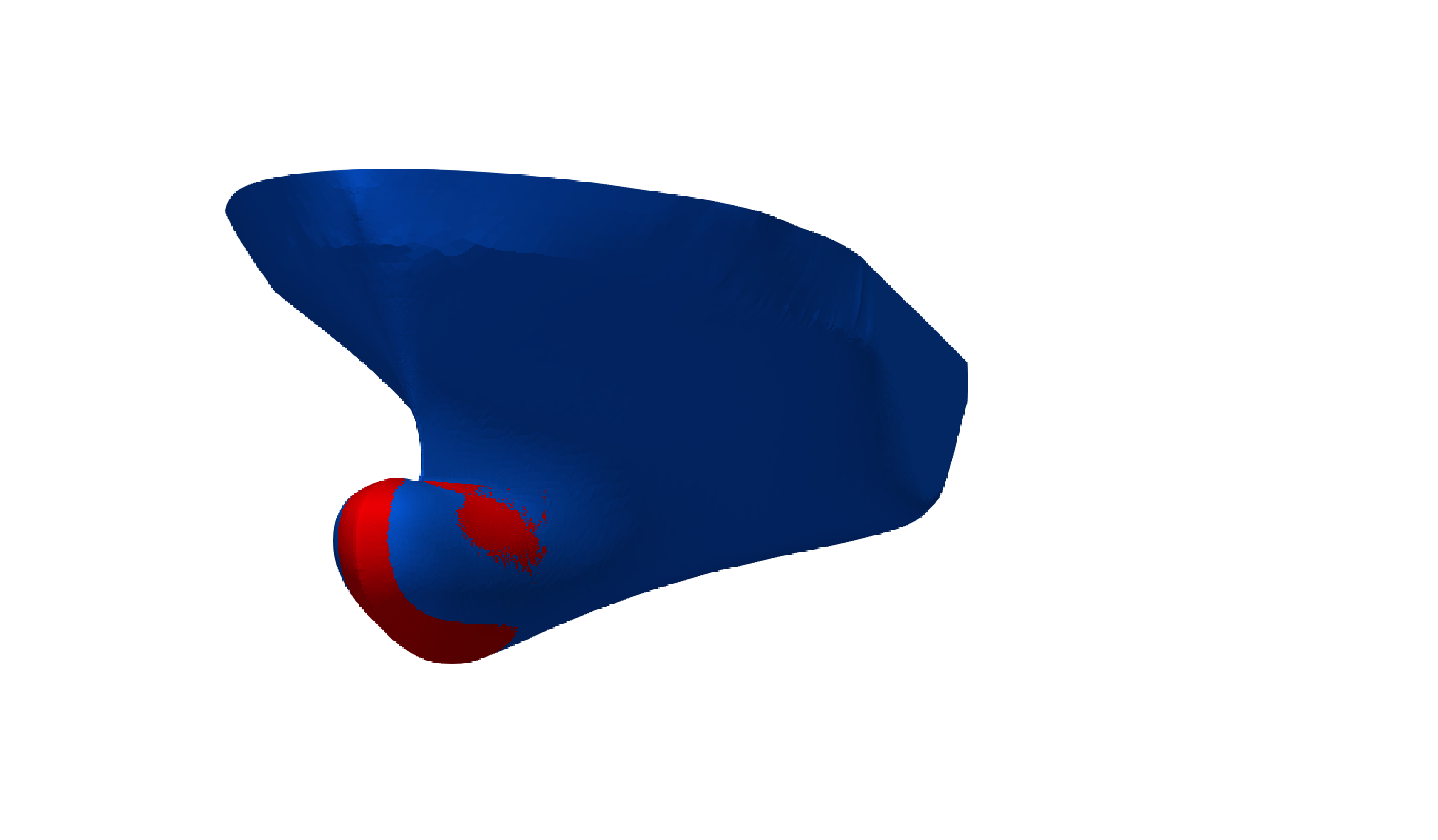}
    \includegraphics[width=0.12\textwidth, trim={200 100 500 100}, clip]{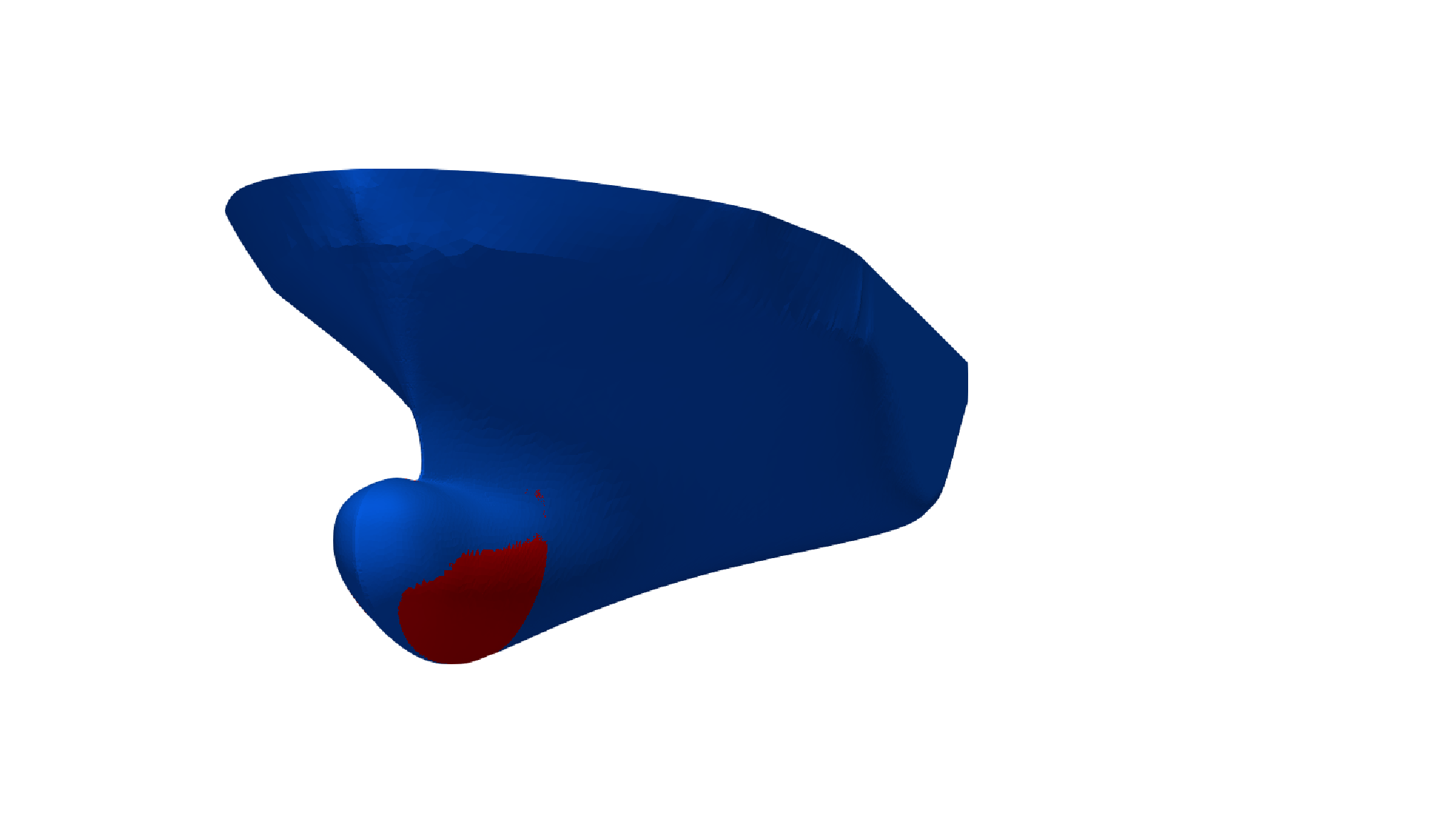}
    \includegraphics[width=0.12\textwidth, trim={200 100 500 100}, clip]{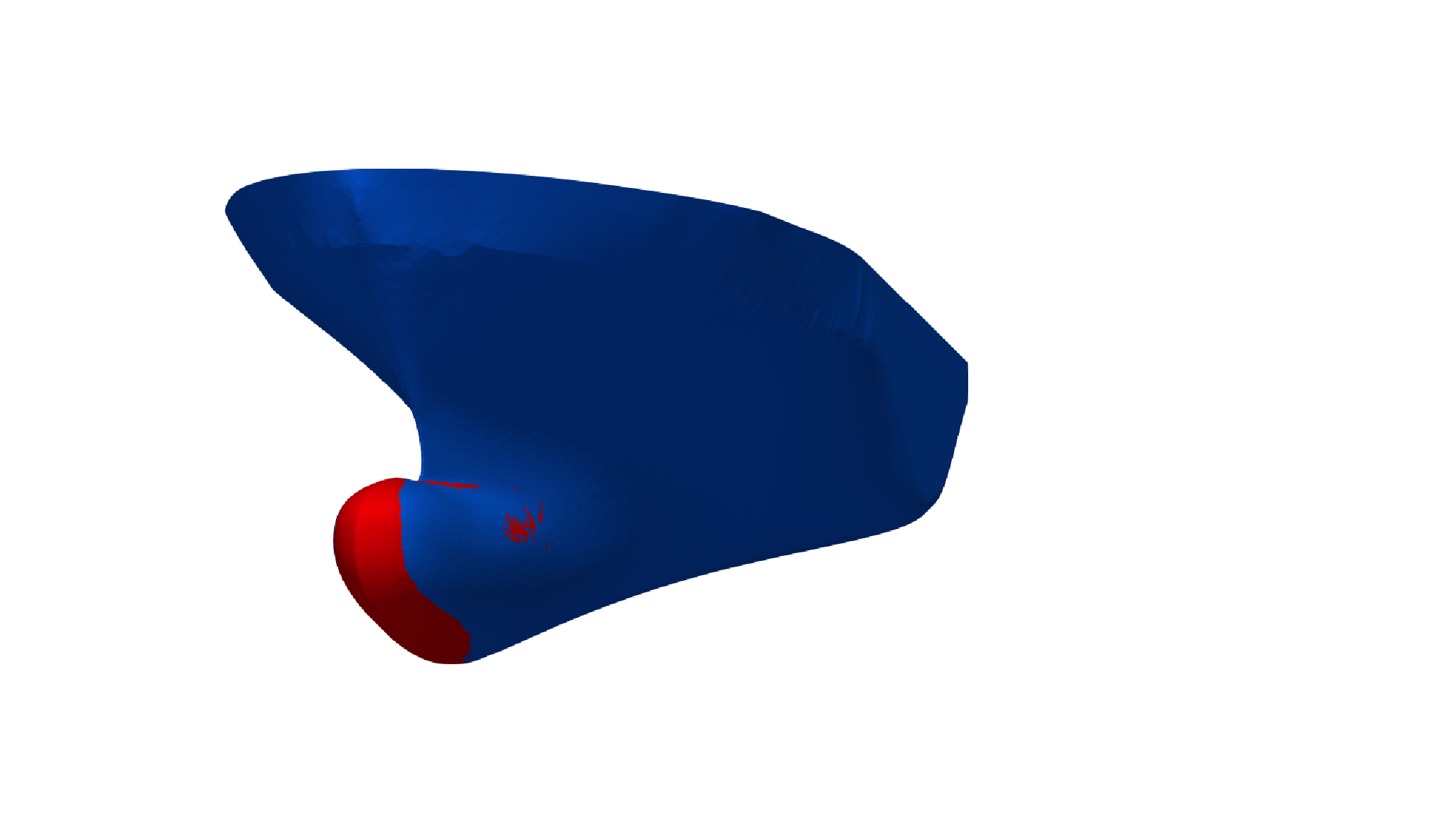}
    \includegraphics[width=0.12\textwidth, trim={200 100 500 100}, clip]{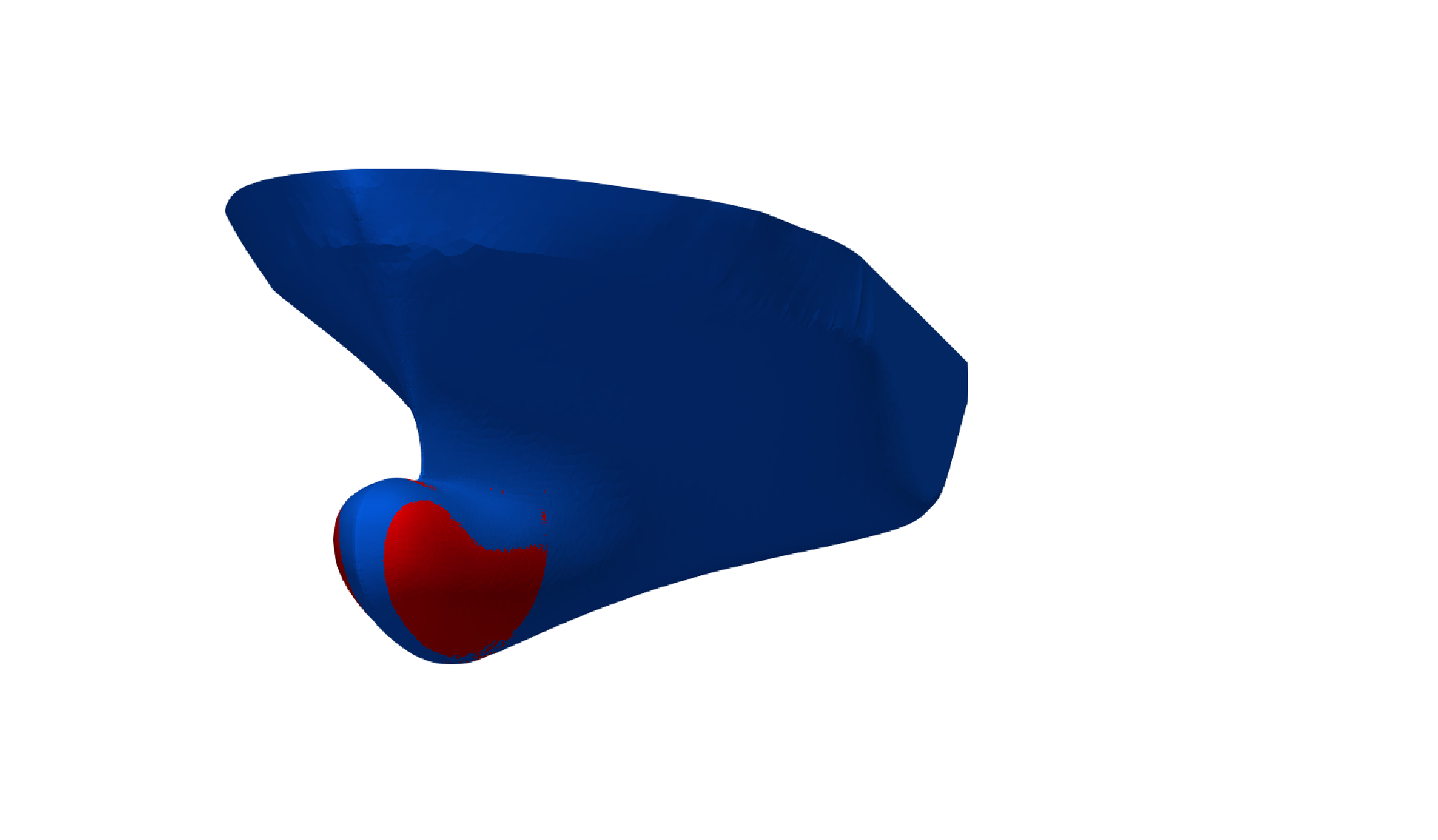}\\
    \includegraphics[width=0.12\textwidth, trim={200 100 500 100}, clip]{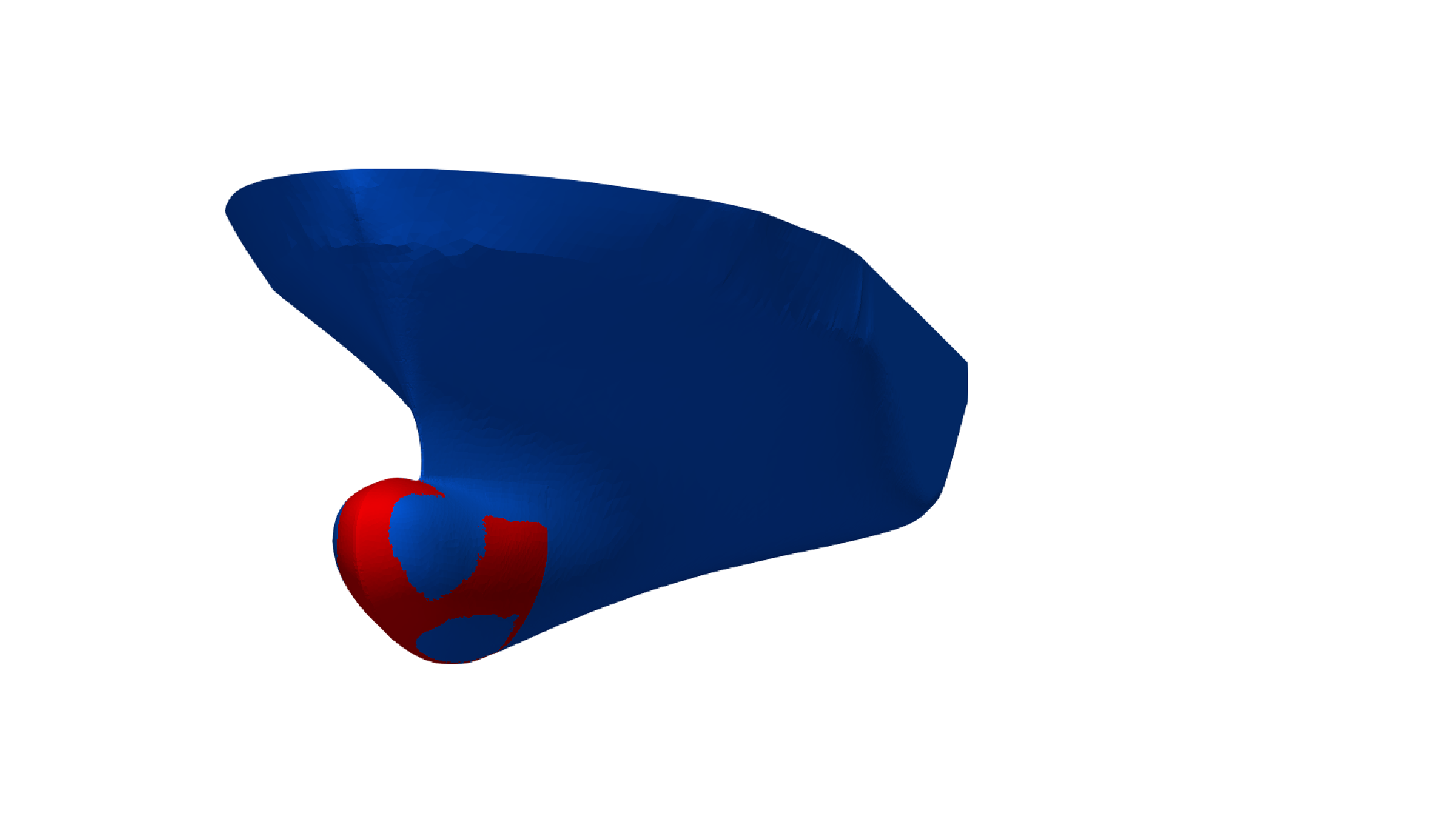}
    \includegraphics[width=0.12\textwidth, trim={200 100 500 100}, clip]{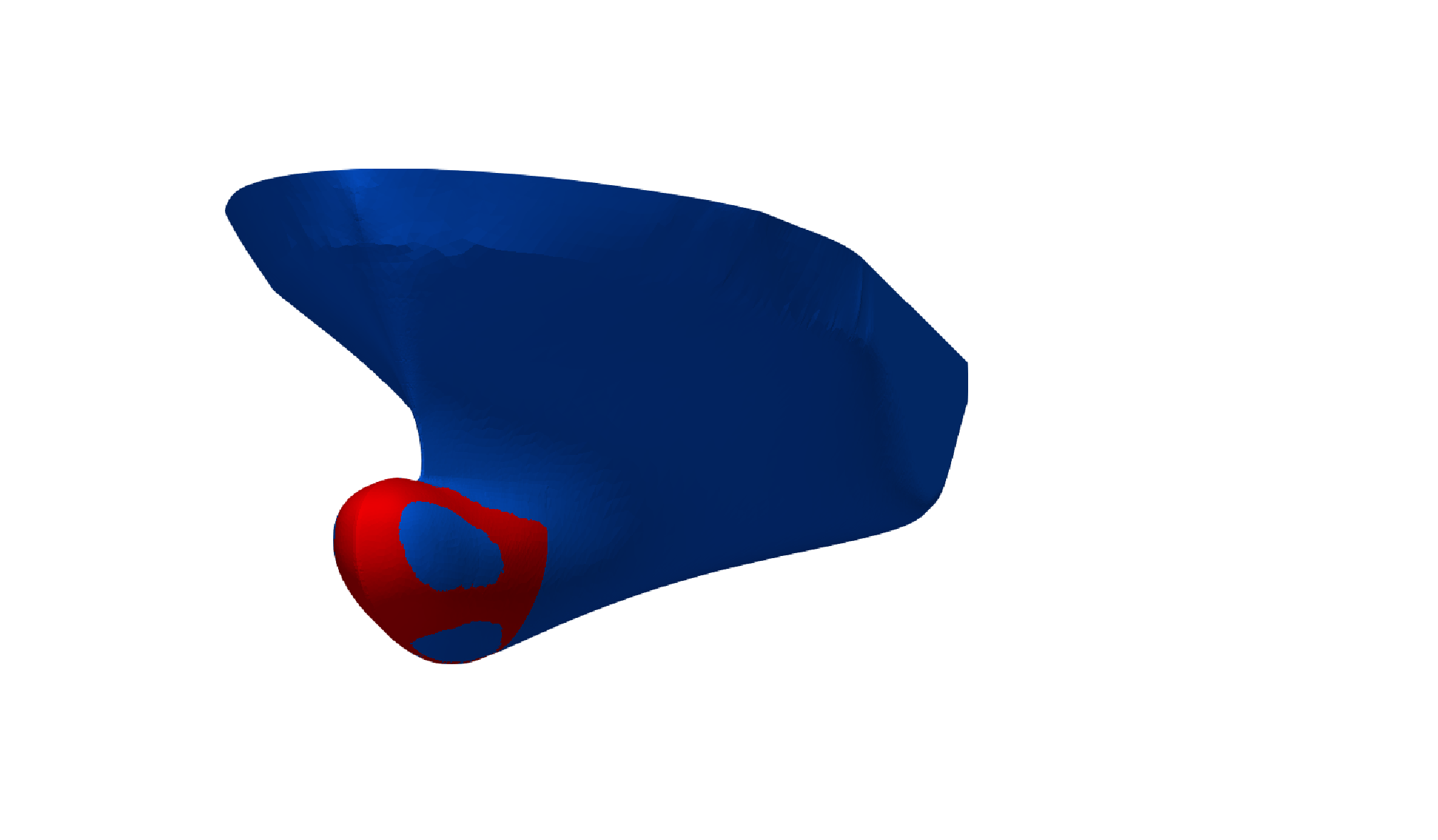}
    \includegraphics[width=0.12\textwidth, trim={200 100 500 100}, clip]{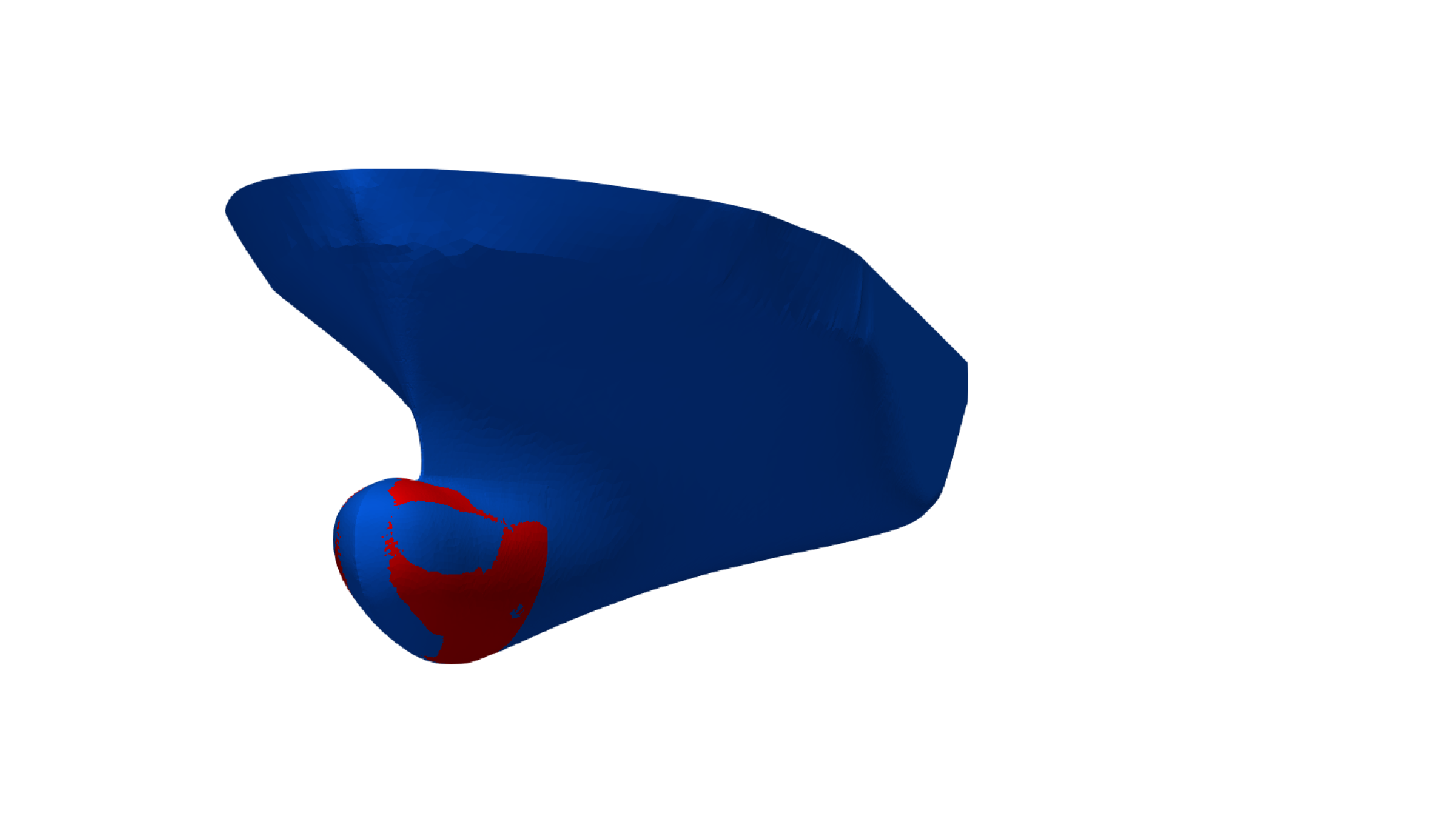}
    \includegraphics[width=0.12\textwidth, trim={200 100 500 100}, clip]{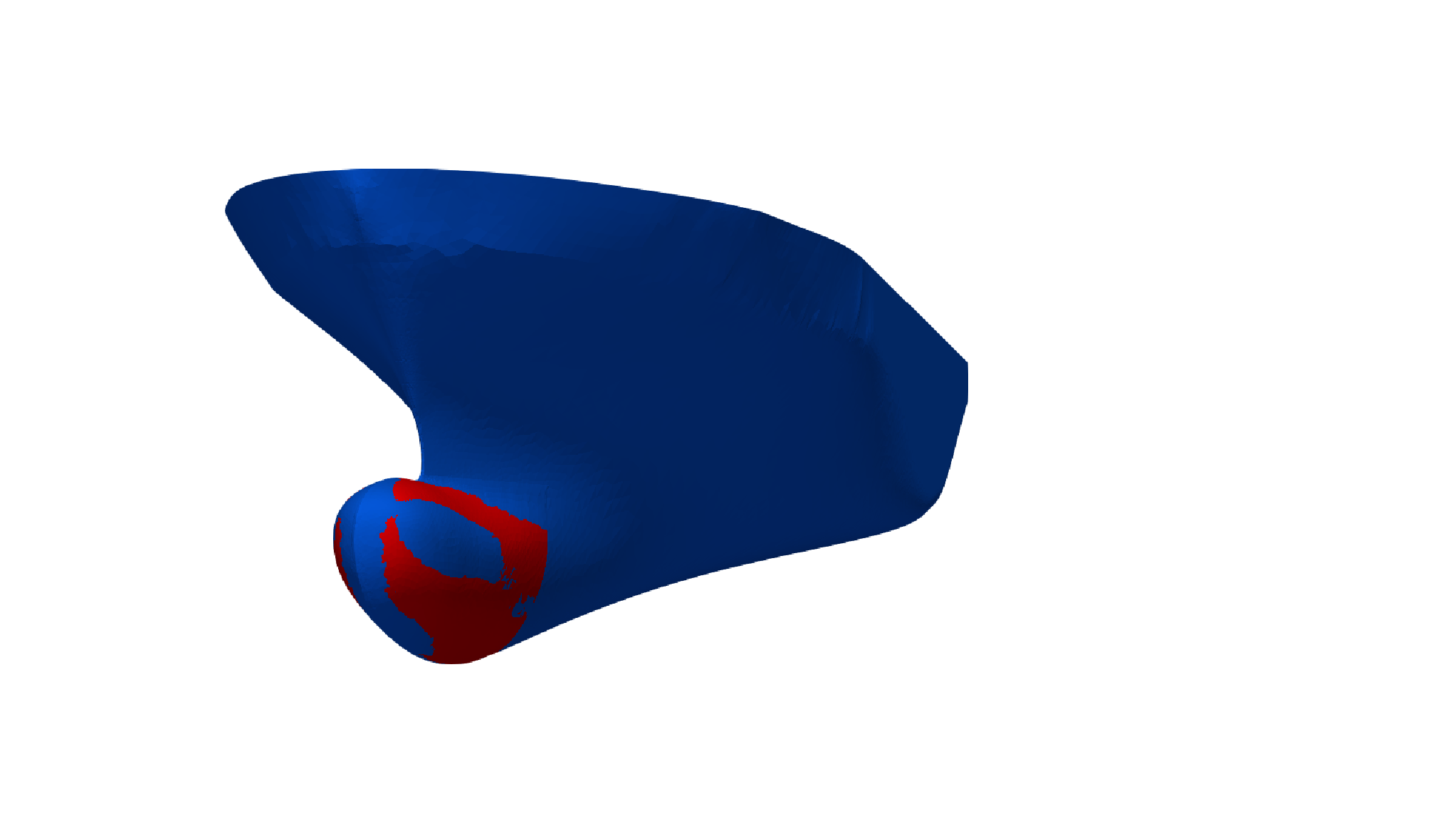}
    \includegraphics[width=0.12\textwidth, trim={200 100 500 100}, clip]{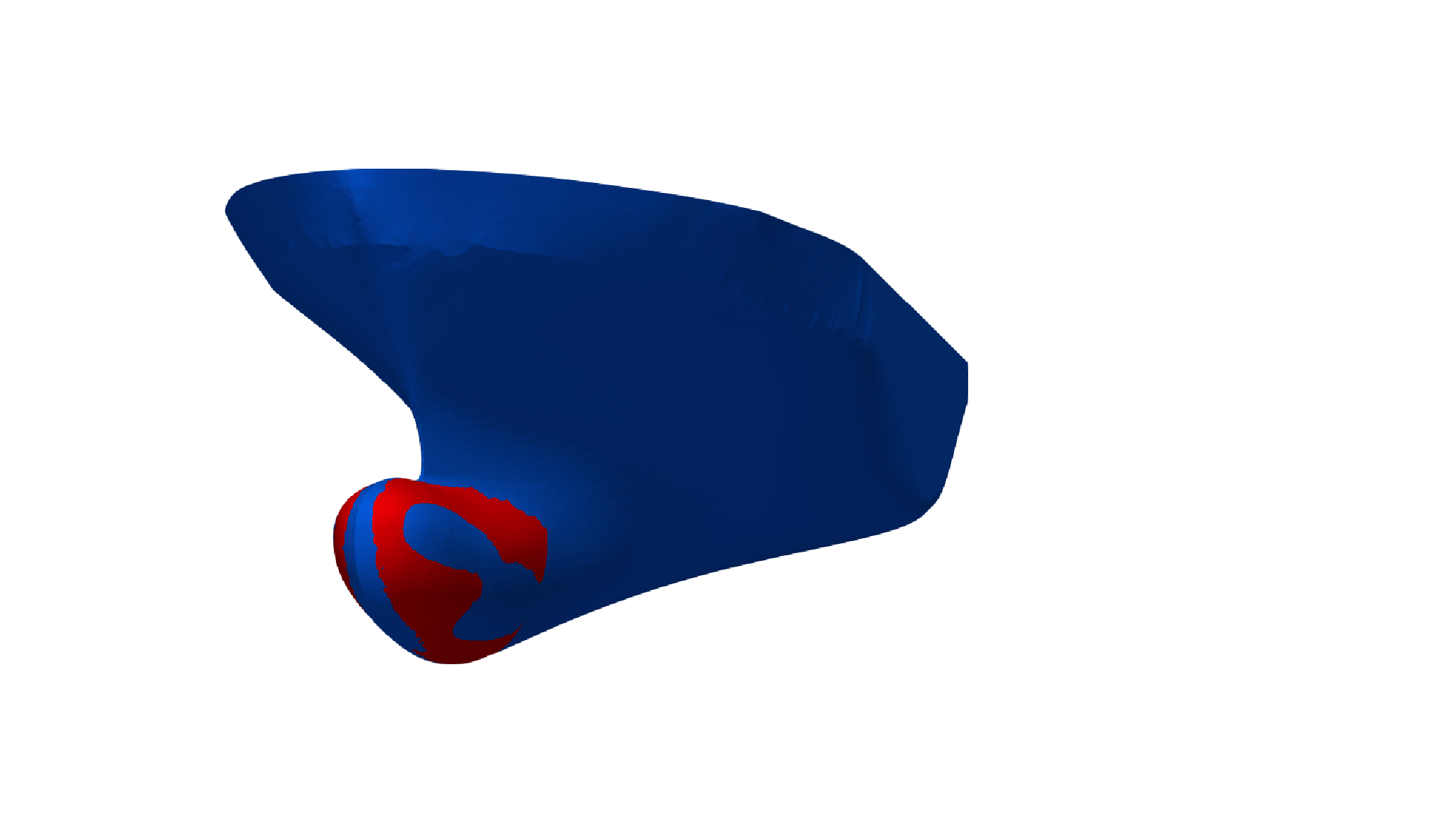}\\
    \includegraphics[width=0.12\textwidth, trim={200 100 500 100}, clip]{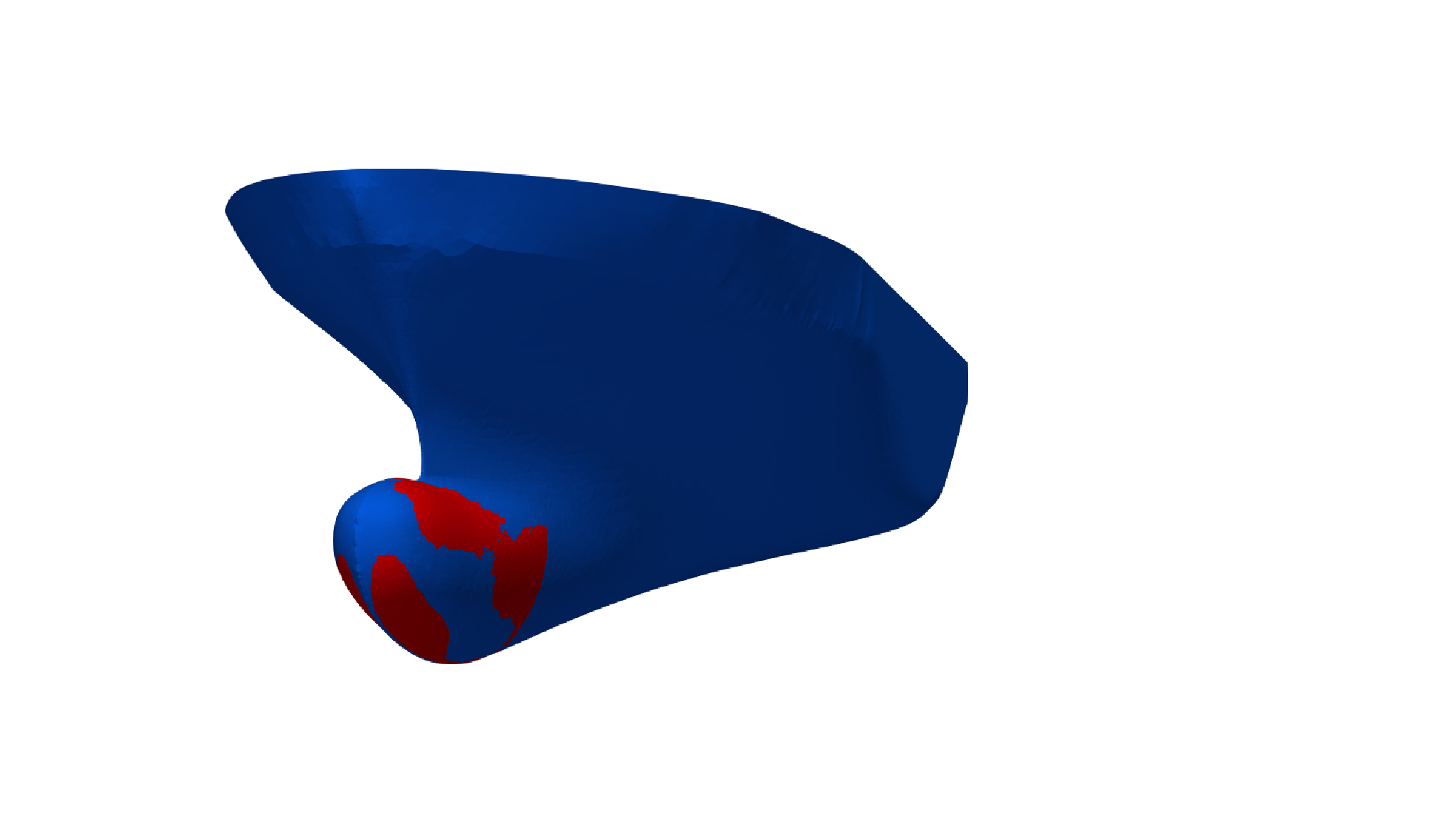}
    \includegraphics[width=0.12\textwidth, trim={200 100 500 100}, clip]{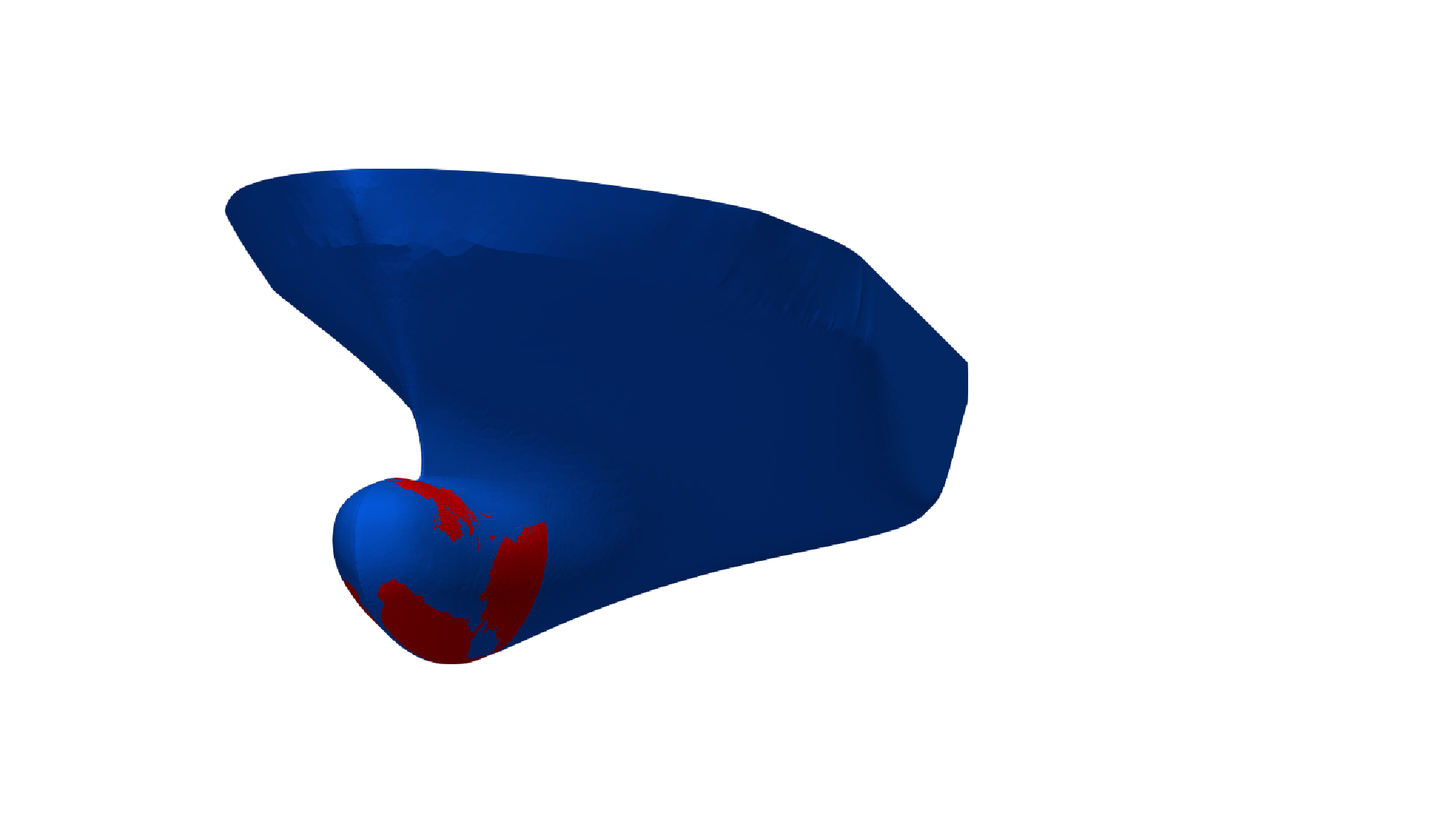}
    \includegraphics[width=0.12\textwidth, trim={200 100 500 100}, clip]{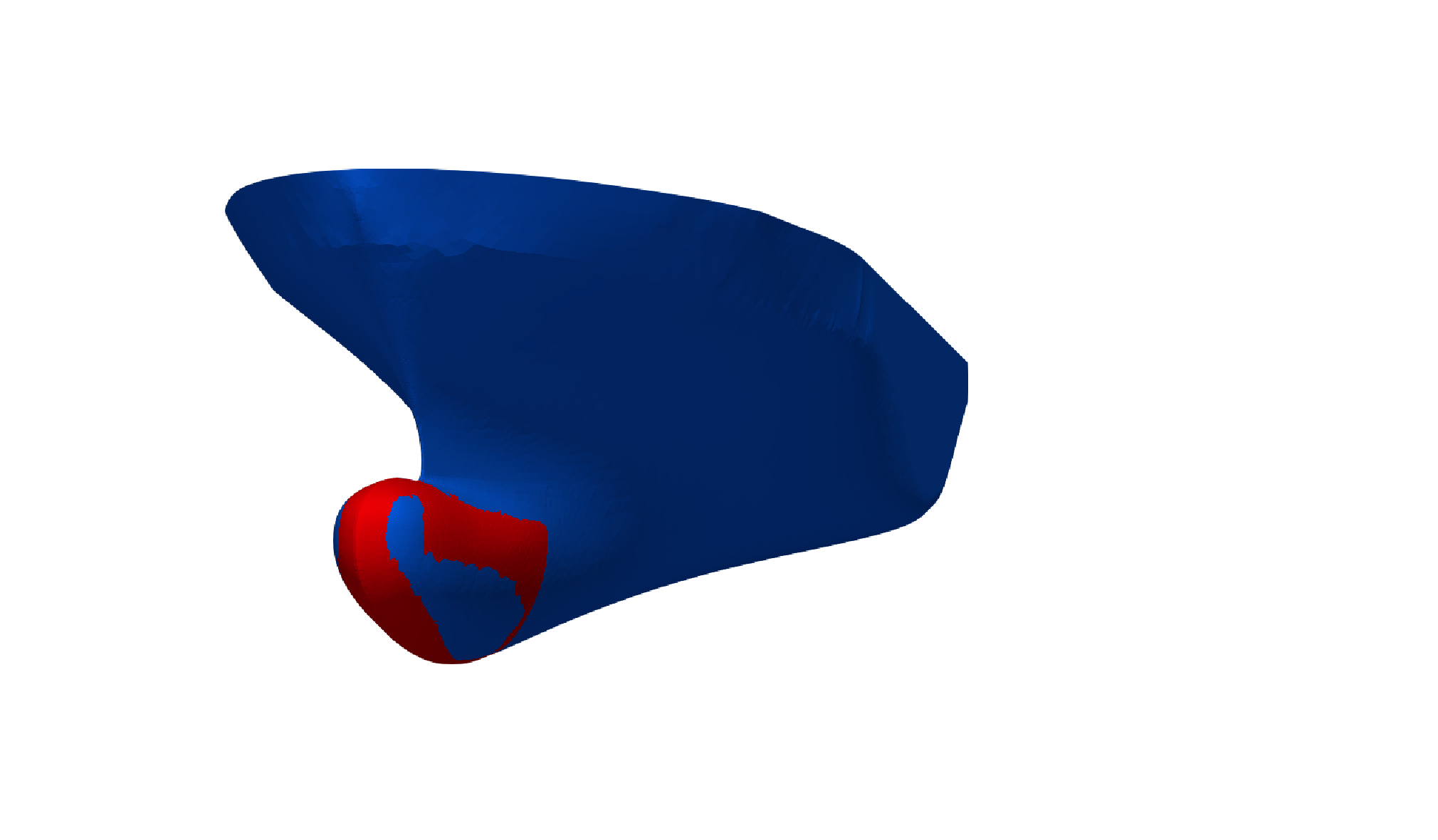}
    \includegraphics[width=0.12\textwidth, trim={200 100 500 100}, clip]{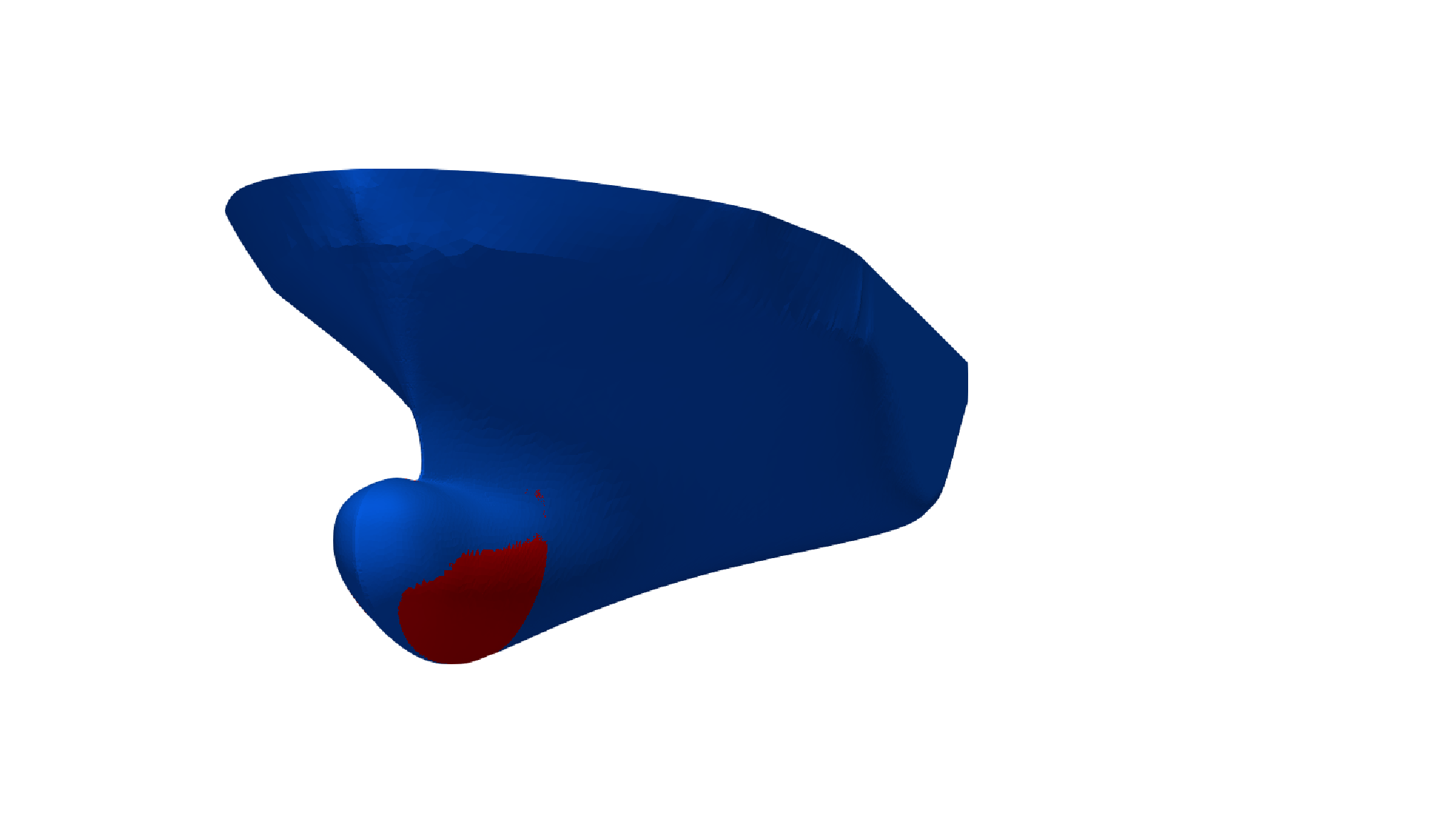}
    \includegraphics[width=0.12\textwidth, trim={200 100 500 100}, clip]{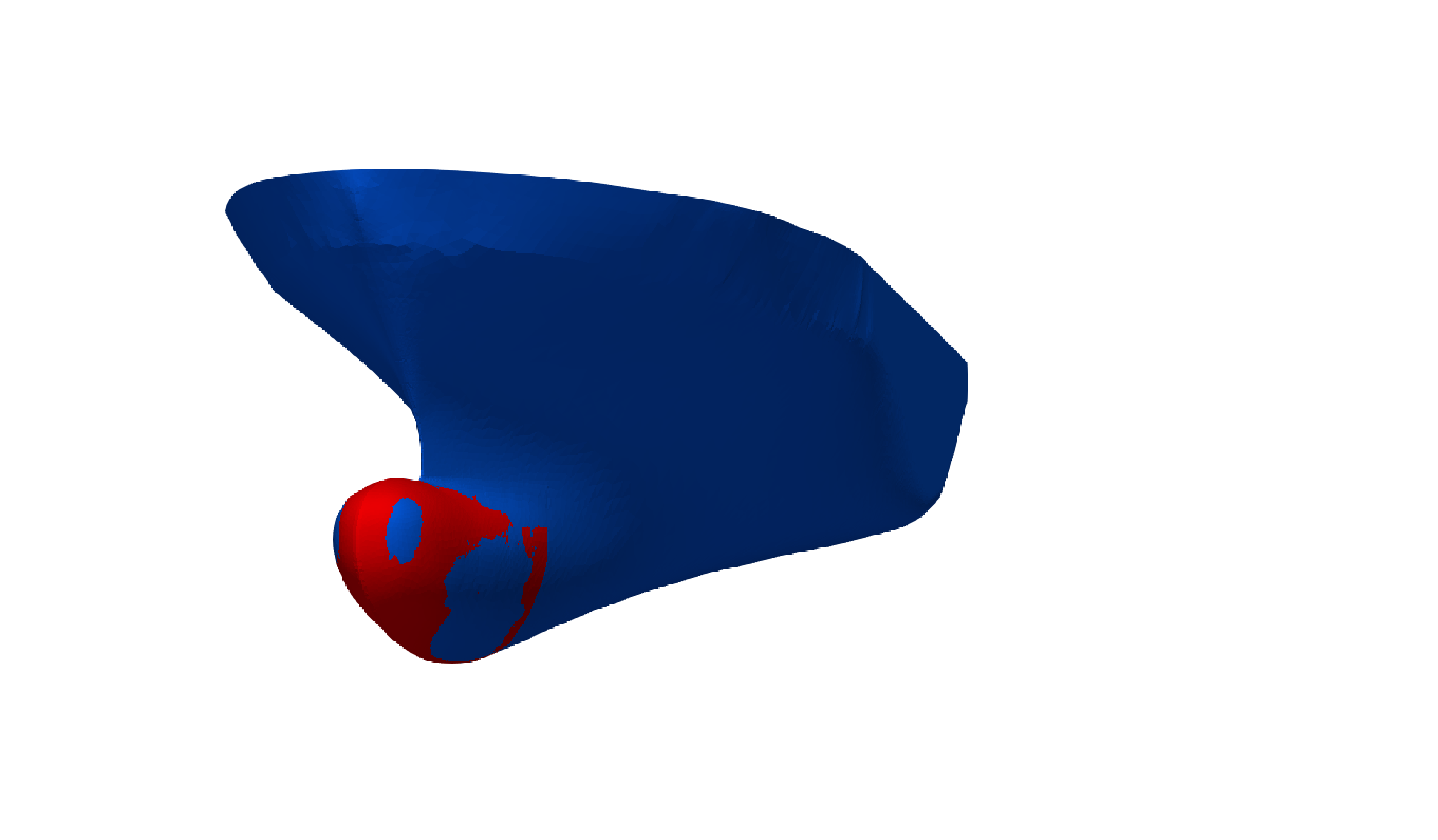}\\
    \includegraphics[width=0.12\textwidth, trim={200 100 500 100}, clip]{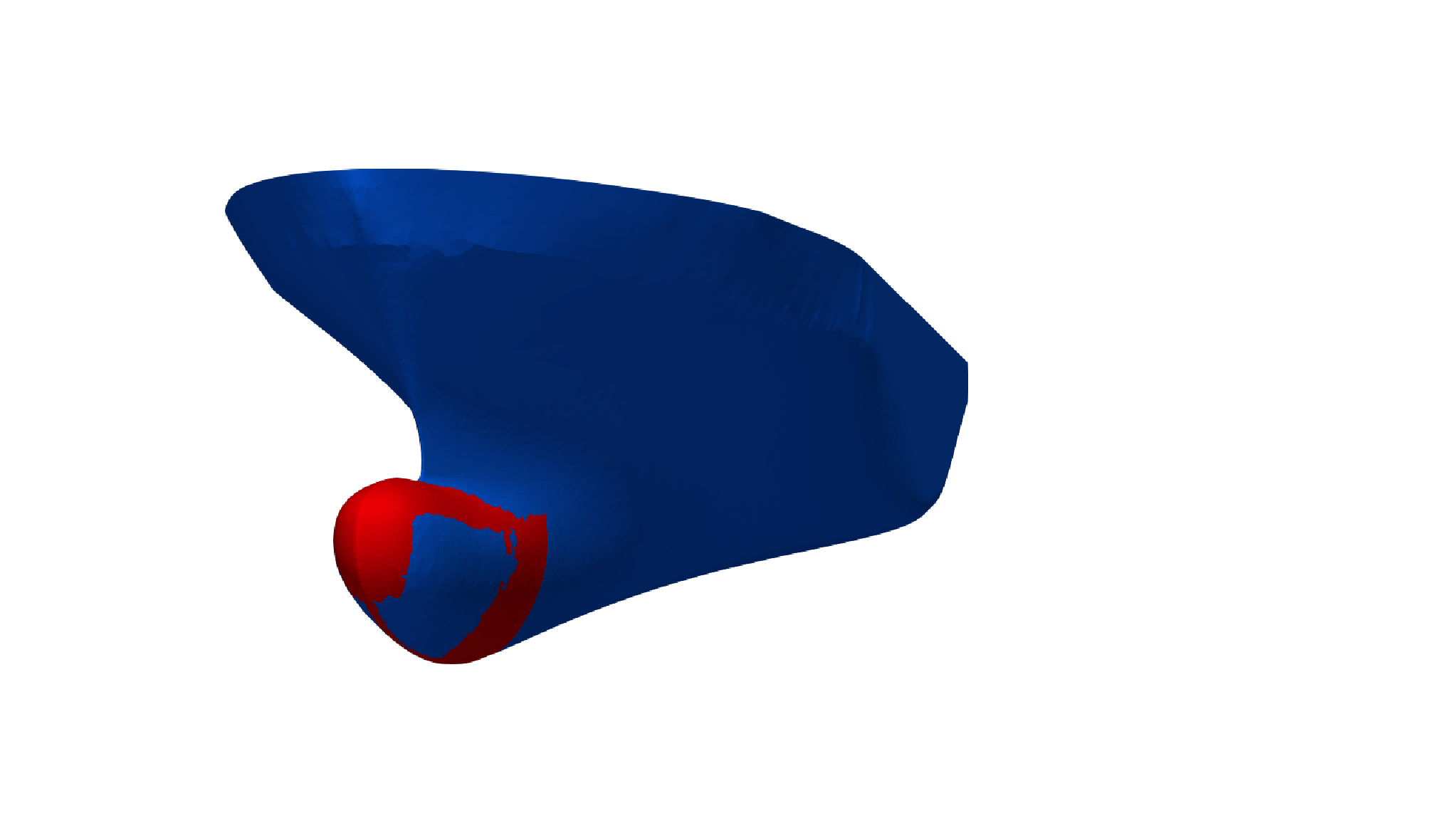}
    \includegraphics[width=0.12\textwidth, trim={200 100 500 100}, clip]{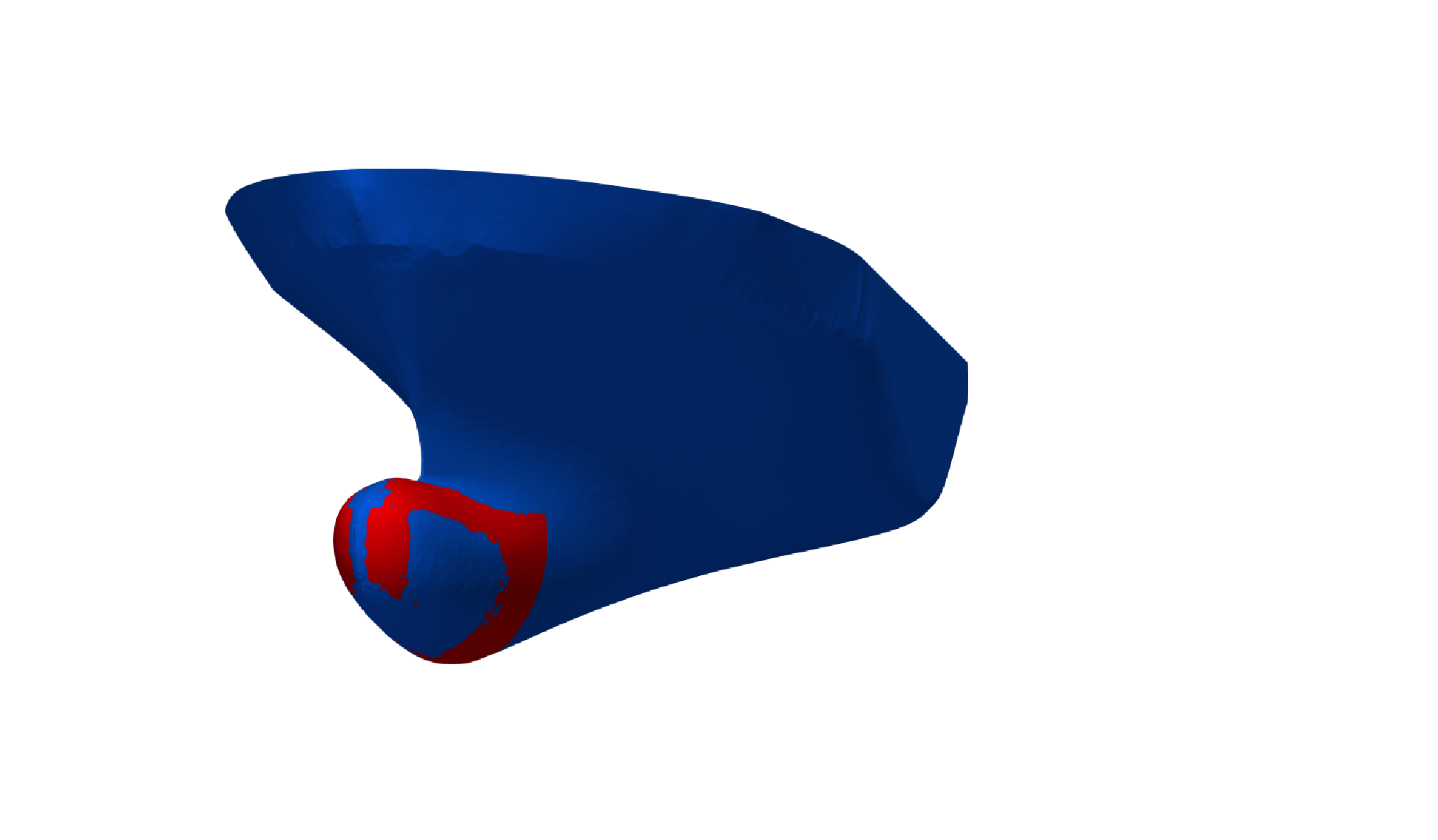}
    \includegraphics[width=0.12\textwidth, trim={200 100 500 100}, clip]{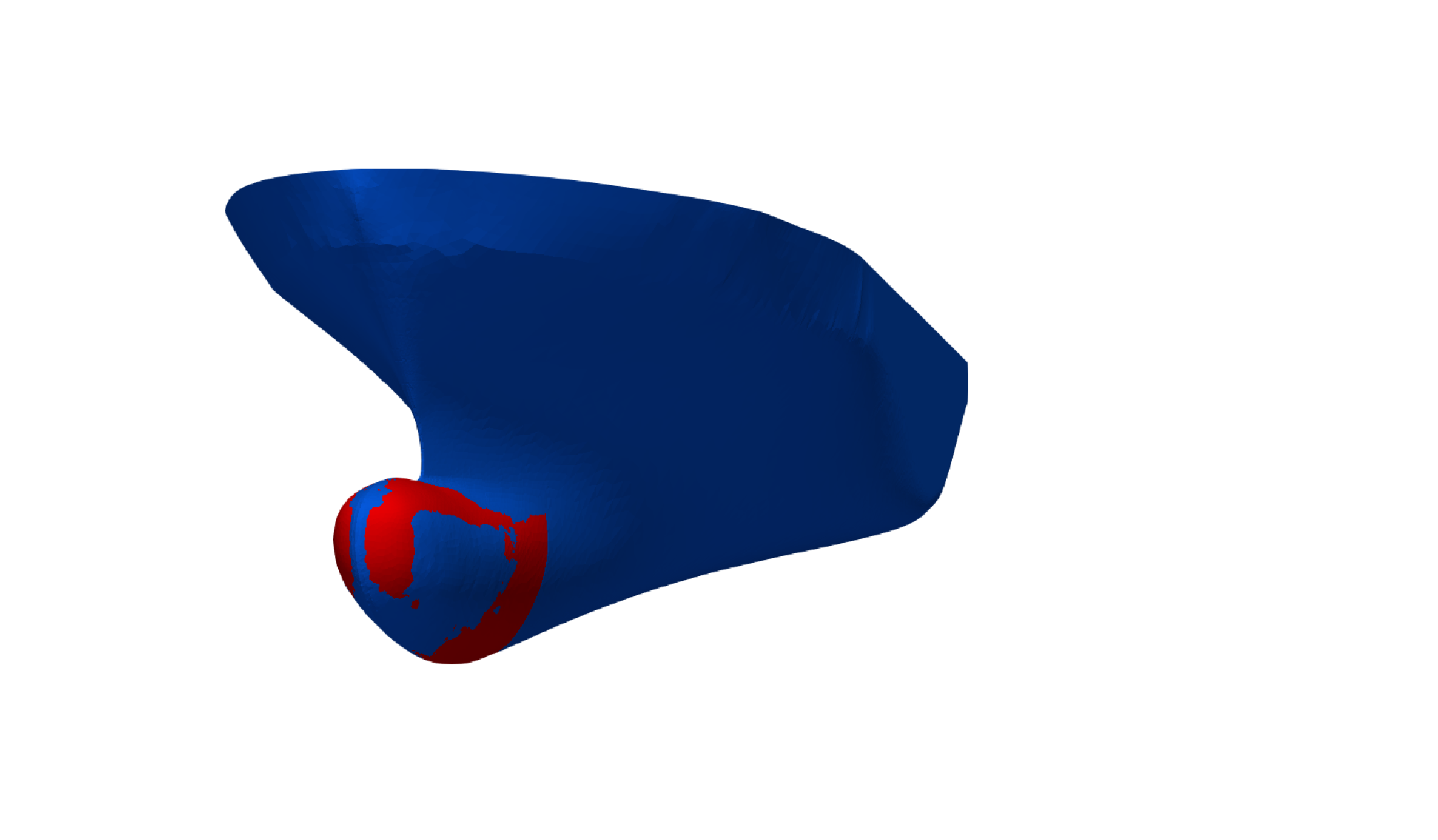}
    \includegraphics[width=0.12\textwidth, trim={200 100 500 100}, clip]{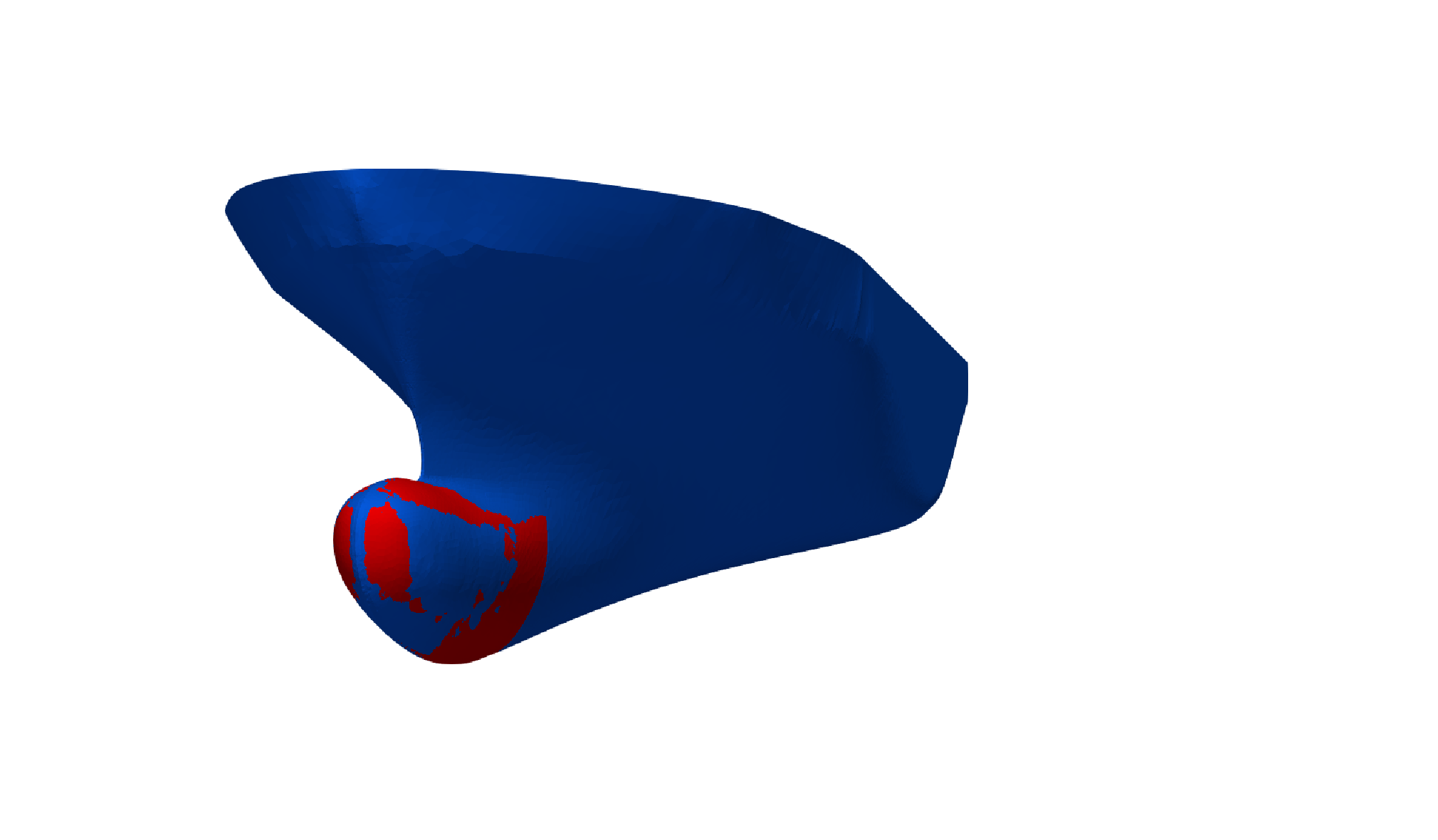}
    \includegraphics[width=0.12\textwidth, trim={200 100 500 100}, clip]{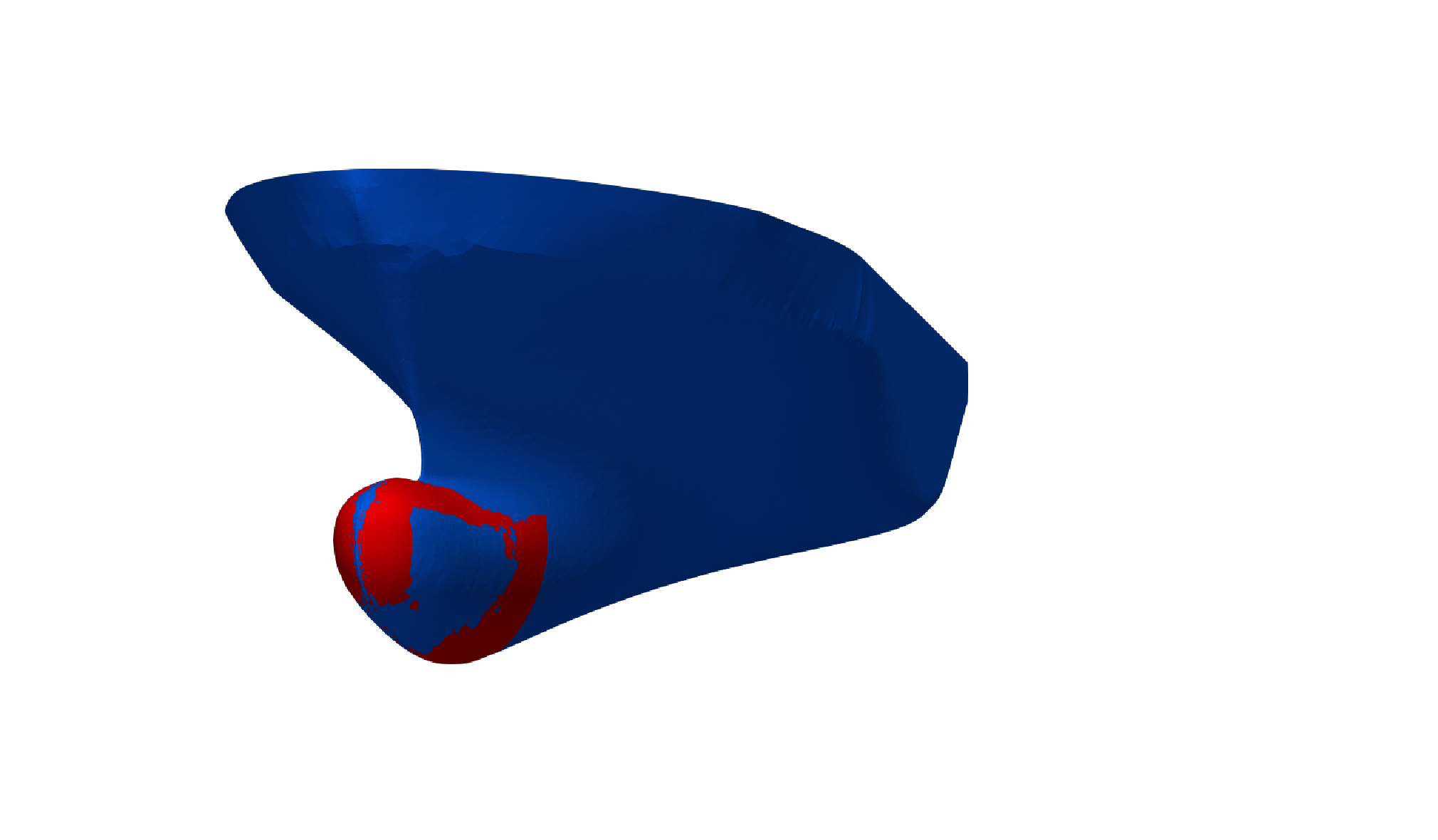}\\
    \includegraphics[width=0.12\textwidth, trim={200 100 500 100}, clip]{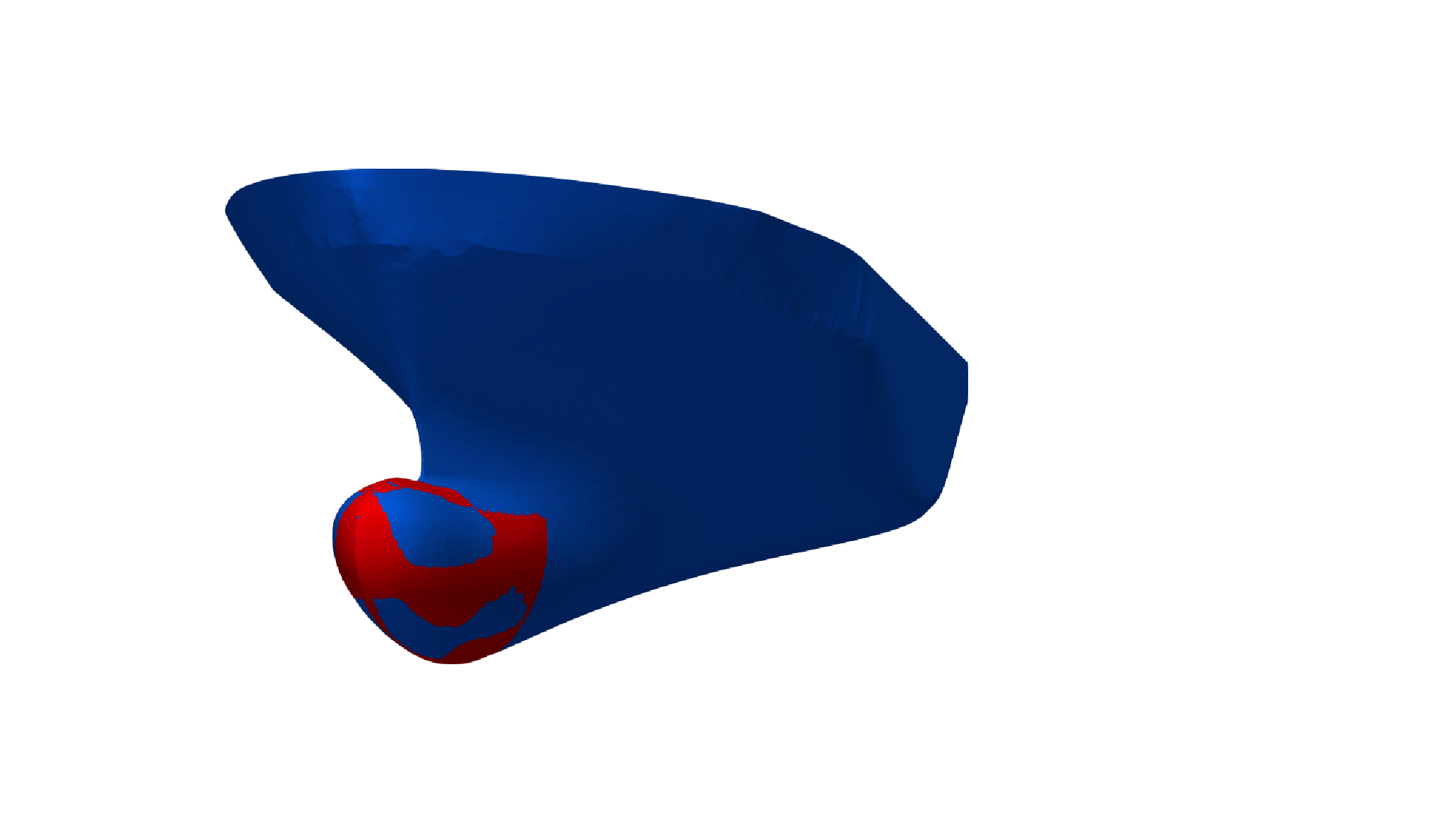}
    \includegraphics[width=0.12\textwidth, trim={200 100 500 100}, clip]{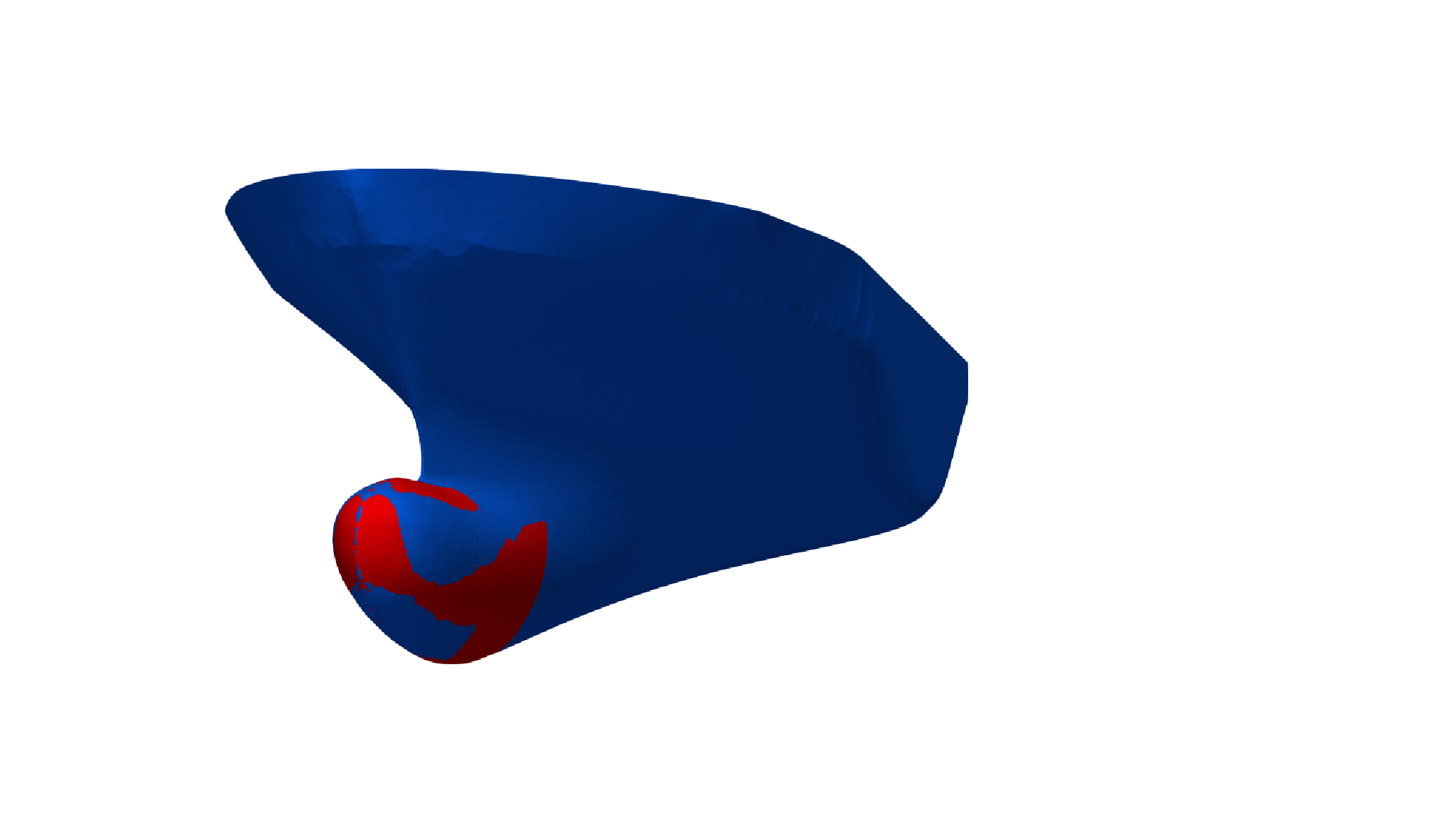}
    \includegraphics[width=0.12\textwidth, trim={200 100 500 100}, clip]{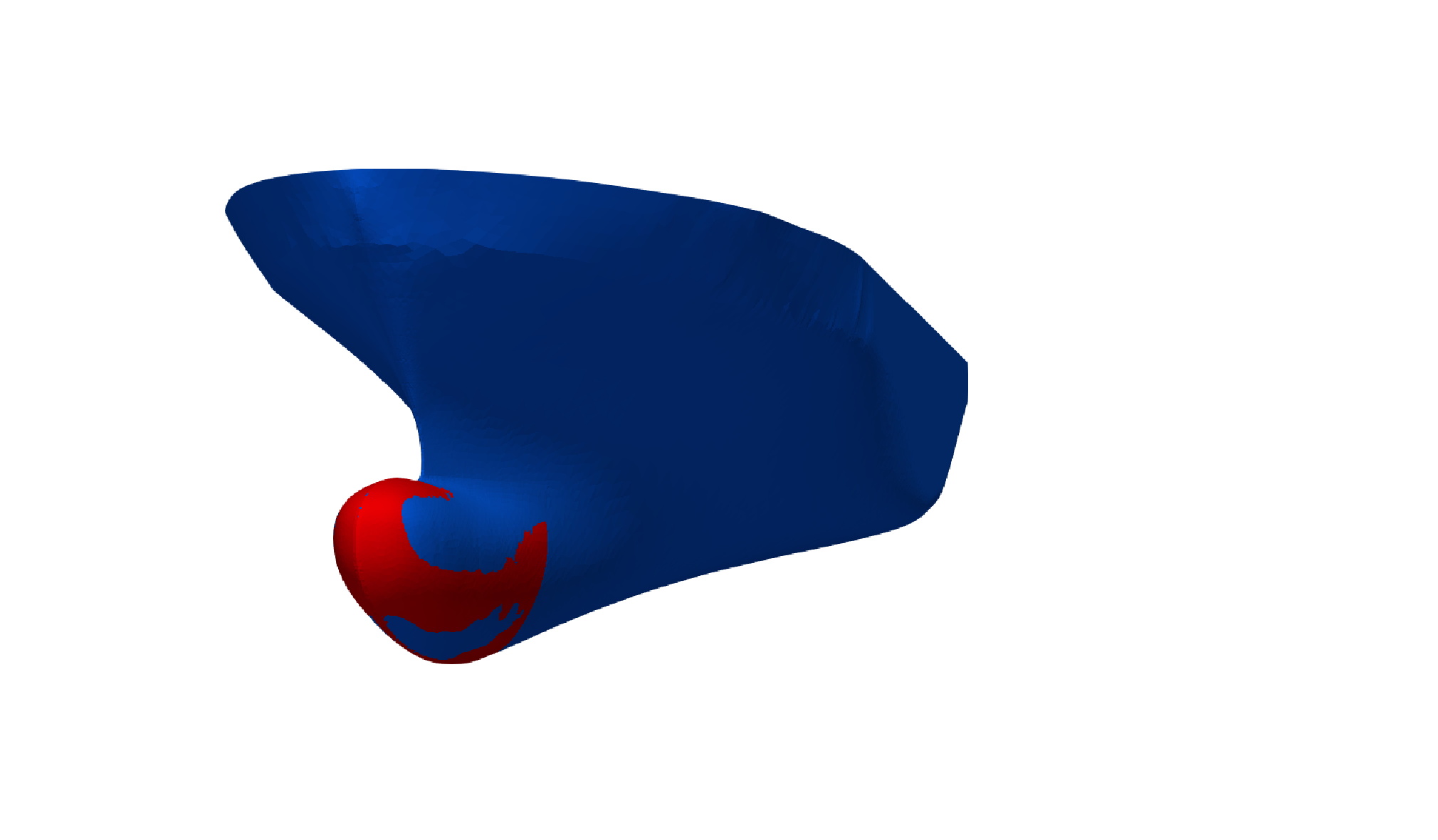}
    \includegraphics[width=0.12\textwidth, trim={200 100 500 100}, clip]{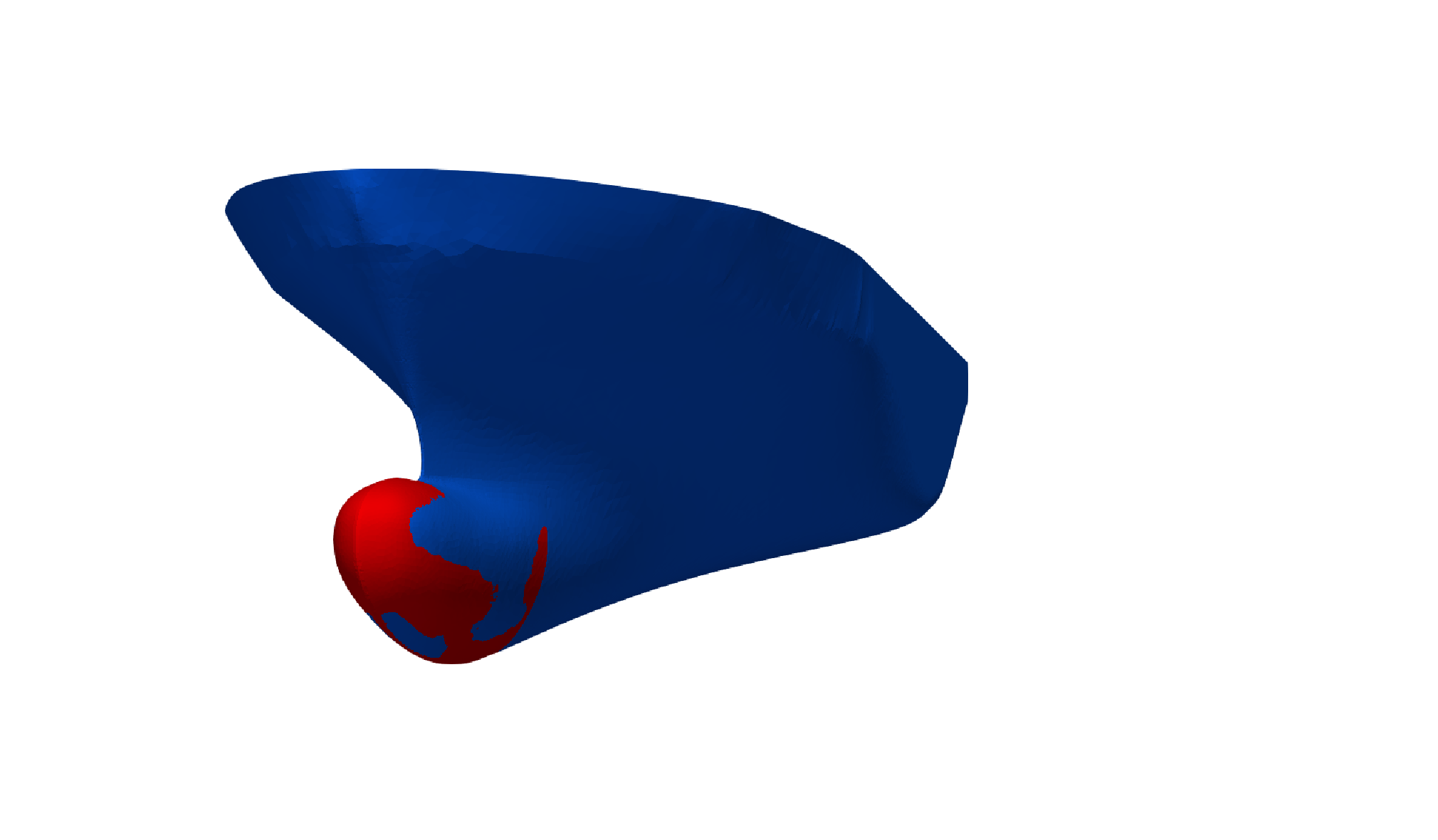}
    \includegraphics[width=0.12\textwidth, trim={200 100 500 100}, clip]{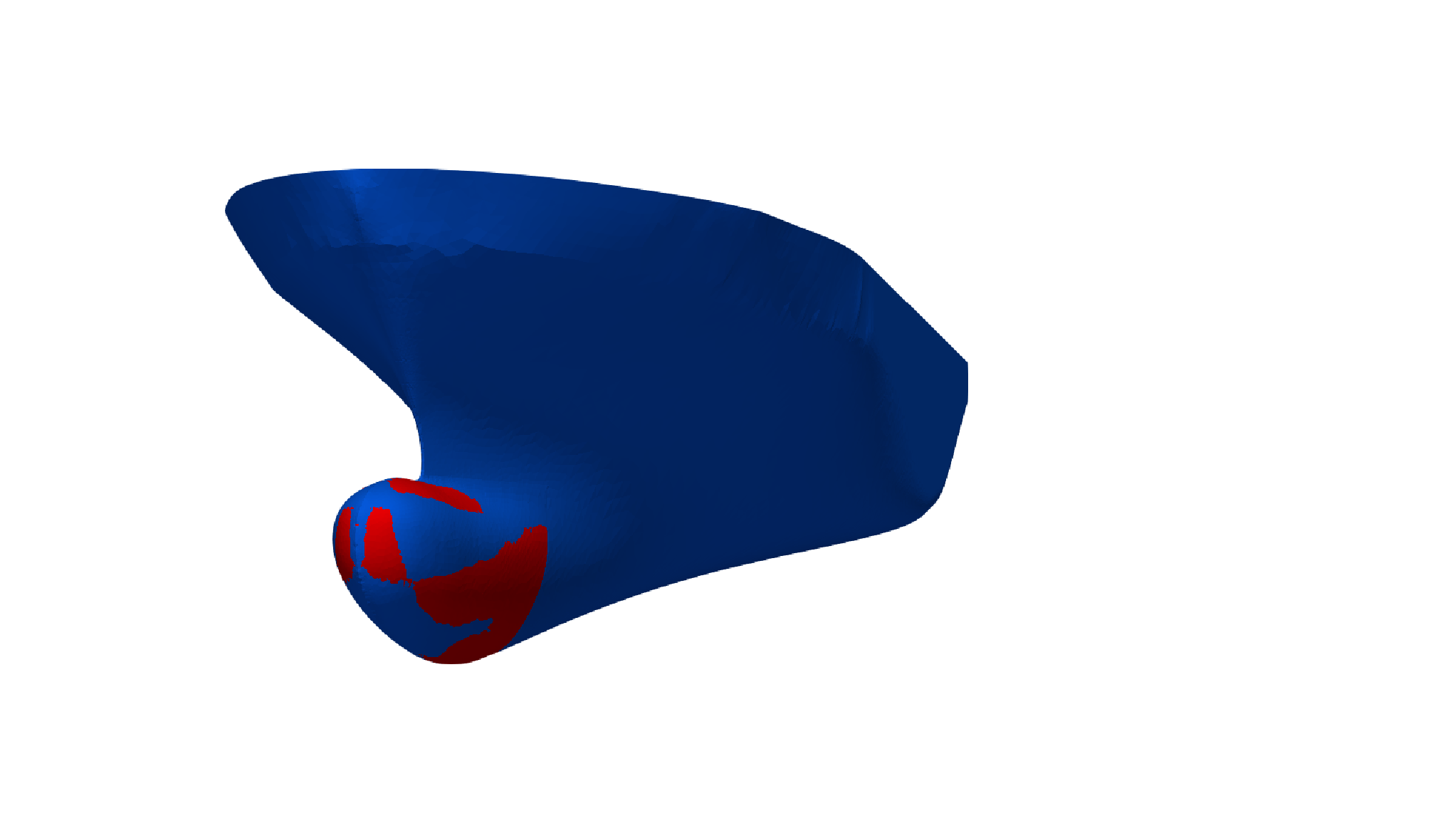}
    \caption{\textbf{HB. }In this figure is shown the overlapping of the newly sampled bulbs from cFFD and cGMs (in red) with respect to the reference STL file used for this test case (in blue). The reference geometry is shown on the top left. Each row corresponds to a different architecture: the first row refers to constrained cFFD of section~\ref{subsec:cffd}, the second to the adversarial autoencoder of paragraph~\ref{par:aae}, the third to the simple autoencoder of paragraph~\ref{par:ae}, the fourth to the beta variational autoencoder of paragraph~\ref{par:vae} and the last to the boundary equilibrium generative adversarial networks of paragraph~\ref{par:began}. All the generative models implement the multilinear constraints enforcing layer of section~\ref{sec:constrained generative models} to preserve the volume.}
    \label{ref:def_HB}
\end{figure}

The number of training and test samples are $\mathbf{n_{\text{train}}=400}$ and $n_{\text{train}}=200$, respectively. For model order reduction, the number of training and test components are instead $\mathbf{n_{_{\text{ROM}},\text{train}}=80}$ and $\mathbf{n_{_{\text{ROM}},\text{test}}=20}$, respectively. We use $\mathbf{r_{\text{PCA}}=140}$ PCA modes for preprocessing inside the cGMs and $\mathbf{r_{\text{POD}}=3}$ modes to perform model order reduction. The parameters' space dimension of cFFD is $\mathbf{p=84}$, while the latent space dimensions of the cGMs is $\mathbf{R=10}$. Some deformations of the cGMs introduced are shown in Figure~\ref{ref:def_HB}: it is shown the overlapping hull's bulbs from the reference STL (in blue) and the STL files generated by cFFD and the other cGMs (in red). The results with respect to the geometrical and physical metrics defined previously are shown in Table~\ref{tab:HB}.

\begin{figure}[htpb!]
    \includegraphics[width=\textwidth]{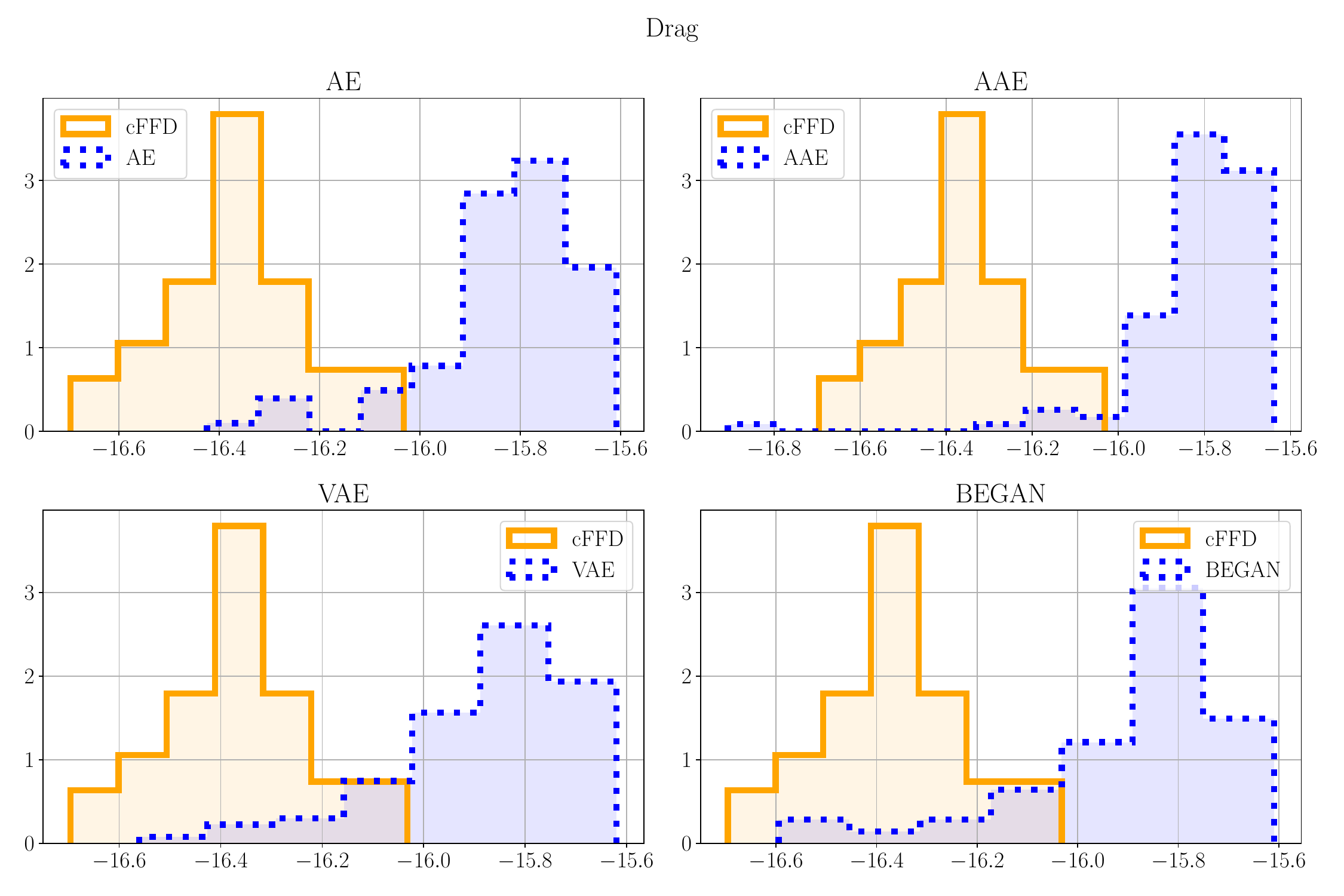}
    \caption{\textbf{HB. }In this figure the histograms of the drag coefficient from cFFD and from the cGMs are shown. It can be seen that the drag coefficient is very sensitive to the deformations of the bulb and that our generated geometries are probably biased towards higher values of the drag. For a quantitative measure see Table~\ref{tab:HB}. The histograms are obtained from the $\mathbf{n_{\text{ROM},\text{test}}=20}$ test samples.}
    \label{fig:hist_cd}
\end{figure}

\begin{figure}[htpb!]
    \includegraphics[width=\textwidth]{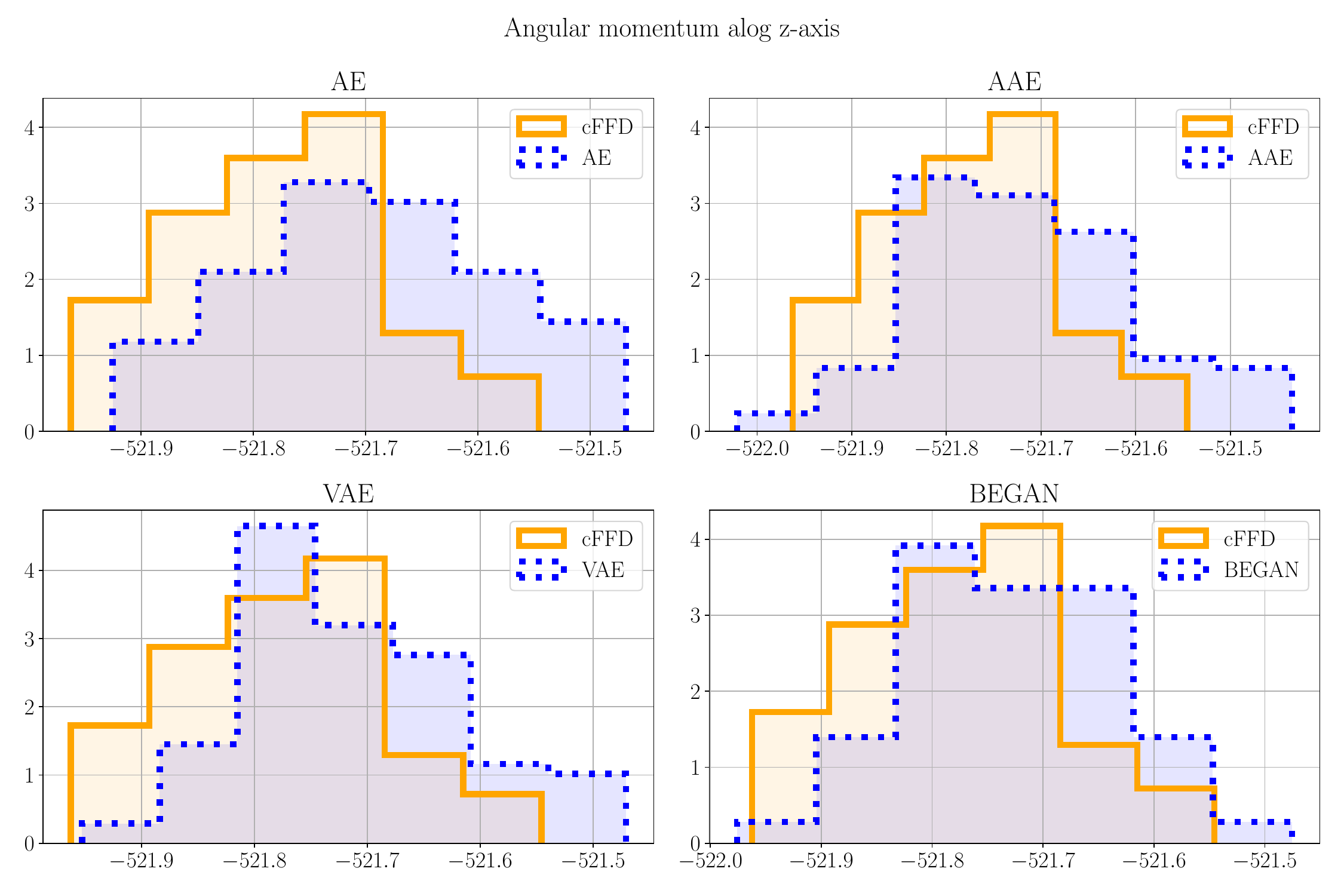} \caption{\textbf{HB. }In this figure the histograms of the angular momentum along the z-axis $M_{z}$ from the mesh obtained from cFFD or cGMs are shown. The histogram area intersection is more than $50\%$. For a quantitative measure see Table~\ref{tab:HB}. The histograms are obtained from the $\mathbf{n_{\text{ROM},\text{test}}=20}$ test samples.}
    \label{fig:hist_M}
\end{figure}

Qualitatively the histograms of each architecture are reported in Figure~\ref{fig:hist_M} for the $M_z$ angular momentum along the z-axis and in Figure~\ref{fig:hist_cd} for the $c_d$ drag coefficient. There is an evident bias of the cGMs distributions with respect to the cFFD's one when considering the drag coefficient. Possible reasons for this bias are the sharp edges introduced with the cFFD deformations on the gluing sites of the bulb on the hull and the coarse computational mesh employed for the numerical simulations with respect to the resolution of the STL files used to train the cGMs.

\begin{figure}[htpb!]
    \includegraphics[width=\textwidth]{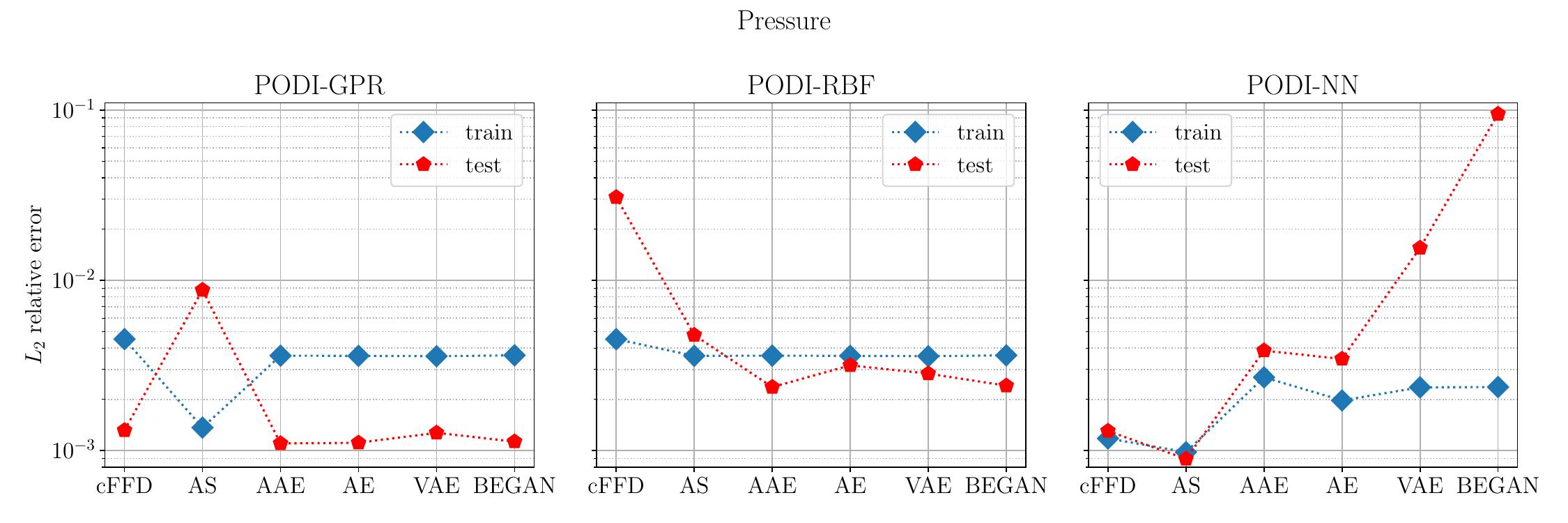}
    \caption{\textbf{HB. }Here we show the ROM performance on the pressure on the whole domain over various datasets, using different interpolation and regression techniques combined with proper orthogonal decomposition with interpolation (PODI): Gaussian process regression (GPR), radial basis functions interpolation (RBF) and feed-forward neural networks (NN). For every method, there is at least a cGM that performs slightly better than the cFFD in terms of accuracy, while reducing the parameters' space dimension from $\mathbf{p=84}$ to $\mathbf{R=10}$. The active subspace (AS) dimension chosen is $\mathbf{r_{\text{AS}}=1}$. The training and test errors are evaluated on $\mathbf{n_{\text{ROM},\text{test}}=80}$ and $\mathbf{n_{\text{ROM},\text{test}}=20}$ training and test samples.}
    \label{fig:HB_rom_p}
\end{figure}

\begin{figure}[htpb!]
    \includegraphics[width=\textwidth]{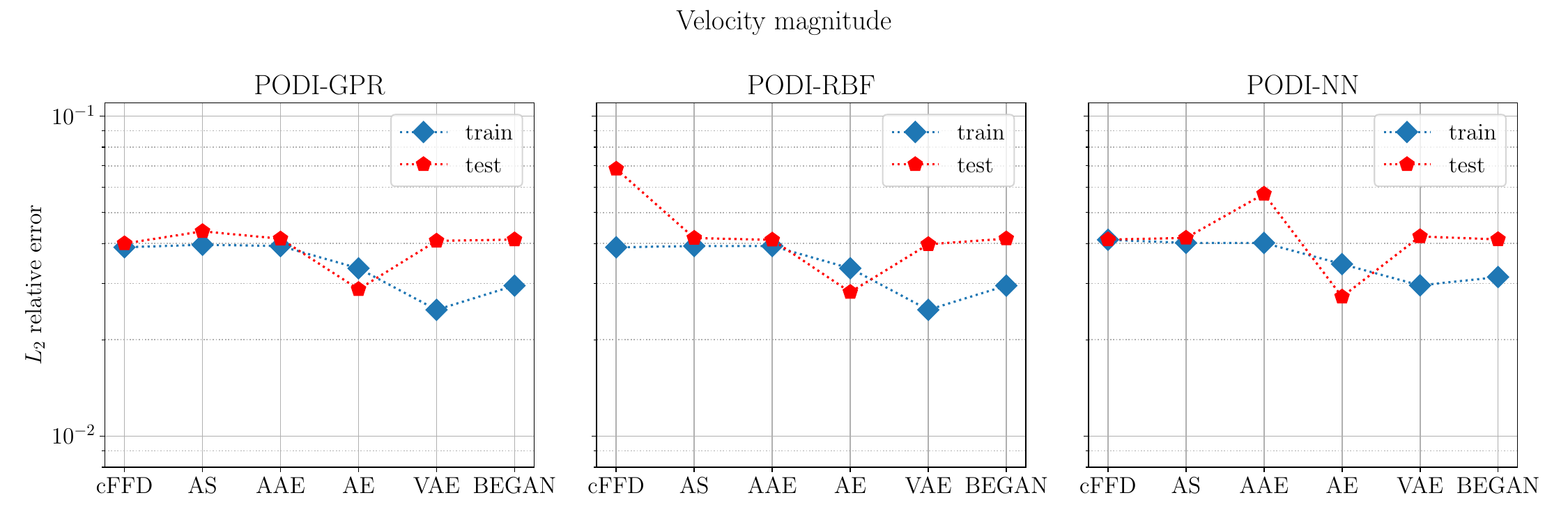}
    \caption{\textbf{HB. }Here we show the ROM performance on the velocity magnitude on the whole domain over various datasets, using different interpolation and regression techniques combined with proper orthogonal decomposition with interpolation (PODI): Gaussian process regression (GPR), radial basis functions interpolation (RBF) and feed-forward neural networks (NN). For every method, there is at least a cGM that performs slightly better than the cFFD in terms of accuracy, while reducing the parameters' space dimension from $\mathbf{p=84}$ to $\mathbf{R=10}$. The active subspace (AS) dimension chosen is $\mathbf{r_{\text{AS}}=1}$. The training and test errors are evaluated on $\mathbf{n_{\text{ROM},\text{test}}=80}$ and $\mathbf{n_{\text{ROM},\text{test}}=20}$ training and test samples.}
    \label{fig:HB_rom_u}
\end{figure}

\begin{figure}[htpb!]
    \centering
    \includegraphics[width=0.4\textwidth]{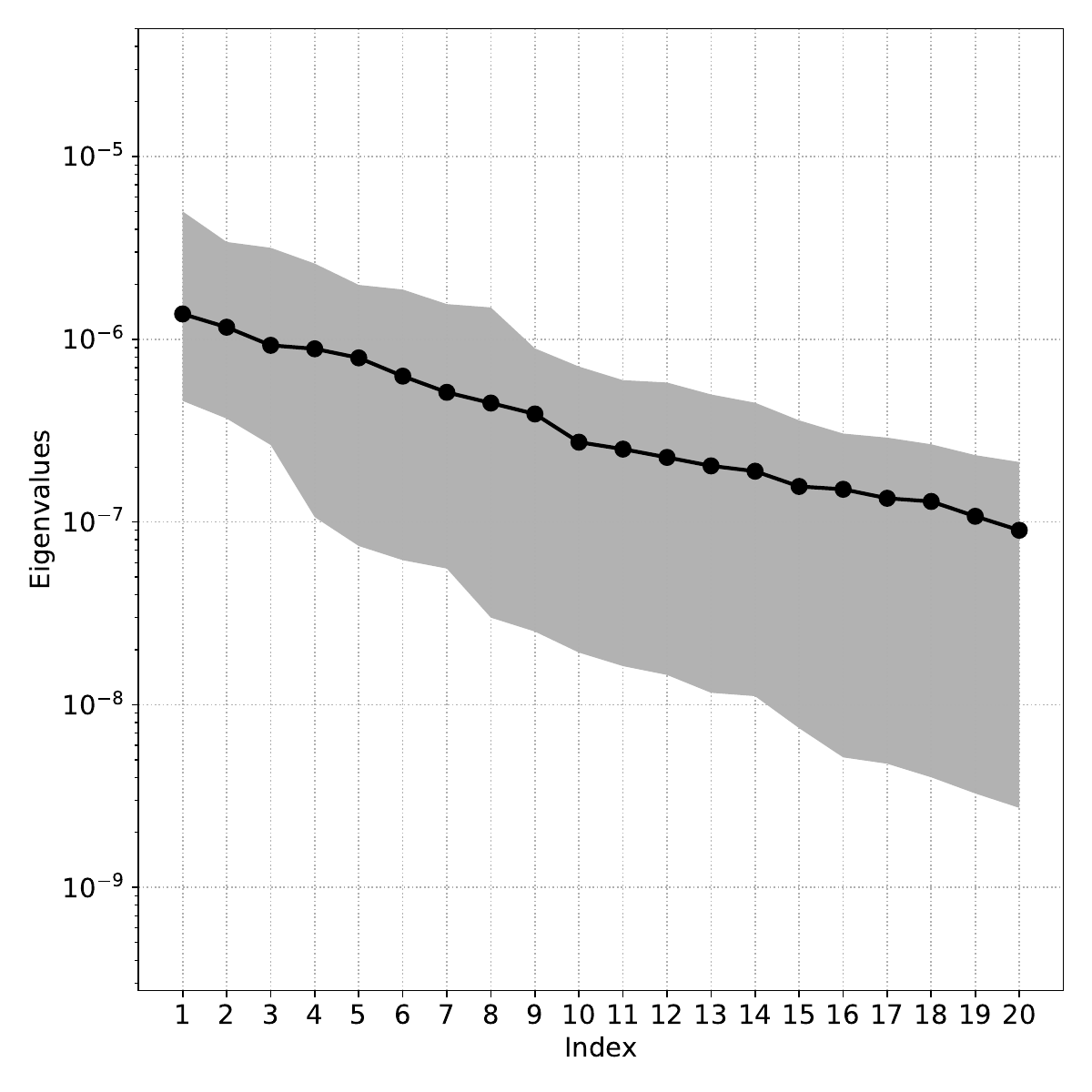}
    \includegraphics[width=0.4\textwidth]{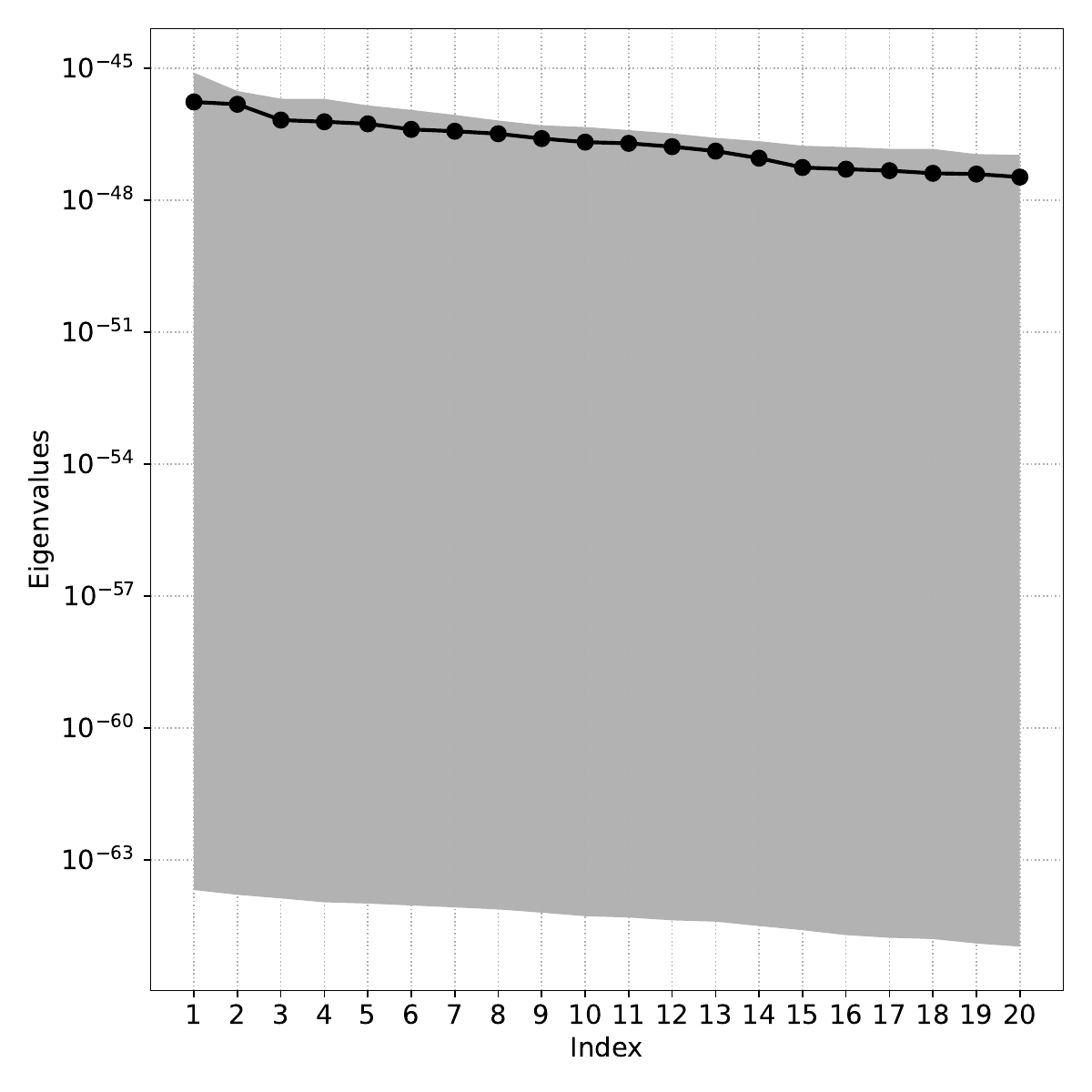}
    \caption{\textbf{SB and HB.} In this figure we show the first $20$ eigenvalues of the Active Subspace method for the \textbf{SB} test case on the left and the \textbf{HB} test case on the right. The shaded grey areas are delimited by the minimum and maximum values of the $100$ bootstrap replicates. The black line corresponds to the average value. The spectral gaps suggest an AS dimension equal to $r_{\text{AS}}=9$ and $r_{\text{AS}}=2$, for the \textbf{SB} and \textbf{HB} test cases, respectively.}
    \label{fig:AS_evals}
\end{figure}

The non-intrusive reduced order models implemented with GPR-PODI, RBF-PODI and NN-PODI are shown in Figures~\ref{fig:HB_rom_p} for the effective pressure and~\ref{fig:HB_rom_u} for the velocity magnitude on the whole computational domain $\Omega$. The original parameters space's dimension of cFFD changes from $\mathbf{p=84}$ to $\mathbf{R=10}$ for the cGMs, while keeping more or less the same accuracy on the surrogate models for the pressure and velocity fields. We apply a further level of reduction in the space of parameters with the AS method on the latent space of $\mathbf{R=10}$ coordinates of the cGMs: the new parameter space is a one-dimensional active subspace for each variable. The accuracy of the AS response surfaces is almost the same for the pressure and the velocity magnitude while reducing the parameter space's dimension from $\mathbf{R=10}$ to $\mathbf{r_{\text{AS}}=1}$. For simplicity, we show the AS response surface with dimension $\mathbf{r_{\text{AS}}=1}$, even if from the plot in Figure~\ref{fig:AS_evals} of the first $20$ eigenvalues of the uncentered covariance matrix from equation~\eqref{eq:covSA}, the spectral gap~\cite{constantine2015active} suggests $\mathbf{r_{\text{AS}}=2}$.

The speedup for the generation of the geometries employing the cGMs instead of cFFD is around $360$. The speedup of PODI non-intrusive model order reduction with respect to the full-order simulations is around $432000$.

\begin{figure}
    \centering
    \includegraphics[width=0.20\textwidth, trim={300 50 450 50}, clip]{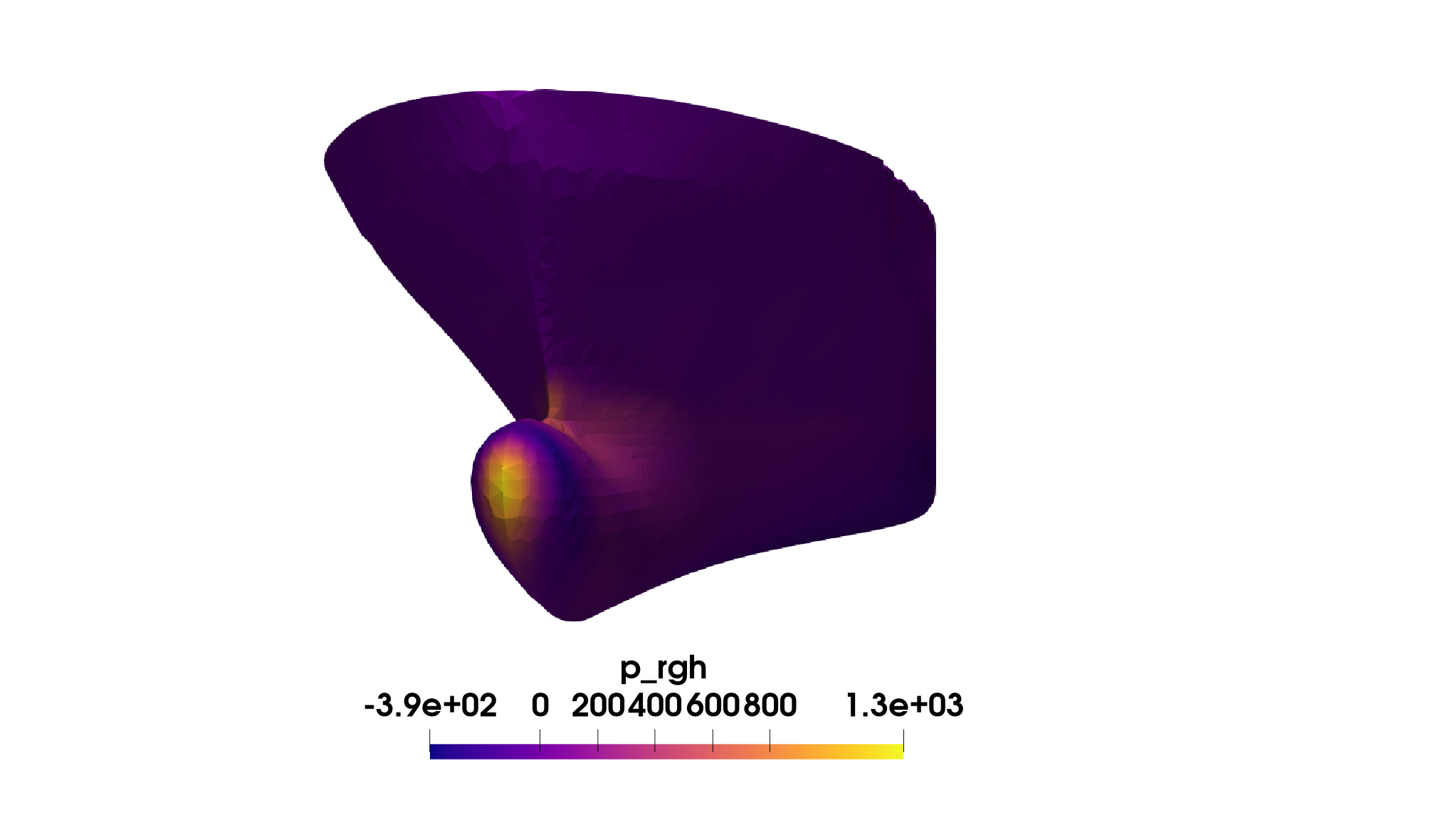}
    \includegraphics[width=0.20\textwidth, trim={300 50 450 50}, clip]{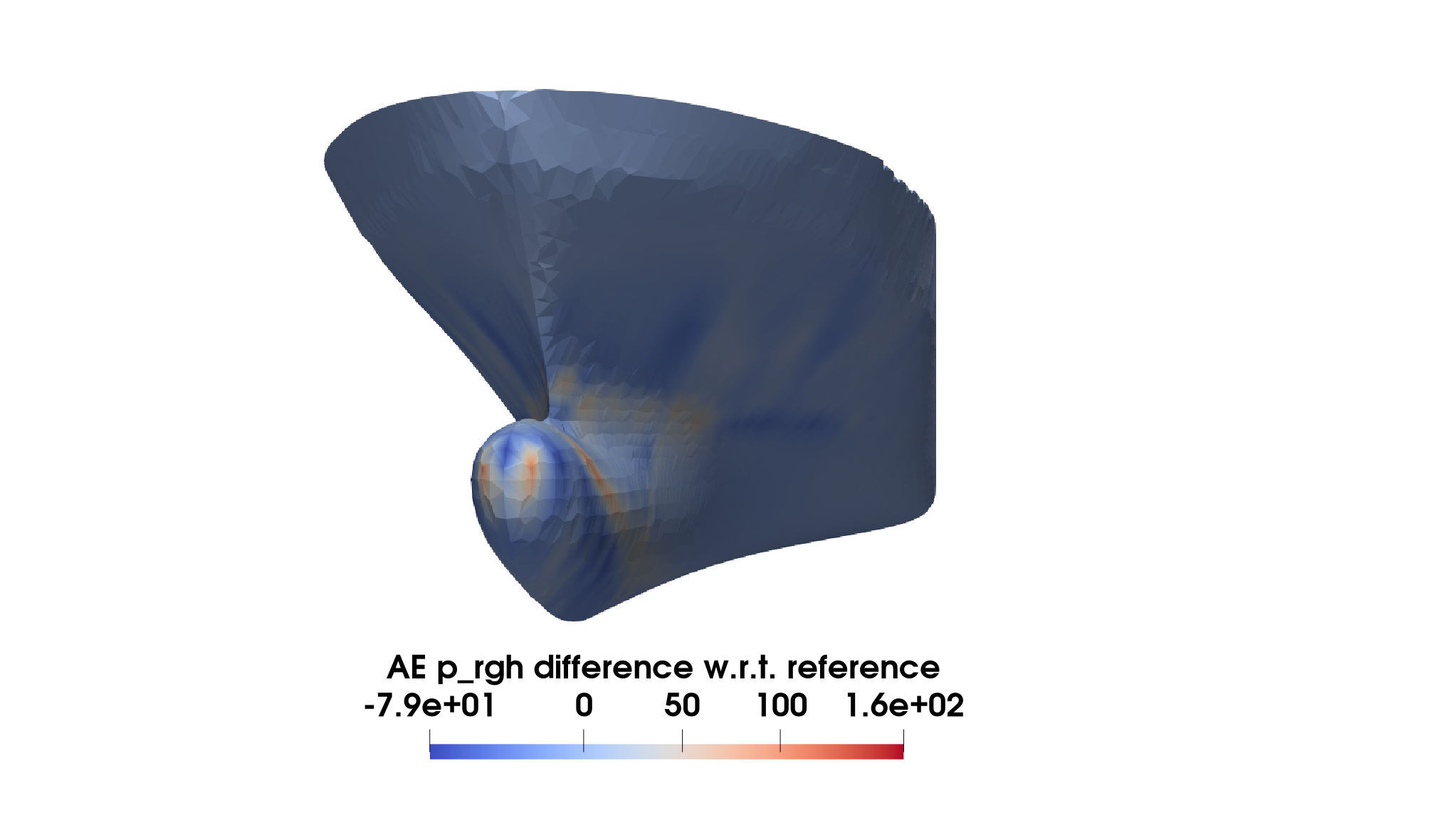}
    \includegraphics[width=0.20\textwidth, trim={300 50 450 50}, clip]{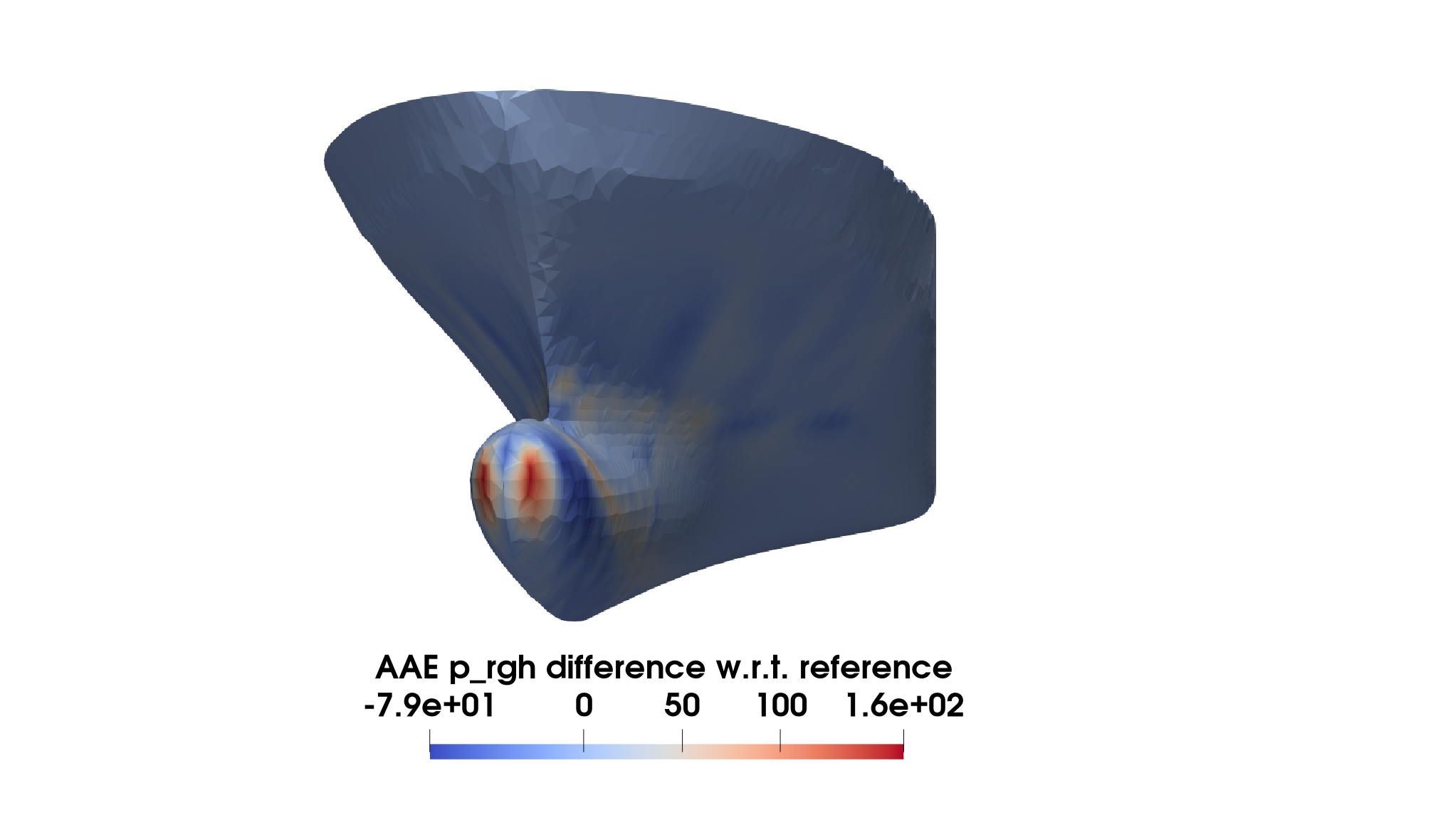}\\
    \includegraphics[width=0.20\textwidth, trim={300 50 450 50}, clip]{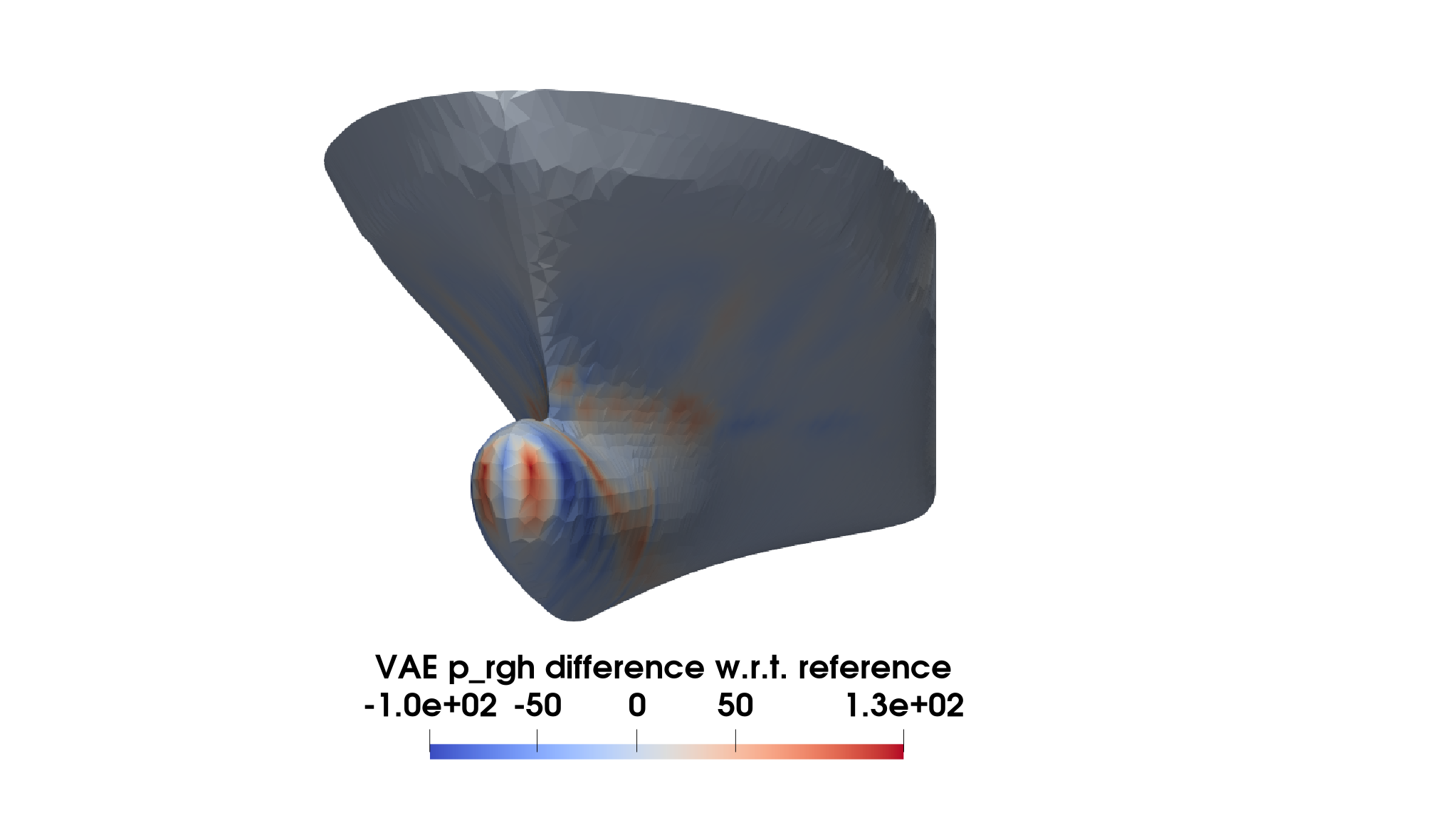}
    \includegraphics[width=0.20\textwidth, trim={300 50 450 50}, clip]{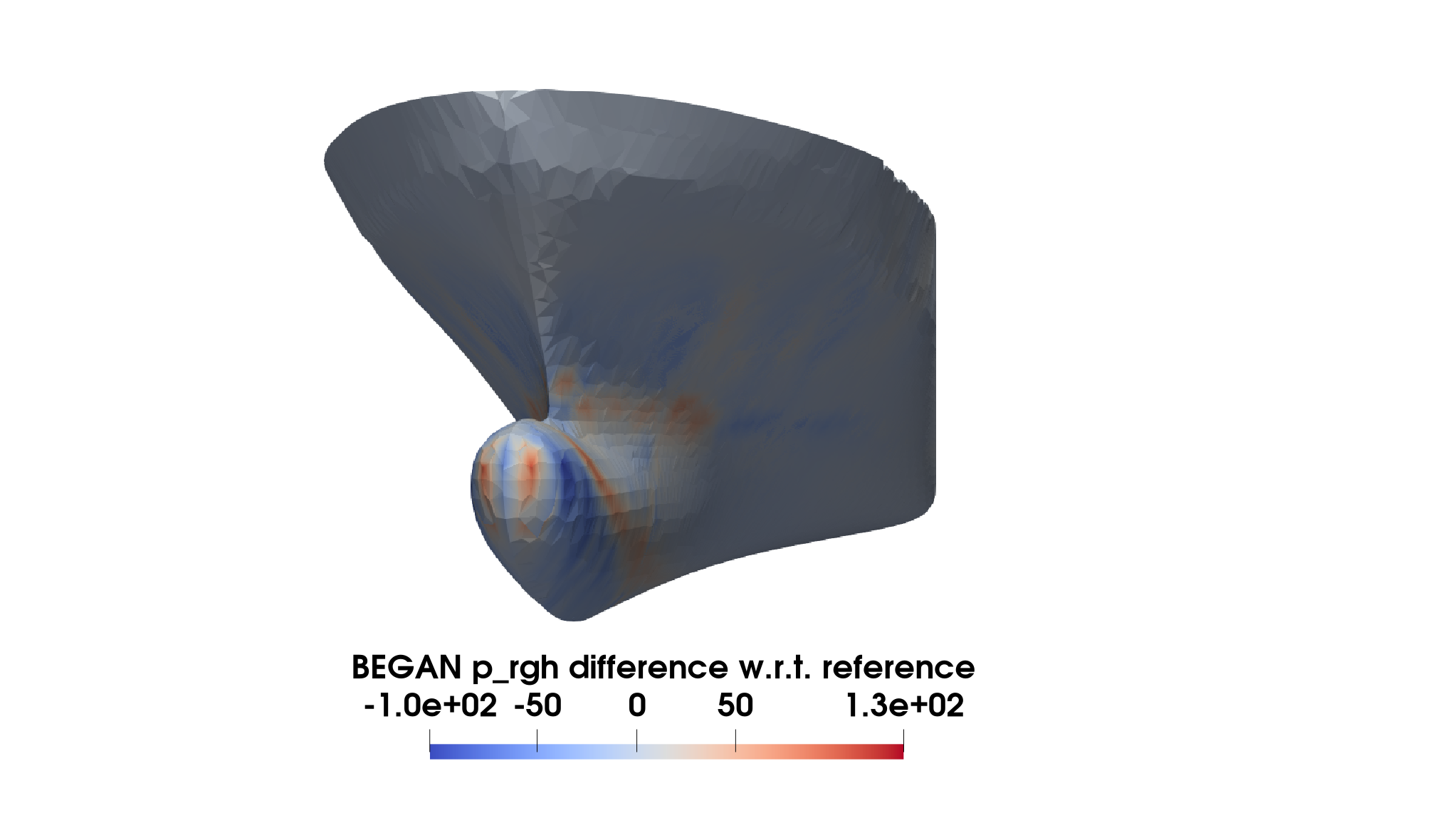}
    \caption{\textbf{HB. }Difference of the effective pressure field on the reference computational mesh with the effective pressure fields on different geometries sampled from constrained generative models, namely from top left to bottom right from the simple autoencoder, from the adversarial autoencoder, from the beta variational autoencoder and from the boundary equilibrium generative adversarial network. The geometrical constraint preserved is the bulb's and thus the hull's volume.}
\end{figure}

\begin{table}[htpb!]
    \centering
    \caption{\textbf{HB. }In this table we show the evaluation metrics of the Duisburg test case defined in equations~\eqref{eq:HB_met}:$I_{xx},I_{xy},I_{xz},I_{yy},I_{yz},I_{zz}$ are the components of the inertia tensor, $A_{\text{hull}}$ area of the surface of the bulb, $c_d$ is the drag coefficient evaluated on the hull, $M_z$ is the angular momentum along the z-axis. The metrics $I_{xx},I_{xy},I_{xz},I_{yy},I_{yz},I_{zz},A_{\text{bulb}}$ are obtained from $\mathbf{n_{\text{test}}=200}$ test samples, while the metrics $c_d, M_z$ are based on  $n_{\text{ROM},\text{test}}=20$ test samples. The model that has overall the best performance is the AE, since the BEGAN produces too many out-of-training distribution samples as can be seen comparing the values of the total variance: the variance of the training data from cFFD is \textbf{1.2e-03}.}
    \begin{tabular}{l|l|l|l|l}
        \hline
        \hline
        &   AE &  AAE &  VAE &  BEGAN   \\
        \hline
        \hline
        $JSD(A_{\text{bulb}})$  & \textbf{1.2e-01} & 2.3e-01 & 3.6e-01  & 5.2e-01  \\ \hline
        $JSD(I_{xx})$ & \textbf{1.6e-01} & 3.0e-01 & 3.8e-01 & 6.0e-01  \\ \hline
        $JSD(I_{xy})$ & 5.5e-01 & 6.9e-01 & \textbf{4.7e-01}  & {6.7e-01}\\ \hline
        $JSD(I_{xz})$ & 2.1e-01 & 2.1e-01 & 5.0e-01   & \textbf{2.0e-01}\\ \hline
        $JSD(I_{yy})$ & \textbf{2.4e-01} & 4.9e-01 & 8.8e-01 & 9.1e-01\\ \hline
        $JSD(I_{yz})$ & 8.7e-01& {5.6e-01} & 6.6e-01  & \textbf{4.8e-01}\\ \hline
        $JSD(I_{zz})$ & \textbf{1.5e-01} & 4.5e-10 & 5.6e-1 & 3.2e-01\\ \hline
        $JSD(c_{d})$ & 9.1e-01 & 9.4e-01 & 8.5e-01 & \textbf{8.1e-01} \\ \hline
        $JSD(M_{z})$ & 3.4e-01 & 2.9e-01 & 2.9e-01 & \textbf{2.6e-01} \\ \hline
        $Var$ & 1.2e-03 & 2.4e-03 & 4e-04 & 7.0e-03 \\
    \end{tabular}
    \label{tab:HB}
\end{table}
\begin{figure}[htpb!]
    \centering
    \includegraphics[width=0.49\textwidth]{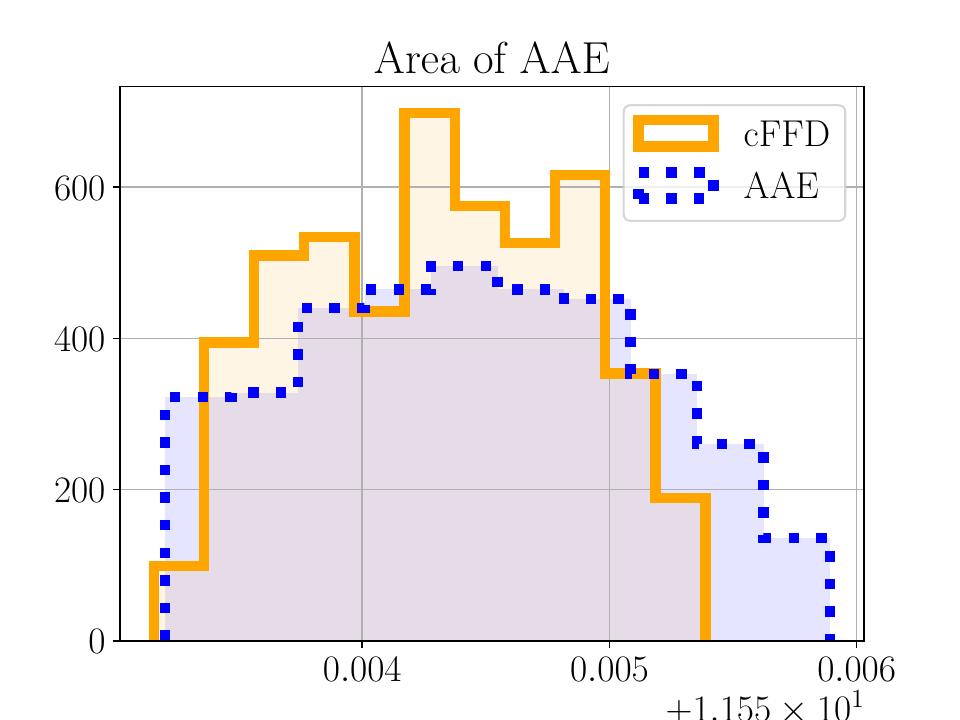}
    \includegraphics[width=0.49\textwidth]{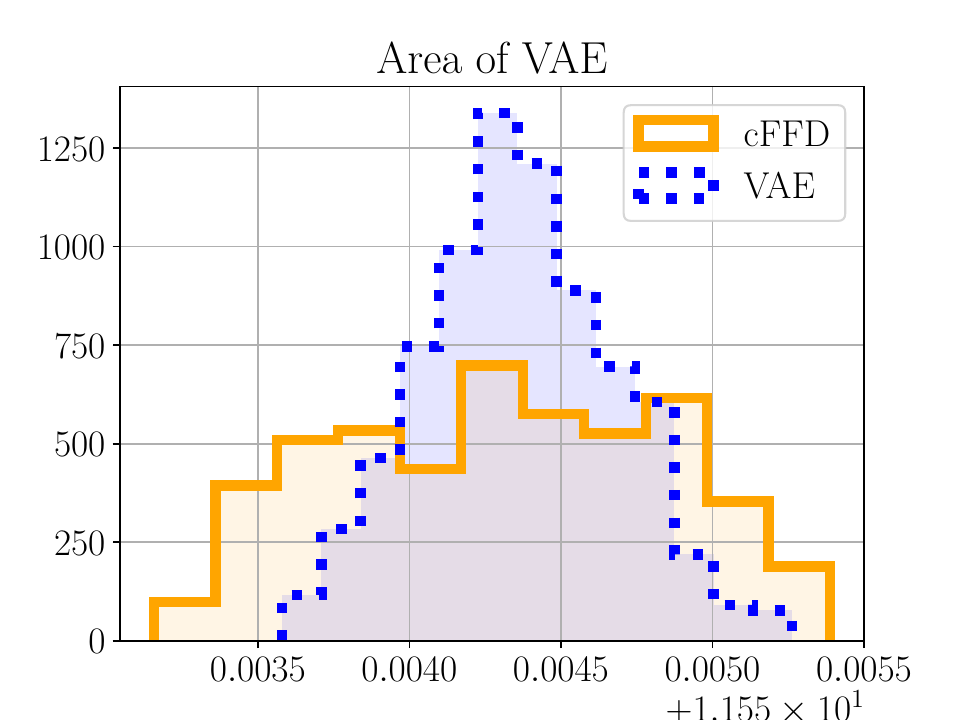}\\
    \includegraphics[width=0.49\textwidth]{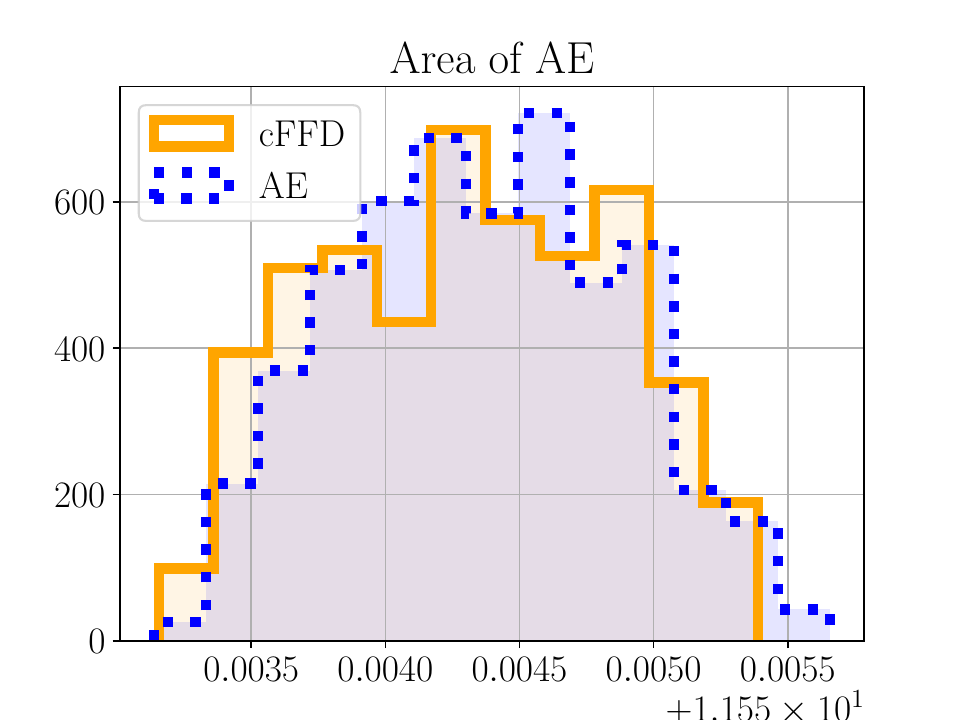}
    \includegraphics[width=0.49\textwidth]{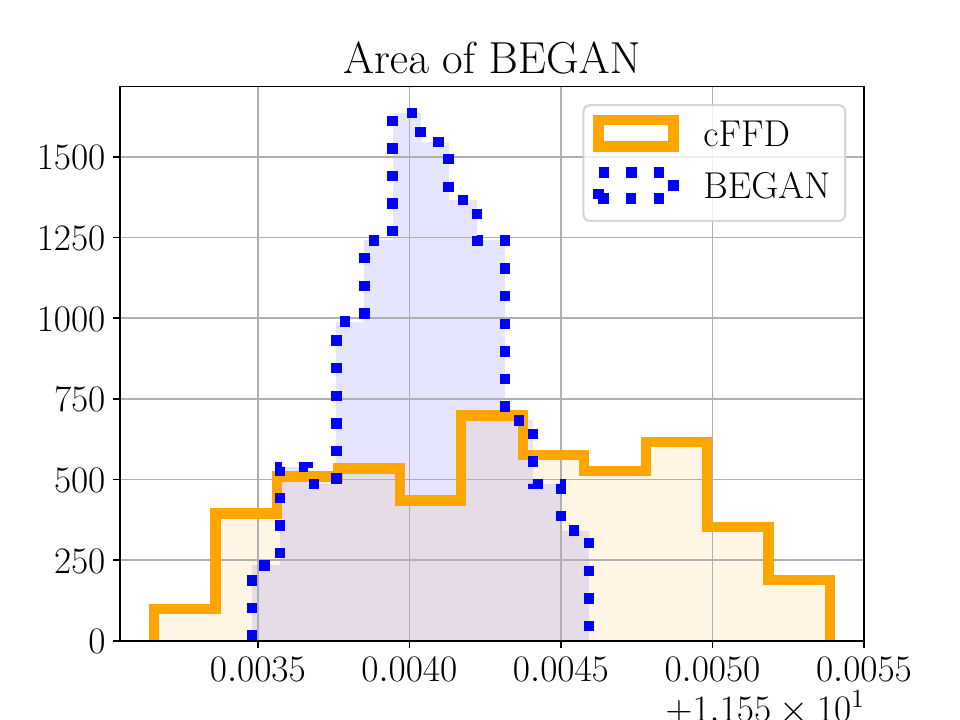}
    \caption{\textbf{HB. }In this figure the histograms of the area of the bulb's deformed with cFFD or with cGMs are shown. The histogram area intersection is more than $50\%$. For a quantitative measure see Table~\ref{tab:HB}. The histograms are obtained from the $\mathbf{n_{\text{test}}=200}$ test samples.}
\end{figure}

\begin{figure}[htpb!]
    \includegraphics[width=\textwidth]{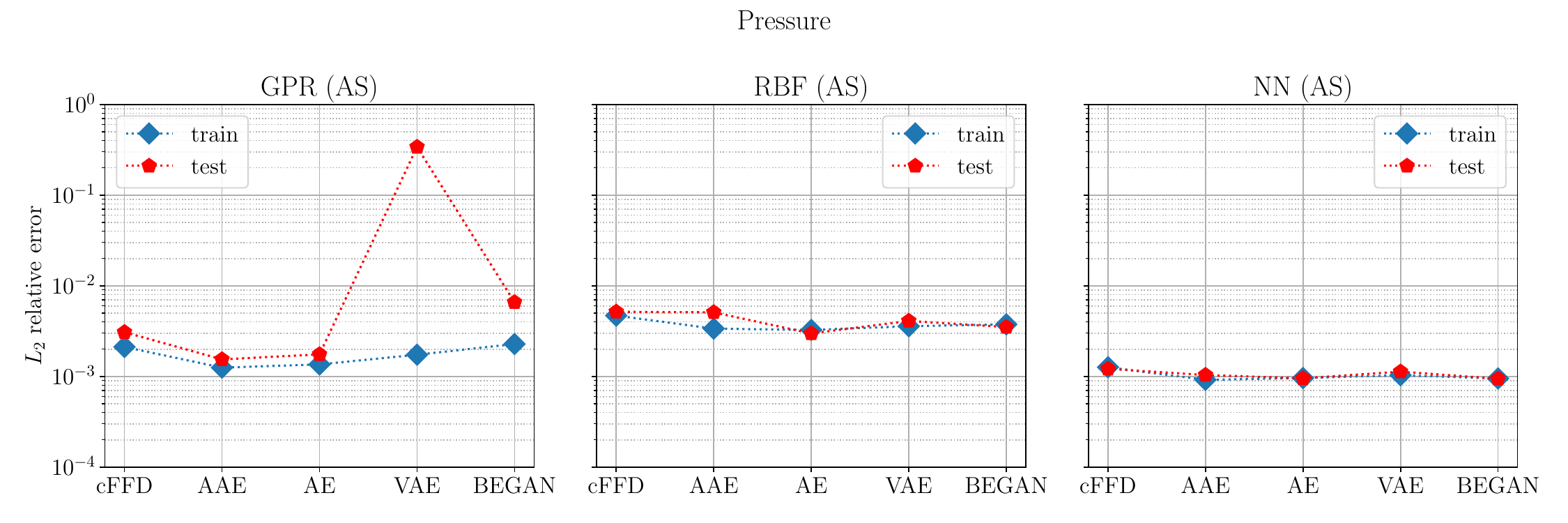}
    \caption{\textbf{HB. }In this figure we show the ROM-PODI performance, for the effective pressure on the whole domain, on the training and test datasets coupled with AS dimension reduction in the space of parameters. It can be seen that the latent dimension can be reduced further with AS $\mathbf{r_{\text{AS}}=1}$ without compromising too much the accuracy. The AS response surface is built over the GPR, RBF, NN interpolations/regressions with an additional GPR from the active one-dimensional variables to the same outputs. The training and test errors are evaluated on $\mathbf{n_{\text{ROM},\text{test}}=80}$ and $\mathbf{n_{\text{ROM},\text{test}}=20}$ training and test samples.}
\end{figure}
\begin{figure}[htpb!]
    \includegraphics[width=\textwidth]{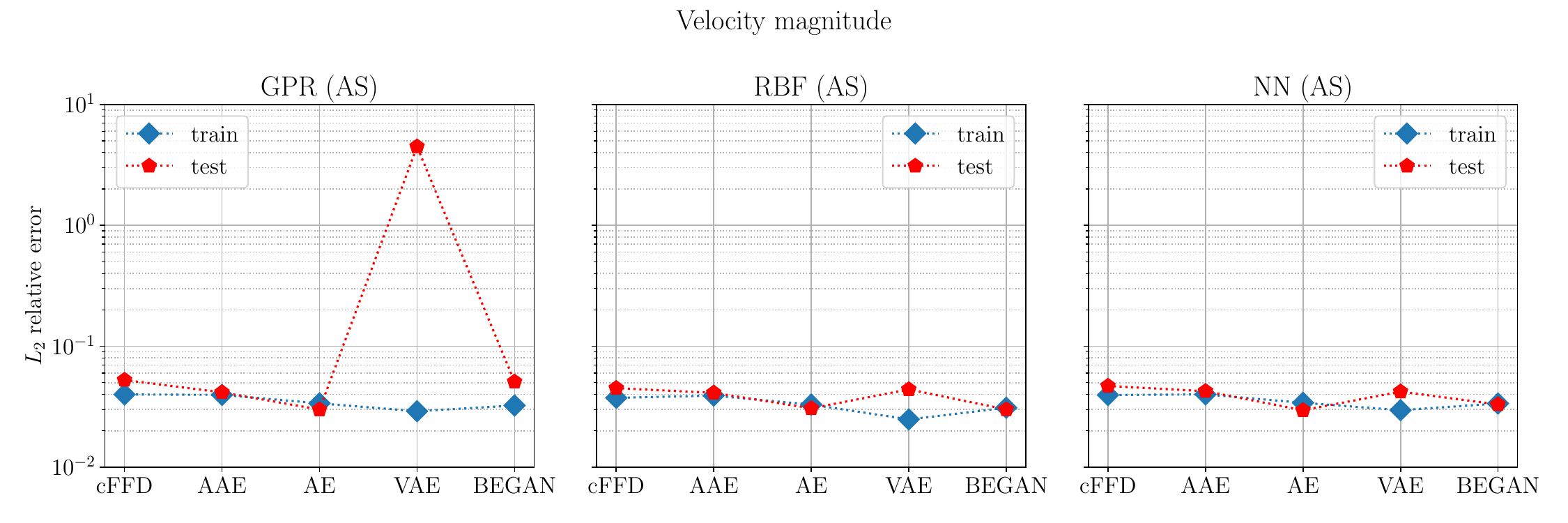}
    \caption{\textbf{HB. }In this figure we show the ROM-PODI performance, for the velocity magnitude on the whole domain, on the training and test datasets coupled with AS dimension reduction in the space of parameters. It can be seen that the latent dimension can be reduced further with AS $\mathbf{r_{\text{AS}}=1}$ without compromising too much the accuracy. The AS response surface is built over the GPR, RBF, NN interpolations/regressions with an additional GPR from the active one-dimensional variables to the same outputs. The training and test errors are evaluated on $\mathbf{n_{\text{ROM},\text{test}}=80}$ and $\mathbf{n_{\text{ROM},\text{test}}=20}$ training and test samples.}
\end{figure}

%% file: sections/discussion.tex
We address some questions that may arise, the critical parts of the developed methodology and other observations:
\begin{itemize}
  \item more complex 3d objects. The respective distributions of the test cases \textbf{SB} and \textbf{HB} are relatively easy to approximate with GMs, in fact, we need only $n_{\text{train}}=400$ training data. Provided that enough data and a sufficient computational budget are available, our methodology can be extended to more complex distributions of 3d objects. For example, points clouds with different numbers of points for each geometry or topological changes can still be approximated by GMs. Eventually, Graph Neural networks can be employed~\cite{bronstein2017geometric}.
  \item validate the results with ad hoc metrics. Especially when cGMs are employed for real-world applications, guarantees that the new geometries are not biased or satisfy certain requirements must be verified. For our test cases, we devised several geometrical and physical metrics relevant to validate our cGMs. In particular, from the drag coefficient histograms in Figure~\ref{fig:hist_cd} we could detect a bias that changed the mean value of the distributions of the cGMs. In our case, it probably depends on the application of RBF interpolation to deform the reference computational mesh, on the way the bulb was extracted from the hull and on the coarseness of the mesh. As mentioned in section~\ref{sec:mor}, RBF interpolation can be avoided and the CGMs can be directly trained on specific regions of the computational mesh as it was done for the \textbf{SB} test case.
  \item the constrained enforcing layer is independent of the numerical methods used to generate the training datasets. For our test cases, we used cFFD because it is relatively simple to implement, even if the constrained deformations it provides are not particularly large. It is crucial to observe that the cGMs and their multilinear constraints enforcing layers do not depend on how the training data are obtained. For example, as long linear or multilinear and not nonlinear geometrical constraints are enforced, also other volume-preserving techniques like~\cite{vonfunck, eisemberg} can be used.
  \item intrusive model order reduction. It is clear that cGMs with their nonlinear layers introduce nonlinearities also in the numerical models at hand through the mesh deformations. From the point of view of intrusive model order reduction~\cite{rozza2022advanced}, this adds an ulterior complexity especially when the numerical model itself is linear. In this case, hyper-reduction methods must be implemented to completely lose the dependency on the number of degrees of freedom. Non-intusive reduced order models are instead not affected.
  \item principal component analysis as preprocessing for GMs. We remark that our choice to employ PCA modes to project the training datasets of the GMs onto lower dimensional spaces is valid only when the training distribution can be well approximated with linear subspaces. If this is not the case, an evident result is a high test error even though the training error is low because the GMs cannot generalize well. Moreover, with the projection of the training data through PCA modes, the spatial correlations are lost in favor of the correlations in the space of frequencies of PCA. These lost spatial correlations could instead be exploited with ad hoc architectures like convolutional layers on Cartesian grids or graph neural networks for general meshes.
\end{itemize}

%% file: sections/conclusions.tex
A novel approach to impose linear or multilinear constraints on classical generative models is tested successfully on two academic benchmarks: the Stanford bunny and the Duisburg test case's naval hull. We proposed to validate the results through specific geometrical and physical metrics of interest for the problem at hand. How to perform non-intrusive model order reduction with proper orthogonal decomposition with interpolation was shown. The benefits of coupling dimension reduction in the space of parameters with cGMs and model order reduction are clear. 

Future directions of research include the approximation of more complex distributions with possibly different numbers of degrees of freedom per geometry or topological changes. A relevant task is to enforce nonlinear constraints without weighing too much on the computational cost of the training of the cGMs. We think that our methodology could speed up the generation of geometrically constrained 3d shapes and solve the problem of the lack of experimental data typical of some real-world applications.

%% file: sections/appendix.tex
\section{Generative models' training specifics and architectures}
\label{sec:appendix}
The following specifics are the same for the test cases \textbf{SB} of section~\ref{subsec:bunny} and \textbf{HB} of section~\ref{subsec:bulb}. Every hidden unit of our models is composed of a linear layer, a normalization layer (typically batch normalization), an Activation Layer (typically ReLU), and a dropout Layer. The AdamW~\cite{loshchilov2017decoupled} optimizer with $500$ epochs is employed with a learning rate of $1e-3$ without any scheduler. The batch size is $200$. The entire dataset size is $600$. We use $n_{\text{train}}=400$ samples for the training set and the remaining for the test set. An Nvidia RTX 3050 has been used for the training. In the following Tables~\ref{tab:aae}, \ref{tab:ae}, \ref{tab:vae} and \ref{tab:began}, \textit{Act.} stands for activation and \textit{Norm.} for normalization: they summarize in order, the architectures of the AAE, AE, VAE and BEGAN GMs we employed for both the test cases in section~\ref{sec:results}.

\begin{table}[htp!]
  \caption{Adversarial Autoencoder Structure. The number of PCA modes used is $r_{\text{PCA}}=140$ for the rabbit and $r_{\text{PCA}}=30$ for the hull, $M$ is the number of mesh points and it is equal to $145821$ for the \textbf{SB} test case and $5000$ for the \textbf{HB} test case.
}
  \centering
  \footnotesize
  \begin{tabular}{ l | c | c | c | c}
      \hline
      \hline
      Encoder & Act. & Weights & Norm. & Dropout \\
      \hline
      \hline
      None & PCA   & [ M, $r_{\text{PCA}}$] & None & None \\
      \hline
      Linear & ReLU   & [$r_{\text{PCA}}$, 500] & Batch (Affine) & 0.1 \\
      \hline
      Linear & ReLU   & [500, 500] & Batch (Affine) & 0.1 \\
      \hline
      Linear & ReLU   & [500, 500] & Batch (Affine) & 0.1 \\
      \hline
      Linear & ReLU   & [500, 500] & Batch (Affine) & 0.1 \\
      \hline
      Linear & ReLU   & [500, 500] & Batch (Affine) & 0.1 \\
      \hline
      Linear & ReLU   & [500, 500] & Batch (Affine) & 0,1 \\
      \hline
       Linear & None   & [500, 10] & Batch (NonAffine) & None \\
      \hline
      \hline
  \end{tabular}
  \hspace{1mm}
    \begin{tabular}{ l | c | c | c | c}
      \hline
      \hline
      Decoder & Act. & Weights & Norm. & Dropout \\
      \hline
      \hline
      Linear & ReLU   & [10, 500] & Batch & 0.1 \\
      \hline
      Linear & ReLU   & [500, 500] & Batch & 0.1 \\
      \hline
      Linear & ReLU   & [500, 500] & Batch & 0.1 \\
      \hline
      Linear & ReLU   & [500, 500] & Batch & 0.1 \\
      \hline
      Linear & ReLU   & [500, 500] & Batch & 0.1 \\
      \hline
      Linear & ReLU   & [500, 500] & Batch & 0.1 \\
      \hline
       Linear & None   & [500, $r_{\text{PCA}}$] & None & None \\
       \hline
       None & PCA (inverse)   & [$r_{\text{PCA}}$,  M] & None & None \\
      \hline
      \hline
  \end{tabular}\\
  \hfill\\
      \begin{tabular}{ l | c | c | c | c}
      \hline
      \hline
      Discriminator & Act. & Weights & Norm. & Dropout \\
      \hline
      \hline
      Linear & ReLU   & [10, 500] & Batch & 0.95 \\
      \hline
      Linear & ReLU   & [500, 500] & Batch & 0.95 \\
      \hline
      Linear & ReLU   & [500, 500] & Batch & 0.95 \\
      \hline
      Linear & ReLU   & [500, 500] & Batch & 0.95 \\
      \hline
      Linear & ReLU   & [500, 500] & Batch & 0.95 \\
      \hline
      Linear & ReLU   & [500, 500] & Batch & 0.95 \\
      \hline
       Linear & Sigmoid   & [500, 1] & None & None \\
      \hline
      \hline
  \end{tabular}
  \label{tab:aae}
\end{table}

\begin{table}[htp!]
  \caption{Autoencoder Structure.  The number of PCA modes used is $r_{\text{PCA}}=140$ for the rabbit and $r_{\text{PCA}}=30$ for the hull, $M$ is the number of mesh points and it is equal to $145821$ for the \textbf{SB} test case and $5000$ for the \textbf{HB} test case.
}
  \centering
  \footnotesize
  \begin{tabular}{ l | c | c | c | c}
      \hline
      \hline
      Encoder & Act. & Weights & Norm. & Dropout \\
      \hline
      \hline
      None & PCA   & [ M, $r_{\text{PCA}}$] & None & None \\
      \hline
      Linear & ReLU   & [$r_{\text{PCA}}$, 500] & Batch (Affine) & 0.1 \\
      \hline
      Linear & ReLU   & [500, 500] & Batch (Affine) & 0.1 \\
      \hline
      Linear & ReLU   & [500, 500] & Batch (Affine) & 0.1 \\
      \hline
      Linear & ReLU   & [500, 500] & Batch (Affine) & 0.1 \\
      \hline
      Linear & ReLU   & [500, 500] & Batch (Affine) & 0.1 \\
      \hline
      Linear & ReLU   & [500, 500] & Batch (Affine) & 0,1 \\
      \hline
       Linear & None   & [500, 10] & Batch (NonAffine) & None \\
      \hline
      \hline
  \end{tabular}
  \hspace{1mm}
    \begin{tabular}{ l | c | c | c | c}
      \hline
      \hline
      Decoder & Act. & Weights & Norm. & Dropout \\
      \hline
      \hline
      Linear & ReLU   & [10, 500] & Batch & 0.1 \\
      \hline
      Linear & ReLU   & [500, 500] & Batch & 0.1 \\
      \hline
      Linear & ReLU   & [500, 500] & Batch & 0.1 \\
      \hline
      Linear & ReLU   & [500, 500] & Batch & 0.1 \\
      \hline
      Linear & ReLU   & [500, 500] & Batch & 0.1 \\
      \hline
      Linear & ReLU   & [500, 500] & Batch & 0.1 \\
      \hline
       Linear & None   & [500, $r_{\text{PCA}}$] & None & None \\
       \hline
       None & PCA (inverse)   & [$r_{\text{PCA}}$,  M] & None & None \\
      \hline
      \hline
  \end{tabular}\\
  \hfill\\
  \label{tab:ae}
\end{table}

\begin{table}[htp!]
  \caption{Beta Variational Autoencoder Structure.  The number of PCA modes used is $r_{\text{PCA}}=140$ for the rabbit and $r_{\text{PCA}}=30$ for the hull, $M$ is the number of mesh points and it is equal to $145821$ for the \textbf{SB} test case and $5000$ for the \textbf{HB} test case.
}
  \centering
  \footnotesize
  \begin{tabular}{ l | c | c | c | c}
      \hline
      \hline
      Encoder & Act. & Weights & Norm. & Dropout \\
      \hline
      \hline
      None & PCA   & [ M, $r_{\text{PCA}}$] & None & None \\
      \hline
      Linear & ReLU   & [$r_{\text{PCA}}$, 500] & Batch (Affine) & 0.1 \\
      \hline
      Linear & ReLU   & [500, 500] & Batch (Affine) & 0.1 \\
      \hline
      Linear & ReLU   & [500, 500] & Batch (Affine) & 0.1 \\
      \hline
      Linear & ReLU   & [500, 500] & Batch (Affine) & 0.1 \\
      \hline
      Linear & ReLU   & [500, 500] & Batch (Affine) & 0.1 \\
      \hline
      Linear & ReLU   & [500, 500] & Batch (Affine) & 0,1 \\
      \hline
       Linear & None   & [500, 10] & Batch (NonAffine) & None \\
      \hline
      \hline
  \end{tabular}
  \hspace{1mm}
    \begin{tabular}{ l | c | c | c | c}
      \hline
      \hline
      Decoder & Act. & Weights & Norm. & Dropout \\
      \hline
      \hline
      Linear & ReLU   & [10, 500] & Batch & 0.1 \\
      \hline
      Linear & ReLU   & [500, 500] & Batch & 0.1 \\
      \hline
      Linear & ReLU   & [500, 500] & Batch & 0.1 \\
      \hline
      Linear & ReLU   & [500, 500] & Batch & 0.1 \\
      \hline
      Linear & ReLU   & [500, 500] & Batch & 0.1 \\
      \hline
      Linear & ReLU   & [500, 500] & Batch & 0.1 \\
      \hline
       Linear & None   & [500, $r_{\text{PCA}}$] & None & None \\
       \hline
       None & PCA (inverse)   & [$r_{\text{PCA}}$,  M] & None & None \\
      \hline
      \hline
  \end{tabular}\\
  \hfill\\
  \label{tab:vae}
\end{table}
\begin{table}[htp!]
  \caption{Boundary Equilibrium Generative Adversarial Network Structure.  The number of PCA modes used is $r_{\text{PCA}}=140$ for the rabbit and $r_{\text{PCA}}=30$ for the hull, $M$ is the number of mesh points and it is equal to $145821$ for the \textbf{SB} test case and $5000$ for the \textbf{HB} test case. The terms \textit{DisEnc} and \textit{DisDec} stand for the discriminator's encoder and decoder.
}
  \centering
  \footnotesize
  \begin{tabular}{ l | c | c | c | c}
      \hline
      \hline
      DisEnc & Act. & Weights & Norm. & Dropout \\
      \hline
      \hline
      None & PCA   & [ M, $r_{\text{PCA}}$] & None & None \\
      \hline
      Linear & ReLU   & [$r_{\text{PCA}}$, 500] & Batch (Affine) & 0.1 \\
      \hline
      Linear & ReLU   & [500, 500] & Batch (Affine) & 0.1 \\
      \hline
      Linear & ReLU   & [500, 500] & Batch (Affine) & 0.1 \\
      \hline
      Linear & ReLU   & [500, 500] & Batch (Affine) & 0.1 \\
      \hline
      Linear & ReLU   & [500, 500] & Batch (Affine) & 0.1 \\
      \hline
      Linear & ReLU   & [500, 500] & Batch (Affine) & 0,1 \\
      \hline
       Linear & None   & [500, 10] & Batch (NonAffine) & None \\
      \hline
      \hline
  \end{tabular}
  \hspace{1mm}
  \begin{tabular}{ l | c | c | c | c}
    \hline
    \hline
    DisDec & Act. & Weights & Norm. & Dropout \\
    \hline
    \hline
    Linear & ReLU   & [10, 500] & Batch & 0.1 \\
    \hline
    Linear & ReLU   & [500, 500] & Batch & 0.1 \\
    \hline
    Linear & ReLU   & [500, 500] & Batch & 0.1 \\
    \hline
    Linear & ReLU   & [500, 500] & Batch & 0.1 \\
    \hline
    Linear & ReLU   & [500, 500] & Batch & 0.1 \\
    \hline
    Linear & ReLU   & [500, 500] & Batch & 0.1 \\
    \hline
     Linear & None   & [500, $r_{\text{PCA}}$] & None & None \\
     \hline
     None & PCA (inverse)   & [$r_{\text{PCA}}$,  M] & None & None \\
    \hline
    \hline
  \end{tabular}\\
  \hfill\\

  \begin{tabular}{ l | c | c | c | c}
    \hline
    \hline
    Generator & Act. & Weights & Norm. & Dropout \\
    \hline
    \hline
    Linear & ReLU   & [10, 500] & Batch & 0.1 \\
    \hline
    Linear & ReLU   & [500, 500] & Batch & 0.1 \\
    \hline
    Linear & ReLU   & [500, 500] & Batch & 0.1 \\
    \hline
    Linear & ReLU   & [500, 500] & Batch & 0.1 \\
    \hline
    Linear & ReLU   & [500, 500] & Batch & 0.1 \\
    \hline
    Linear & ReLU   & [500, 500] & Batch & 0.1 \\
    \hline
      Linear & None   & [500, $r_{\text{PCA}}$] & None & None \\
      \hline
      None & PCA (inverse)   & [$r_{\text{PCA}}$,  M] & None & None \\
    \hline
    \hline
  \end{tabular}
  \hfill\\
  \label{tab:began}
\end{table}

\section{PODI-NN specifics}
The training set consists of $\mathbf{n_{\text{ROM}, \text{train}}=80}$ samples and the test set of $\mathbf{n_{\text{ROM}, \text{test}}=20}$ samples. The AdamW~\cite{loshchilov2017decoupled} optimizer with $1000$ epochs has been used, with a learning rate of $1e-3$ without any scheduler. 

\begin{table}[htp!]
  \caption{PODI-NN interpolation: $r_{\text{POD}}=3$ for both test cases \textbf{SB} and \textbf{HB}, $M$ is the number of degrees of freedom and $P$ the number of parameters that can vary if cFFD ($P=p=54$ for \textbf{SB} or $P=p=84$ for \textbf{HB}), or cGMs ($P=R=15$ for \textbf{SB} or $P=R=10$ for \textbf{HB}) are considered.
}
  \centering
  \footnotesize
  \begin{tabular}{ l | c | c | c | c}
      \hline
      \hline
      Regressor & Activation & Weights & Norm. & Dropout \\
      \hline
      \hline
      Linear & ReLU   & [P, 2000] & None & 0 \\
      \hline
      Linear & ReLU   & [2000, 2000] & None & 0 \\
      \hline
      Linear & ReLU   & [2000, $r_{\text{POD}}$] & None & 0 \\
      \hline
      POD reconstruction & -   & [$r_{\text{POD}}$, M] & None & 0 \\
      \hline
      \hline
  \end{tabular}

\end{table}